\documentclass[12pt]{amsart}
\usepackage{amsfonts}
\usepackage{amsthm}
\usepackage{amsmath,amsthm,amssymb}
\usepackage{amscd}
\usepackage[latin2]{inputenc}
\usepackage{t1enc}
\usepackage[mathscr]{eucal}
\usepackage{indentfirst}
\usepackage{graphicx}
\usepackage{graphics}
\usepackage{pict2e}
\usepackage{epic}
\usepackage{cleveref}
\usepackage{color, colortbl} 
\usepackage{diagbox}

\graphicspath{{figs/}}

\numberwithin{equation}{section}
\usepackage[margin=2.9cm]{geometry}
\usepackage{epstopdf} 
\usepackage{MnSymbol}
\usepackage{multicol}
\usepackage{mathtools}
\usepackage{url}
\usepackage{xcolor}
\usepackage[font=small,labelfont=bf]{caption}
\usepackage[linesnumbered,ruled,algosection]{algorithm2e}

\usepackage{mathtools}
\usepackage{blkarray, bigstrut} %

\usepackage{todonotes}

 \usepackage{subcaption}
\usepackage{tikz}
\usepackage{tikz-cd}
\usetikzlibrary{arrows.meta}
\usetikzlibrary{bending}
 
 \usepackage{MnSymbol}

\theoremstyle{plain}
\newtheorem{Th}{Theorem}[section]
\newtheorem{Lemma}[Th]{Lemma}
\newtheorem{Cor}[Th]{Corollary}
\newtheorem{Prop}[Th]{Proposition}

\newtheorem{ThmDef}[Th]{Theorem/Definition}

\theoremstyle{definition}
\newtheorem{Def}[Th]{Definition}

\newtheorem{Rem}[Th]{Remark}
\newtheorem{?}[Th]{Problem}
\newtheorem{Ex}[Th]{Example}

\newcommand{\antiexch}{\operatorname{AE}}

\newcommand{\co}{\operatorname{co}}

\newcommand{\C}{\mathbb{C}}

\newcommand{\M}{\mathcal{M}}

\newcommand{\QQ}{\mathcal{Q}}

\newcommand{\SSS}{\mathcal{S}}
\newcommand{\uv}{\bigl[ \begin{smallmatrix}v \\ u\end{smallmatrix}\bigr]}

\newcommand{\uprimevprime}{\bigl[ \begin{smallmatrix}v' \\ u'\end{smallmatrix}\bigr]}

\newcommand{\doubleloop}{\begin{tikzpicture}[scale=1.5,>=latex]
  \def\Radius{.1cm}

  \draw (0cm,0cm) circle[radius=\Radius];

  \begin{scope}[
    -{Stealth[length=3pt, width=3pt]},
    shorten >=1pt,
    very thin,
  ]
    \draw (\Radius, 0) arc(-3:3:\Radius);
    \draw (-\Radius, 0) arc(180+3:180-3:\Radius);
  \end{scope}
  
\end{tikzpicture}}

\newcommand\tinydoubleloop{\vcenter{\hbox{\scalebox{0.5}{\doubleloop}}}}

\newcommand\cw{\vcenter{\hbox{\scalebox{0.8}{$\circlearrowright$}}}}
\newcommand\ccw{\vcenter{\hbox{\scalebox{0.8}{$\circlearrowleft$}}}}

\newcommand{\w}{w^{\hspace{0.01in} \tinydoubleloop}}

\newcommand{\z}{z^{\hspace{0.01in} \tinydoubleloop}}
\newcommand{\SSn}{S_n^{\hspace{0.01in} \tinydoubleloop}}
\newcommand{\SSnk}{S_{n,k}^{\hspace{0.01in} \tinydoubleloop}}
\newcommand{\SSnminusk}{S_{n,n-k}^{\hspace{0.01in} \tinydoubleloop}}

\newcommand{\Gkn}{Gr(k,n)}
\newcommand{\GIkn}{Gi(k,n)}

\newcommand{\aligns}{Alignments}
\newcommand{\necklace}{\mathcal{N}}
\newcommand{\clockwise}[1]{\overrightarrow{#1}}
\newcommand{\counterclockwise}[1]{\overleftarrow{#1}}

\definecolor{LightCyan}{rgb}{0.88,1,1}
\definecolor{yellow}{rgb}{1,1,0}

\DeclareMathOperator{\GL}{GL}

\DeclareMathOperator{\rank}{rank}

\usepackage{lipsum}

\begin{document}

\title[Criteria for smoothness of Positroid varieties]{Criteria for smoothness of Positroid varieties \\ via pattern avoidance, Johnson graphs, and spirographs}

\author[S. Billey, J. Weaver]{Sara C. Billey and Jordan E. Weaver}

\date{\today}

\thanks{Both authors
were partially supported by the National Science Foundation grant
DMS-1764012. Email:billey@math.washington.edu and jeweaver@uw.edu.}

\address{Department of Mathematics, University of Washington, Seattle, WA, USA}

\begin{abstract}{ Positroids are certain representable matroids
originally studied by Postnikov in connection with the totally
nonnegative Grassmannian and now used widely in algebraic
combinatorics.  The positroids give rise to determinantal equations
defining positroid varieties as subvarieties of the Grassmannian
variety. Rietsch, Knutson-Lam-Speyer, and Pawlowski studied geometric
and cohomological properties of these varieties.  In this paper, we
continue the study of the geometric properties of positroid varieties
by establishing several equivalent conditions characterizing smooth
positroid varieties using a variation of pattern avoidance defined on
decorated permutations, which are in bijection with positroids.  This
allows us to give two formulas for counting the number of smooth
positroids along with two $q$-analogs.  Furthermore, we give a
combinatorial method for determining the dimension of the tangent
space of a positroid variety at key points using an induced subgraph
of the Johnson graph.  We also give a Bruhat interval characterization
of positroids.}
\end{abstract}

\maketitle

\tableofcontents

\section{Introduction}
\label{sec: intro}

\textit{Positroids} are an important family of realizable matroids
originally defined by Postnikov in \cite{Postnikov.2006} in the
context of the totally nonnegative part of the Grassmannian
variety. These matroids and the totally positive part of the
Grassmannian variety have played a critical role in the theory of
cluster algebras and soliton solutions to the KP equations and have
connections to statistical physics, integrable systems, and scattering
amplitudes
\cite{AHBCGPT.2016,Brown-Goodearl-Yakimov,Fomin-Williams-Zelevinsky,Lusztig.1998,Rietsch.2006,Williams.2007,williams2021positive}.
Positroids are closed under restriction, contraction, duality, and
cyclic shift of the ground set, and furthermore they have particularly
elegant matroid polytopes \cite{Ardila-Rincon-Williams}.

\textit{Positroid varieties} were studied by Knutson, Lam, and Speyer
in \cite{KLS}, building on the work of Lusztig, Postnikov and Rietsch.
They are homogeneous subvarieties of the complex Grassmannian variety
$\Gkn$ which are defined by determinantal equations determined by the
nonbases of a positroid.  They can also be described as projections of
Richardson varieties in the complete flag manifold
to $\Gkn$.  These varieties have beautiful geometric, representation
theory, and combinatorial connections \cite{KLS2,Paw}.  See the
background section for notation and further background.

Positroids and positroid varieties can be bijectively associated with
many different combinatorial objects \cite{Oh,Postnikov.2006}. For the
purposes of this paper, we will need to use the bijections associating
each of the following types of objects with each other:

\begin{enumerate}
\item positroids $\M$ of rank $k$ on a ground set of size $n$, 
\item decorated permutations $\w$ on $n$ elements with $k$ anti-exceedances, 
\item Grassmann necklaces $\necklace=(I_{1},\dotsc , I_{n}) \in 
\binom{[n]}{k}^{n}$,  and 
\item Grassmann intervals $[u,v]$ in $\GIkn$.  
\end{enumerate}
Here, a decorated permutation $\w$ on $n$ elements is a permutation $w
\in S_n$ together with an orientation clockwise or counterclockwise,
denoted $\clockwise{i}$ or $\counterclockwise{i}$ respectively, assigned to each
fixed point of $w$.  A \textit{Grassmann interval} $[u,v] \in \GIkn$
is an interval in Bruhat order on permutations in $S_{n}$ such that
$v$ has at most one descent, specifically in position $k$.  In
\Cref{sec:background}, we will sketch the relevant bijections and
remaining terminology. In addition to these, there are bijections to
juggling sequences, \reflectbox{L}-diagrams, equivalence classes of
plabic graphs, and bounded affine permutations
\cite{Ardila-Rincon-Williams,KLS,Postnikov.2006}.

Many of the properties of positroid varieties can be ``read off'' from
one or more of these bijectively equivalent definitions.  Thus, we
will index a positroid variety
\begin{equation}\label{eq:equivalent.defs}
\Pi_{\M}=\Pi_{\w}=\Pi_{\necklace}=\Pi_{[u,v]}
\end{equation}
using any of the associated objects, depending on the relevant context.
For example, the codimension is easy to read off from the decorated
permutation using the notions of the chord diagram and its
alignments.

Let $\SSnk$ be the set of decorated permutations on $n$ elements with
$k$ anti-exceedances.  The \textit{chord diagram} $D(\w)$ of $\w \in
\SSnk$ is constructed by placing the numbers $1, 2,\ldots, n$ on $n$
vertices around a circle in clockwise order, and then, for each $i$,
drawing a directed arc from $i$ to $w(i)$ with a minimal number of
crossings between distinct arcs while staying completely inside the
circle.  The arcs beginning at fixed points should be drawn clockwise
or counterclockwise according to their orientation in $\w$.

An \textit{alignment} in $D(\w)$ is a pair of directed edges $(i
\mapsto w(i), j \mapsto w(j)) $ which can be drawn as distinct
noncrossing arcs oriented in the same direction.  A pair of directed
edges $(i \mapsto w(i), j \mapsto w(j)) $ which can be drawn as
distinct noncrossing arcs oriented in opposite directions is called a
\textit{misalignment}.  A pair of directed edges which must cross if
both are drawn inside the cycle is called a \textit{crossing}
\cite[Sect. 5]{Postnikov.2006}.  Let $\aligns(\w)$ denote the set of
alignments of $\w$.

\begin{Ex} \label{Example:intro} Let $\w =895\counterclockwise{4}7\clockwise{6}132$
be the decorated permutation with a
counterclockwise fixed point at $4$ and a clockwise fixed point at
$6$.  The chord diagram for $\w$ is the following.
\begin{center}
\includegraphics[height=4cm]{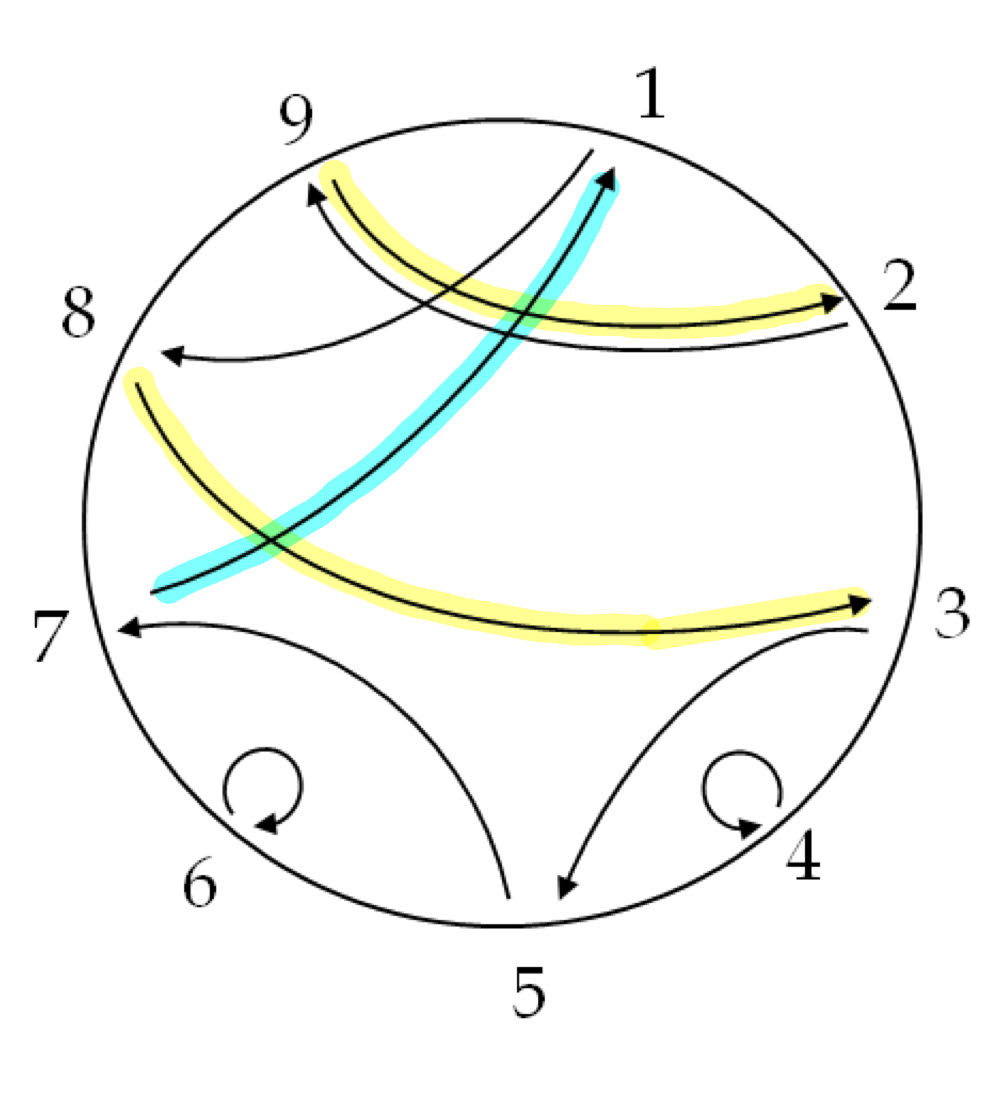}
\end{center}
Here for example, $(9\mapsto 2, 8 \mapsto 3)$ highlighted in yellow is an alignment,
$(2\mapsto 9, 8 \mapsto 3)$ is a misalignment, and both $(7\mapsto 1,
8 \mapsto 3)$ and $(7\mapsto 1, 5 \mapsto 7)$ are crossings.  Note,
$(7\mapsto 1, 6 \mapsto 6)$ is an alignment and $(7\mapsto 1, 4
\mapsto 4 )$ is a misalignment.  
\end{Ex}

\begin{Th}\label{thm:Alignments} \cite{KLS,Postnikov.2006} For any
decorated permutation $\w \in \SSnk$ and associated Grassmann interval
$[u,v]$, the codimension of $\Pi_{\w}$ in $\Gkn$ is
\begin{equation}\label{eq:codimPI}
\mathrm{codim}(\Pi_{\w}) =  \#\aligns(\w) = k(n-k) - [\ell(v) - \ell(u)].
\end{equation}
\end{Th}

We use the explicit equations defining a positroid
variety in $\Gkn$ to determine whether the variety is smooth or singular.
In general, a variety $X$ defined by polynomials $f_1, \ldots, f_s$ is
\textit{singular} if there exists a point $x \in X$ such that the
Jacobian matrix, $Jac$, of partial derivatives of the $f_i$ satisfies
$\text{rank}(Jac|_x) < \text{codim}\ X$. It is \textit{smooth} if no
such point exists.  The value $\text{rank}(Jac|_x)$ is the codimension
of the tangent space to $X$ at the point $x$.  Thus,
$\text{rank}(Jac|_x) < \text{codim}\ X$ implies the dimension of the
tangent space to the variety $X$ at $x$ is strictly larger than the
dimension of the variety $X$, hence $x$ is a singularity like a cusp on a
curve.  In the case of a positroid variety $\Pi_{\w}$,
\Cref{thm:Alignments} implies that a point $x \in \Pi_{\w}$ is a
singularity of $\Pi_{\w}$ if
\begin{equation} \label{eq: sing condition Jac algn}
\text{rank}(Jac|_x) \; < \; \text{codim} \; \Pi_{\w} \; = \; \#\aligns(\w).
\end{equation}

Our first main theorem reduces the problem of finding singular points
in a positroid variety to checking the rank of the Jacobian only at a
finite number of $T$-fixed points.  For any $J = \{j_1, \ldots, j_k\} \subseteq [n]$, let $A_{J}$ be the element in
$\Gkn$ spanned by the elementary row vectors $e_{i}$ with $i \in J$,
or equivalently the subspace represented by a $k\times n$ matrix with
a $1$ in cell $(i,j_i)$ for each $j_i \in J$ and zeros everywhere
else.  These are the $T$-fixed points of $\Gkn$, where $T$ is the set
of invertible diagonal matrices over $\mathbb{C}$. The reduction follows from the
decomposition of $\Pi_{[u,v]}$ as a projected Richardson variety. Every point $A \in \Pi_{[u,v]}$ lies in the projection of some
intersection of a Schubert cell with an opposite Schubert variety
$C_{y}\cap X^{v}$ for $y\in [u,v]$. In particular,
if  $y=y_{1} y_{2}\cdots y_{n} \in [u,v]$ in one-line notation and we define
$y[k] :=\{y_{1},y_{2},\dots, y_{k}\}$ to be an \textit{initial set} of $y$, then $A_{y[k]}$ is in the
projection of $C_y \cap X^v$.

\begin{Th}\label{thm:bounded.below}
Assume $A \in \Pi_{[u,v]}$ is the image of a point in $C_{y} \cap X^v$
projected to $\Gkn$ for some $y \in [u,v]$. Then the codimension of the tangent space to
$\Pi_{[u,v]}$ at $A$ is bounded below by $\rank(Jac|_{A_{y[k]}})$.
\end{Th}

\Cref{thm:bounded.below} indicates that the $T$-fixed points of the
form $A_{y[k]}$ such that $y \in [u,v]$ are key to understanding the
singularities of $\Pi_{[u,v]}$.  In fact, the equations determining
$\Pi_{[u,v]}$ and the bases of the positroid $\M$ associated to the
interval $[u,v]$ can be determined from the permutations in the
interval by the following theorem.  Our proof of the following theorem
depends on Knutson-Lam-Speyer's \Cref{Def: positroid as projection} of
a positroid variety as the projection of a Richardson variety and
\Cref{Lemma: BS}. It also follows from \cite[Lemma 3.11]{KW2015}.
Both groups, Knutson-Lam-Speyer and Kodama-Williams, were aware of
this result in the context of positroid varieties for Grassmannians,
but it does not appear to be in the literature in the form we needed,
hence we prove the result in \Cref{sec:initial.sets}.

\begin{Th}  \label{Th: positroid = initial sets} \normalfont Let
$\w \in \SSnk$ with associated Grassmann interval $[u,v]$ and
positroid $\M$. Then $\M$ is exactly the collection of initial sets of permutations in the Grassmann interval $[u,v]$,
\[
\M = \{y[k] \; : \; y \in [u,v]\}.
\]
\end{Th}

Our next theorem provides a method to compute the rank of the Jacobian
of $\Pi_{{[u,v]}}$ explicitly at the $T$-fixed points.  Therefore, we
can also compute the dimension of the tangent space of a positroid
variety at those points.  Comparing that with the number of alignments
gives a test for singularity of points in positroid varieties by \eqref{eq:
sing condition Jac algn}.  

\begin{Th} \label{thm:tangent.space.dim} Let
$\w \in \SSnk$ with associated Grassmann interval $[u,v]$ and
positroid $\M$. For $y\in [u,v]$, the codimension of the tangent space to $\Pi_{[u, v]}
\subseteq \Gkn$ at
$A_{y[k]}$ is
\begin{equation}\label{eq:tangent.space.codim}
\rank (Jac|_{A_{y[k]}}) \; = \; \# \Big\{I \in \binom{[n]}{k} \setminus \M \; : \; |I \cap y[k]| = k-1
\Big\}.
\end{equation}
\end{Th}

The formula in \eqref{eq:tangent.space.codim} is reminiscent of the
\textit{Johnson graph} $J(k,n)$ with vertices given by the $k$-subsets
of $[n]$ such that two $k$-subsets $I,J$ are connected by an edge precisely
if $|I\cap J|=k-1$.  For a positroid $\M \subseteq \binom{[n]}{k}$, let $J(\M)$ denote the induced subgraph of the
Johnson graph on the vertices in $\M$.  We call $J(\M)$ the
\textit{Johnson graph of} $\M$.  Note, the Johnson graph is closely
related to bases of matroids by the Basis Exchange Property.
\Cref{thm:tangent.space.dim} implies $J(\M)$ encodes aspects of the
geometry of the positroid varieties like the Bruhat graph in the
theory of Schubert varieties \cite{carrell94}. 

To state our main theorem characterizing smoothness of positroid
varieties, we need to define two types of patterns that may occur in a
chord diagram.  First, given an alignment $(i \mapsto w(i), j \mapsto
w(j)) $ in $D(\w)$, if there exists a third arc $(h \mapsto w(h))$
which forms a crossing with both $(i \mapsto w(i))$ and $(j \mapsto
w(j))$, we say $(i \mapsto w(i), j \mapsto w(j)) $ is a
\textit{crossed alignment} of $\w$. In the example above, $(9\mapsto
2, 8 \mapsto 3)$ is a crossed alignment; this alignment is crossed for
instance by $(7 \mapsto 1)$, highlighted in blue.

\begin{figure}[h]
\begin{center}
   \includegraphics[scale=0.23]{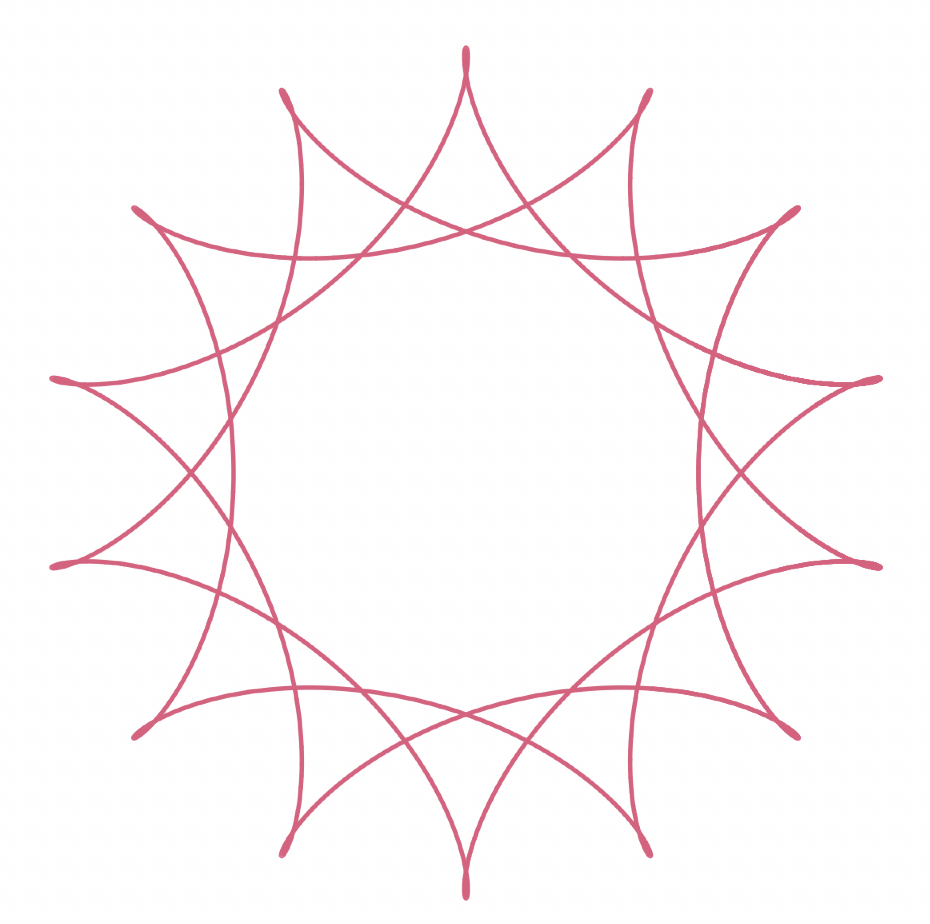}
\hspace{.2in}
   \includegraphics[scale=0.24]{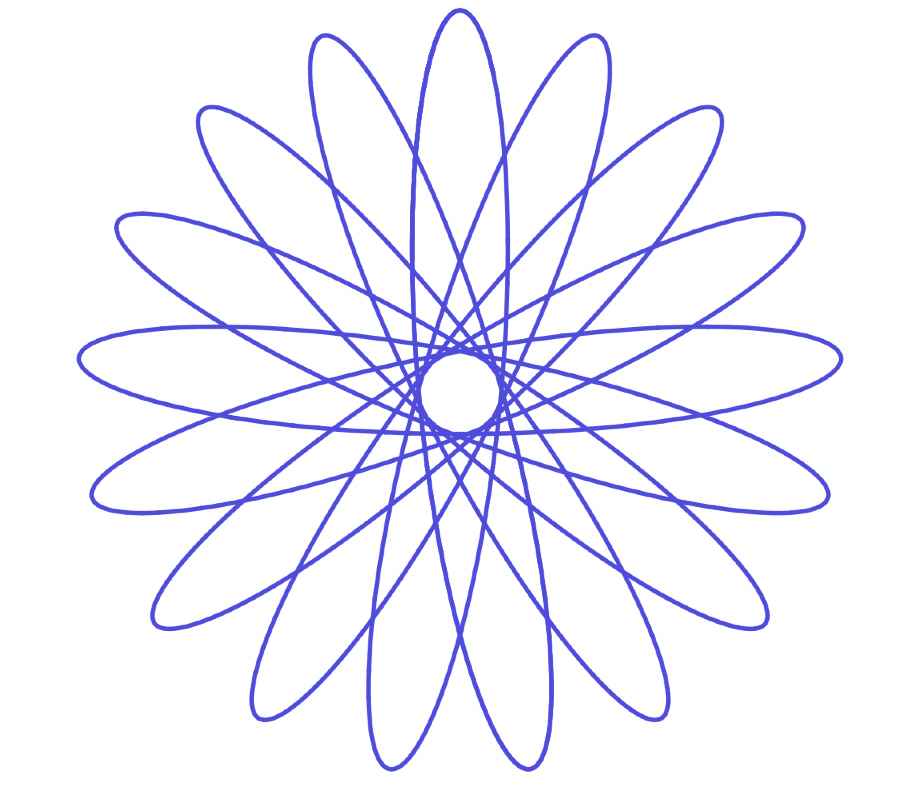}
\hspace{.2in}
   \includegraphics[scale=0.18]{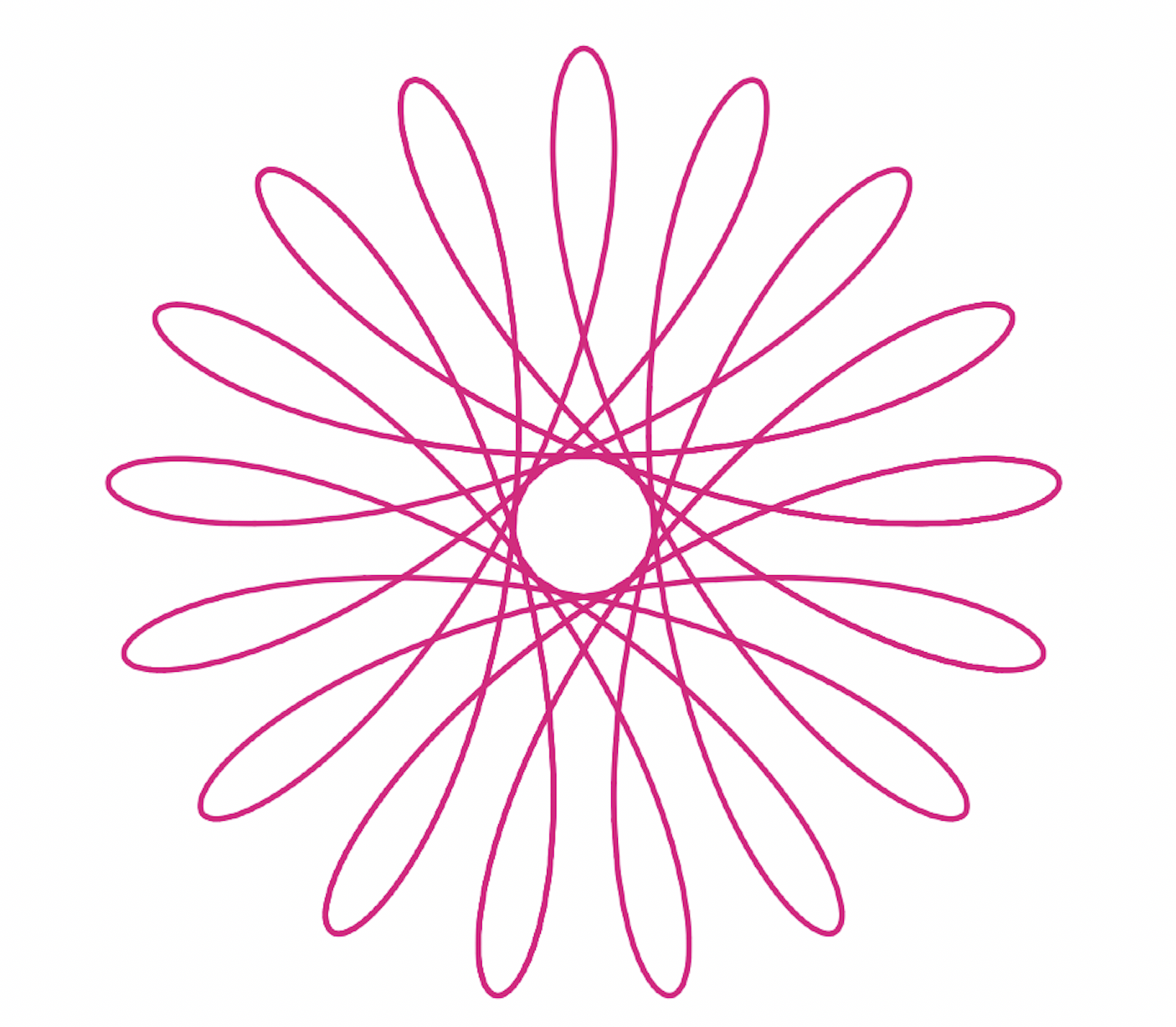}
\end{center}
\caption{Spirographs made by the Spirograph Maker app for the iphone.}  \label{Fig:spirographs}
\end{figure}

The second type of pattern is related to the
Spirograph\textsuperscript{TM} toy, designed by Denys Fisher and
trademarked by Hasbro to draw a variety of curves inside a circle
which meet the circle in a finite number of discrete points. See
\Cref{Fig:spirographs}. Once oriented and vertices are added, such
curves each determine a chord diagram from a special class, which we
will call spirographs.  We think of alignments, crossings, crossed
alignments, and spirographs as \textit{subgraph patterns} for
decorated permutations.

\begin{Def}\label{def:spirograph} A decorated permutation
$\w=(w,\co)\in \SSnk$ is a \textit{spirograph permutation} if there
exists a positive integer $m$ such that $w(i) = i+m$ (mod $n$) for all
$i$ and $\w$ has at most one fixed point.  The chord diagram of a
spirograph permutation will be called a \textit{spirograph}.
\end{Def}

%

\begin{Th} \label{The: Main theorem} Let
$\w \in \SSnk$ with associated Grassmann interval $[u,v]$ and
positroid $\M = \{y[k] \; : \; y \in [u,v]\}$. Then, the following are
equivalent.
\begin{enumerate}
\item \label{MainPart:1}The positroid variety $\Pi_{\w} = \Pi_{[u,v]} = \Pi_\M$ is smooth.
\item \label{MainPart:2} For every $y \in [u,v]$, $\# \big\{I \in \M \; : \; |I \cap y[k]| = k-1 \big\} \; = \; \ell(v)-\ell(u).$
\item \label{MainPart:3} For every $J \in \M$,  $\# \{I \in \M \; : \; |I \cap J| = k-1 \} \; = \; k(n-k) - \# \aligns(\w).$
\item The graph $J(\M)$ is regular, and each vertex
has degree  $\ell(v)-\ell(u).$
\item The decorated permutation $\w$ has no crossed alignments.
\item The chord diagram $D(\w)$ is a disjoint union of spirographs. 
\item The positroid $\M$ is a direct sum of uniform matroids. 
\end{enumerate}
\end{Th}

Schubert varieties in the flag variety are indexed by permutations and
are closely related to positroid varieties, as explained in
\Cref{sec:background}.  Smoothness of Schubert varieties is completely
characterized by pattern avoidance of the corresponding permutations
\cite{Lak-San}. When studying the partially asymmetric exclusion
process and its surprising connection to the Grassmannian, Sylvie
Corteel posed the idea of considering patterns in decorated
permutations \cite{Corteel}.  This suggestion was the foundation for
the present work.  Some of the results in this paper were announced in
\cite{Billey-Weaver-FPSAC2022}.

In \Cref{sec:background}, we provide background material from the
literature and define our notation. In \Cref{sec:initial.sets}, we
study initial sets for Grassmann intervals and prove \Cref{Th:
positroid = initial sets} as well as the equalities (2) $\Leftrightarrow$ (3) $\Leftrightarrow$ (4) of \Cref{The: Main theorem}. In \Cref{Sec: perm matrices suffice}, we
prove \Cref{thm:bounded.below} and \Cref{thm:tangent.space.dim}.
\Cref{thm:bounded.below} and \Cref{thm:tangent.space.dim} together
yield the equality (1) $\Leftrightarrow$ (2) in \Cref{The: Main theorem}.
In \Cref{Sec: rigid transformations}, we consider several rigid transformations of chord diagrams and the associated transformations of the related objects. Furthermore, we show in \Cref{Lemma: invariances} that applying these transformations to the objects indexing positroid varieties preserves the property of being smooth.
In \Cref{sec:connected}, we utilize a decomposition of $\M$ into connected components on a non-crossing partition and show that a positroid variety is smooth if and only if all of the components of the decomposition correspond to smooth positroid varieties.
In \Cref{Sec: spirograph direction}, we identify the special properties of positroid
varieties indexed by spirograph permutations, which leads to the implication (6) $\Rightarrow$ (1) and the equalities (5) $\Leftrightarrow$ (6) $\Leftrightarrow$ (7) in \Cref{The: Main theorem}.  
In \Cref{Sec: anti-exchange pairs and crossed aligns}, we complete the proof of \Cref{The: Main theorem} by showing that (1) $\Rightarrow$ (5). This implication is accomplished by constructing an injective map from the \textit{anti-exchange pairs for} a particular set $J \in \M$ to $\aligns(\w)$. 
In \Cref{sec:enumeration}, we provide some enumerative results for smooth positroid varieties. We conclude with some open questions in \Cref{sec:future}.

\section{Background}\label{sec:background}

We begin by giving notation and some background on several
combinatorial objects and theorems from the literature. These objects
will be used to index the varieties discussed throughout the paper. We
will then introduce notation for several geometrical objects,
including Grassmannian varieties, flag varieties, Schubert
varieties, Richardson varieties, and positroid varieties.

\subsection{Combinatorial objects} 

For integers $i \leq j$, let $[i,j]$ denote the set
$\{i,i+1,\ldots,j\}$, and write $[n] := [1,n]$ for a positive integer
$n$. Let $\binom{[n]}{k}$ be the set of size $k$ subsets of $[n]$ for
$k \in [0,n]$. Call $J \in \binom{[n]}{k}$ a $k$-subset of $[n]$.

Define the \textit{Gale partial order}, $\preceq$, on $k$-subsets of
$[n]$ as follows. Let $I = \{i_1 < \cdots < i_k\}$ and $J = \{j_1 <
\cdots < j_k\}$. Then $I \preceq J$ if and only if $i_h \leq j_h$ for
all $h \in [k]$. This partial order is known by many other names; we
are following \cite{Ardila-Rincon-Williams} for consistency. Gale studied this partial order in the context of matroids in the 1960s \cite{Gale.1968}.

For any $k \times n$
matrix $A$ and any set $J \in \binom{[n]}{k}$, define $\Delta_J(A)$ to be the determinant of the $k \times k$ submatrix of $A$ lying in
column set $J$. The minors $\Delta_I(A)$ for $I \in \binom{[n]}{k}$ are called the \textit{Pl{\"u}cker coordinate} of $A$. We think of $\Delta_J$ as a polynomial function on
the set of all $k \times n$ matrices over a chosen field using
variables of the form $x_{ij}$ indexed by row $i \in [k]$ and column
$j\in J$.




\begin{Def}\label{def:noncrossing.partition}
Let $S$ be a partition $[n] = B_{1} \sqcup \cdots \sqcup B_{t}$ of $[n]$ into
pairwise disjoint, nonempty subsets. We say that $S$ is a
\textit{non-crossing partition} if there are no distinct $a,b,c,d$ in
cyclic order such that $a,c \in B_{i}$ and $b,d \in B_{j}$ for some $i
\neq j$.  Equivalently, place the numbers $1, 2,\ldots, n$ on $n$
vertices around a circle in clockwise order, and then for each
$B_{i}$, draw a polygon on the corresponding vertices. If no two of
these polygons intersect, then $S$ is a non-crossing partition of
$[n]$.  
\end{Def}


A \textit{matroid} of rank $k$ on $[n]$, defined by its bases, is a
nonempty subset $\mathcal{M} \subseteq \binom{[n]}{k}$ satisfying the
following Basis Exchange Property: if $I, J \in \mathcal{M}$ such that
$I \neq J$ and $a \in I \setminus J$, then there exists some $b \in J
\setminus I$ such that $(I \setminus \{a\}) \cup \{b\} \in
\mathcal{M}$. Compare the notion of matroid basis exchange to basis
exchange in linear algebra. We call the sets in $\binom{[n]}{k} \setminus \M$ the \textit{nonbases} of $\M$, and we denote this collection of sets by $\QQ(\M)$. The set $\binom{[n]}{k}$ is the \textit{uniform matroid} of rank $k$ on $[n]$. For more on matroids, see
\cite{ArdilaChapter,Oxley}.

\begin{Ex}\label{ex:matroid}
A notable family of matroids called \textit{representable matroids}
comes from matrices.  Let $A$ be a full rank $k \times n$ matrix. The
\textit{matroid of $A$} is the set
\[
\mathcal{M}_A := \Big\{J \in \binom{[n]}{k} \; : \; \Delta_J(A) \neq 0\Big\}.
\]
The matroid of
\[
A = \begin{bmatrix}
0 & 3 & 1 & 2 &  4 & 0 \\
0 & 0 & 0 &  1 & 2 & 1
\end{bmatrix} 
\]
is \[
\mathcal{M}_A = \{ \{2,4\}, \{2,5\}, \{2,6\}, \{3,4\}, \{3,5\}, \{3,6\}, \{4,6\}, \{5,6\}\} \subseteq \binom{[6]}{2}.
\]
\end{Ex}

Recall, the \textit{Johnson graph} $J(k,n)$ with vertices given by the
$k$-subsets of $[n]$ such that two $k$-subsets $I,J$ are connected by
an edge precisely if $|I\cap J|=k-1$.  For a matroid $\M \subseteq
\binom{[n]}{k}$, $J(\M)$ is the induced subgraph of the Johnson graph
on the vertices in $\M$.  The Basis Exchange Property for $\M$ implies
that $J(\M)$ is connected, and furthermore, between any two vertices $I,J$ in $J(\M)$, there exists a path in $J(\M)$ which is a minimal length path between $I$ and $J$ in $J(k,n)$.  The Johnson graph $J(\M_A)$ from \Cref{ex:matroid} is

\begin{center}
\includegraphics[width=6cm]{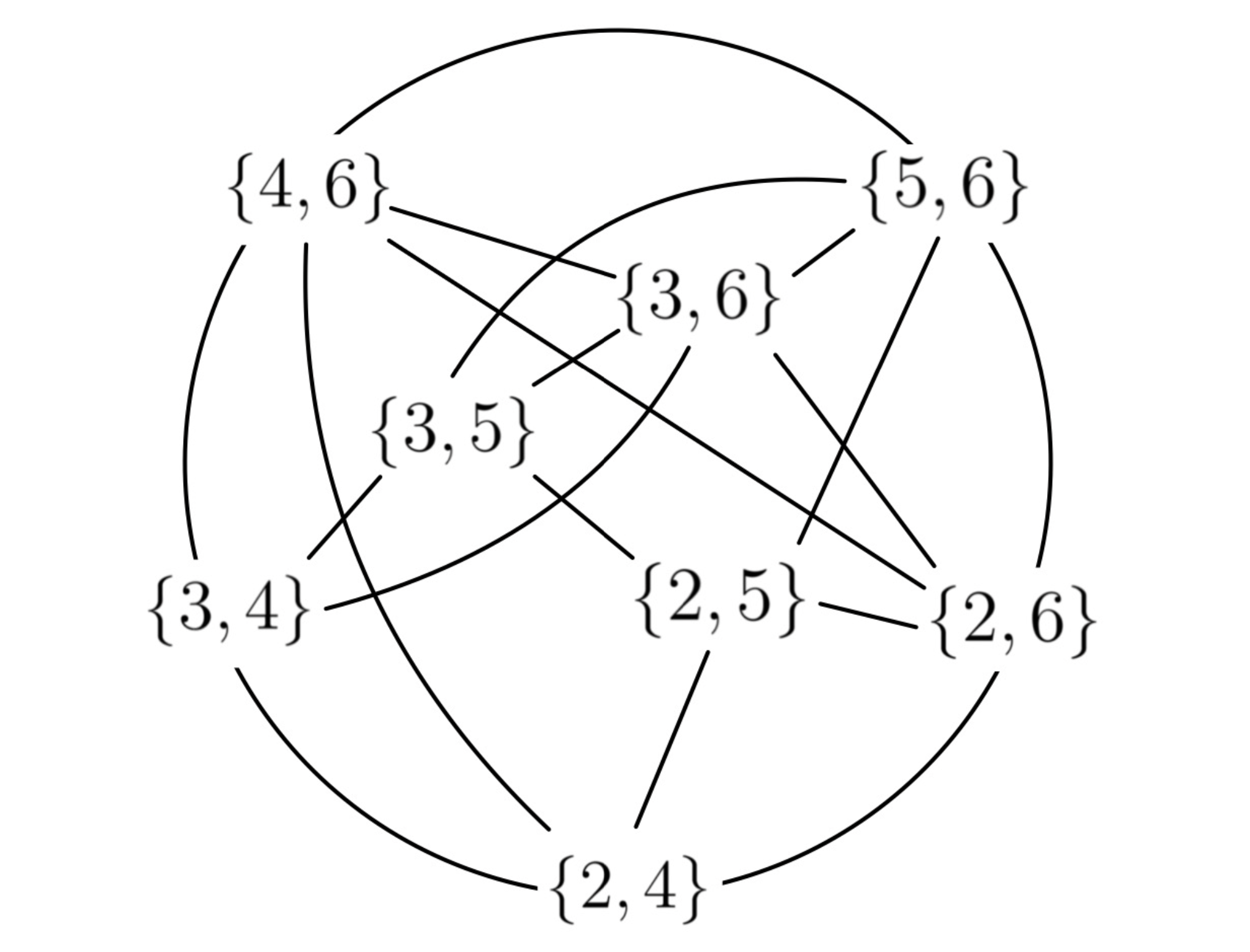}
\end{center}

The direct sum of two matroids on disjoint ground sets, denoted $M_{1} \oplus M_{2}$, is the
matroid with bases given by $\{I\cup J\ : \ I \in M_{1},J \in M_{2}
\}$ on the ground set which is the disjoint union of the ground sets
of the matroids $M_{1}$ and $M_{2}$.  A matroid $M$ on ground set $[n]$
is \textit{connected} if $M$ cannot be expressed as the direct sum of
two matroids. Every matroid can be decomposed into the direct sum of its
connected components.

Given a matroid $\M$ on ground set $[n]$, the \textit{dual matroid}
$M^{*}$ is the matroid on $[n]$ with bases $\{[n]\setminus I \subseteq
[n]: I \in \M \}$.  If $\M$ has rank $k$, then its dual matroid has
rank $n-k$.


Let $S_n$ be
the set of permutations of $[n]$, where we think of a permutation as a
bijection from a set to itself. For $w \in S_n$, let
$w_i=w(i)$, and write $w$ in \textit{one-line notation} as
$w=w_1w_2\cdots w_n$. For an interval of indices $[i,j]$, we denote the image of $[i,j]$ under $w$ by $w[i,j] = \{w_i, w_{i+1}, \ldots, w_j\}$. The
permutation matrix $M_w$ of $w$ is the $n \times n$ matrix that has a
1 in cell $(i,w_i)$ for each $i \in [n]$ and zeros elsewhere. The
\textit{length of} $w \in S_n$ is
\[
\ell(w) := \# \{(i,j) \; : \; i < j \; \text{and} \; w(i) > w(j)\}.
\]
The permutation in $S_n$ of maximal length is $w_0 := n (n-1) \cdots 2 1$.

\begin{Ex} For $w = 3124$, the length of $w = 3124$ is $\ell(w) = 2$, and $M_w$ is the matrix
\[
M_{3124} = 
\begin{bmatrix}
0 & 0 & 1 & 0 \\
1 & 0 & 0 & 0 \\
0 & 1 & 0 & 0 \\
0 & 0 & 0 & 1
\end{bmatrix}.
\]
\end{Ex}

\noindent Note that the permutation matrix of $w^{-1}$ is $M_{w}^{T}$.
Furthermore, permutation multiplication is given by function
composition so that if $wv=u$, then $w(v(i))= u(i) $. By this
definition, $M_w^TM_v^T= M_u^T$.

\begin{Def}\label{th:tableau.criterion} \cite[Chapter 2]{BjornerBrenti}
For $u,v$ in $S_{n}$,\ $u\leq v$ in \textit{Bruhat order} if
$u[i]\preceq v[i]$ for all $i \in [n]$.  Equivalently, Bruhat order is
the ranked poset defined as the transitive closure of the relation
$u<ut_{ij}$ whenever $u_{i}<u_{j}$, where $t_{ij}$ is the
permutation transposing $i$ and $j$ and fixing all other values.
For each $u\leq v$ in Bruhat order, the \textit{interval} $[u,v]$ is
defined to be
\begin{equation}\label{def:interval}
[u,v] :=\{y\in S_{n}\ : \ u \leq y\leq v \}.
\end{equation}
\end{Def}

For $0\leq k \leq n$, write $S_k \times S_{n-k}$ for the subgroup of
$S_n$ consisting of permutations that send $[k]$ to $[k]$ and
$[k+1,n]$ to $[k+1,n]$. If $k=0$, consider $[k]$ to be the empty set.
A permutation $w \in S_n$ is $k$-\textit{Grassmannian} if $w_1 <
\cdots < w_k$ and $w_{k+1} < \cdots < w_n$. This is equivalent to
saying that $w$ is the minimal length element of its coset $w \cdot
(S_k \times S_{n-k})$. For example, the permutation $w=3124$ is
1-Grassmannian.  The set of \textit{Grassmannian permutations} in
$S_{n}$ is the union over $k \in [0,n]$ of all $k$-Grassmannian permutations.

\begin{Def}\label{def:grassmann.internval} Assume $u\leq v$ in Bruhat
order on $S_{n}$.  Then, the interval $[u,v]$ is a \textit{Grassmann
interval} provided $v$ is a $k$-Grassmannian permutation for some $k
\in [0,n]$.  Denote by $\GIkn$ the set of all Grassmann intervals
$[u,v]$ in $S_n$, where $v$ is $k$-Grassmannian.
\end{Def}

The Grassmann intervals $[u,v]$ are key objects for this paper.  Note
that $u$ need not be a Grassmannian permutation.  In the case where $v$ is
Grassmannian, there is a simpler criterion for Bruhat order that
follows closely from work of Bergeron-Sottile \cite[Theorem
A]{BergSott}.  This criterion also appears in \cite[Lemma
20.2]{Postnikov.2006}. We will use this criterion extensively.

\begin{Th}\label{Lemma: BS}  Let $u,v \in
S_n$, where $v$ is $k$-Grassmannian.  Then $u \leq v$ if and only if
\begin{enumerate}
\item[(i)] for every $1 \leq j \leq k$, we have $u(j) \leq v(j)$, and
\item[(ii)] for every $k < m \leq n$, we have  $u(m) \geq v(m)$.  
\end{enumerate}
\end{Th}

\subsection{Grassmannian, Flag, and Richardson Varieties} \label{Sec: Fl(n) and Gr(k,n)}


Schubert varieties in both the flag variety and the Grassmannian are a
well studied class of varieties that have found importance in
algebraic geometry, representation theory, and
combinatorics. Singularities of Schubert varieties, in particular,
have been investigated thoroughly.  See \cite{Billey-Lak, Ful,Kumar-book} for
further background on these varieties.  

For $0 \leq k \leq n$, the points in the \textit{Grassmannian
variety}, $\Gkn$, are the $k$-dimensional subspaces of $\C^n$. Up to
left multiplication by a matrix in $\GL_k$, we may represent $V \in
\Gkn$ by a full rank $k \times n$ matrix $A_V$ such that $V$ is the
row span of $A_V$. Let $\mathrm{Mat}_{kn}$ be the set of full rank $k
\times n$ matrices. The points in $\Gkn$ can be bijectively identified
with the cosets $\GL_k \setminus \mathrm{Mat}_{kn}$.  The
Grassmannian varieties are smooth manifolds via the Pl{\"u}cker coordinate embedding of $\Gkn$ into projective space. This includes the case
when $k=n=0$, in which case $\Gkn$ consists of one point, which is
the $0$-dimensional vector space in $\C^0$.  

Let $\mathcal{F}\ell(n)$ be the \textit{complete flag variety} of
nested subspaces of $\C^n$.  A complete flag $V_\bullet = (0 \subset V_1
\subset \cdots \subset V_n)$ can be represented as an invertible $n \times n$
matrix where the row span of the first $k$ rows is the subspace
$V_{k}$ in the flag.  Throughout the paper, we will often identify a
full rank matrix with the point it represents in $\Gkn$ or
$\mathcal{F}\ell(n)$. For a subset $J \subseteq [n]$, let Proj$_J :
\C^n \rightarrow \C^{|J|}$ be the projection map onto the indices
specified by $J$. Then, for every permutation $w \in S_n$, there is a
\textit{Schubert cell} $C_w$ and an \textit{opposite Schubert cell} $C^w$ in
$\mathcal{F}\ell(n)$ defined by
\begin{align*}
&C_w = \{V_\bullet \in \mathcal{F}\ell(n) \; : \; \dim(\text{Proj}_{[j]}(V_i)) = |w[i] \cap [j]| \;\text{for all $i,j$}\}, \\
&C^w = \{V_\bullet \in \mathcal{F}\ell(n) \; : \; \dim(\text{Proj}_{[n-j+1,n]}(V_i)) = |w[i] \cap [n-j+1,n]| \;\text{for all $i,j$}\}.
\end{align*}
By row elimination and rescaling, we can find a \textit{canonical
matrix representative} $A_{V_\bullet} := (a_{i,j})$ for each
$V_{\bullet} \in C_{w}$ such that $a_{i,w_i} = 1$ for $i\in [n]$,
$a_{i,j} = 0$ for all $1\leq j < w_i$, and $a_{h,w_i} = 0$ for $h>i$.
Canonical matrices for $C^w$ can be found similarly, but so that every leading
1 has all zeros to the right instead of the left.  For example, the
canonical matrices for $C_{3124}$ have the form below, where entries labeled $*$ can be replaced by any element of
$\mathbb{C}$. 
\[
\begin{bmatrix}
0 & 0 & 1 & * \\
1 & * & 0 & * \\
0 & 1 & 0 & * \\
0 & 0 & 0 & 1
\end{bmatrix}
\]

The \textit{Schubert variety} $X_w$ is the closure of $C_w$ in the
Zariski topology on $\mathcal{F}\ell(n)$, and similarly, the
\textit{opposite Schubert variety} $X^w$ is the closure of $C^w$. Specifically, $X_w$ and $X^w$ can be defined by
\begin{align*}
&X_w = \{V_\bullet \in \mathcal{F}\ell(n) \; : \; \dim(\text{Proj}_{[j]}(V_i)) \leq |w[i] \cap [j]| \;\text{for all $i,j$}\}, \\
&X^w = \{V_\bullet \in \mathcal{F}\ell(n) \; : \; \dim(\text{Proj}_{[n-j+1,n]}(V_i)) \leq |w[i] \cap [n-j+1,n]| \;\text{for all $i,j$}\}.
\end{align*}
Bruhat order determines which Schubert cells are in a Schubert variety, 
\begin{equation}\label{eq:schubs}
 X_w = \bigsqcup_{y \geq w} C_{y} \ \text{ and } \ X^{w} = \bigsqcup_{v
\leq w} C^{v}.
\end{equation}

Schubert varieties in the Grassmannian manifold $\Gkn$ are indexed by the sets in $\binom{[n]}{k}$ and can be defined
as the projections of Schubert varieties in $\mathcal{F}\ell(n)$.  Let
\begin{equation}\label{eq:projection}
\pi_k : \mathcal{F}\ell(n) \rightarrow Gr(k,n)
\end{equation}
be the projection map which sends a flag $V_\bullet = (0 \subset V_1
\subset \cdots \subset V_n)$ to the $k$-dimensional subspace
$V_k$. Identifying a full rank $n \times n$ matrix $M$ with the point
it represents in $\mathcal{F}\ell(n)$, then $\pi_k(M)$ denotes the
span of the top $k$ rows of $M$.  For each $J \in \binom{[n]}{k}$,
there exists a $k$-Grassmannian permutation $v$ such that $J=v[k]
:=\{v_{1},v_{2},\dots, v_{k}\}$.  The \textit{Grassmannian Schubert variety}
indexed by $J$ is the projection $\pi_{k}(X_{v})$.   

For permutations $u$ and $v$ in $S_n$, with $u \leq v$, the
\textit{Richardson variety} is a nonempty variety in
$\mathcal{F}\ell(n)$ and is defined as the intersection $X_u^v := X_u
\cap X^v$. Then $\dim X_u^v = \ell(v) - \ell(u)$. 
The decompositions of $X_u$ and $X^v$ into Schubert cells and opposite Schubert cells in \eqref{eq:schubs} yield
\[
X_u^v = {\displaystyle \bigsqcup_{u \leq y \leq v} (C_y \cap X^v)}
=
\Big({\displaystyle \bigsqcup_{y \geq u}} C_y \Big) \cap
\Big({\displaystyle \bigsqcup_{t \leq v}} C^t\Big).
\]
By \Cref{eq:schubs}, one can observe that the permutation matrix
$M_{y} \in (C_y \cap X^v) \subset X_u^v$ for each $y \in [u,v]$.  

In \cite{Lak-San}, Lakshmibai and Sandhya prove that $X_w$ is smooth
if and only if $w$ avoids 1324 and 2143 as permutation
patterns. Related work was also done by Ryan \cite{Ry}, Wolper
\cite{Wol}, and Haiman \cite{Haiman}.  In \cite[Thm 2.4]{BP}, Billey
and Postnikov extend this result by giving a pattern avoidance
criterion for smoothness of Schubert varieties for all Weyl
groups. Given a singular Schubert variety, $X_w$, \cite[Thm 1]{BW},
\cite[Thm 1.3]{KLR}, \cite[Thm 2]{Man}, and \cite{Cortez} developed
criteria based on permutation patterns to determine the singular
points of $X_w$.  These theorems characterizing smooth versus singular
points in Schubert varieties using permutation patterns motivated this work.  

Singularities in Richardson varieties and their projections have also
been studied in the literature
\cite{Brion:flag,BC,KLS2,kreiman-lakshmibai}.  The characterizations
of smooth versus singular Richardson varieties described there are not
based on pattern avoidance but rely on computations in the associated
cohomology rings.


\subsection{Decorated Permutations}\label{Sec: decperms}

A decorated permutation $\w$ is defined by a permutation $w$
together with a circular orientation map called $\co$ from the fixed
points of $w$ to the set of clockwise or counterclockwise orientations, denoted by $\{\cw, \ccw\}$. Therefore, we will sometimes
describe a decorated permutation as a pair $(w,\co)$.

Postnikov made the following definitions in \cite[Sect
16]{Postnikov.2006}.  Given a decorated permutation $\w$, call $i \in
[n]$ an \textit{anti-exceedance} of $\w=(w,\co)$ if $i<w^{-1}(i)$ or
if $\w(i) = \clockwise{i}$ is a clockwise fixed point.  If $i \in [n]$
is not an anti-exceedance, it is an \textit{exceedance}.  Let
$I_{1}(\w)=I_1$ be the \textit{set of anti-exceedances} of $\w$.
For an arc $(i \mapsto w(i))$ in the chord diagram of $\w$, we say the
arc is an \textit{anti-exceedance arc} or \textit{exceedance arc}
depending on whether $w(i)$ is in $I_{1}(\w)$ or not.  Let $k(\w) :=
|I_1(\w)|$.  Recall that $\SSnk$ is the set of decorated permutations
with anti-exceedance set of size $k$.

More generally, let $<_r$ be the \textit{shifted linear order} on $[n]$ given
by $r <_r (r+1) <_r \cdots <_r n <_r 1 <_r \cdots <_r (r-1)$ for $r
\in [n]$.  The \textit{shifted anti-exceedance set $I_r(\w)$} of $\w$
is the anti-exceedance set of $\w$ with respect to the shifted linear
order $<_r$ on $[n]$,
\[
I_r(\w) = \{i \in [n] \; : \; i <_r w^{-1}(i) \; \text{or} \; \w(i) =
\clockwise{i} \}.
\]
An element $i \in I_r(\w)$ is called an \textit{$r$-anti-exceedance}, and an element $i \notin I_r(\w)$ is called an \textit{$r$-exceedance}. Note from the
construction that either $I_{r+1}(\w) = I_{r}(\w)$ or $I_{r+1}(\w) =
I_{r}(\w) \setminus \{r\} \cup \{w(r) \}$, so
$|I_1(\w)|=\cdots=|I_n(\w)|=k(\w)$.  Furthermore, $I_{r+1}(\w) = I_{r}(\w) $
if and only if $r$ is a fixed point of $w$, so clockwise fixed points
will be in all of the shifted anti-exceedance sets and
counterclockwise fixed points will be in none.  Therefore, $\w$ is
easily recovered from $(I_1(\w), \ldots, I_n(\w))$.  The sequence of
$k(\w)$-subsets $(I_1(\w), \ldots, I_n(\w))$ is called the \textit{Grassmann
necklace} associated to $\w$.

Given a decorated permutation $\w=(w,\co)\in \SSnk$ and its
anti-exceedance set, we can also easily identify the Grassmann
interval associated to it.  The $k$-Grassmannian permutation $v$ is
determined by $v[k]=w^{-1}(I_{1}(\w))$, and then $u$ is
determined by $u=wv$, and hence $I_{1}(\w) = u[k]$ by this
construction. Let $u(\w)=u$ and $v(\w)=v$.  The interval $[u,v]\in
\GIkn$ is the \textit{Grassmann interval} associated to $\w \in
\SSnk$.  To identify a decorated permutation from a Grassmann interval
$[u,v] \in \GIkn$, simply let
\begin{equation}\label{eq:uv.to.dp.bijection}
w = u v^{-1} \quad \text{with} \quad
co(j) := \begin{cases}
\cw & \text{if  } j \in u[k] \\
\ccw &  \text{if  } j \notin u[k]
\end{cases}.
\end{equation}

We visualize the bijection from decorated permutations to Grassmann
intervals as a shuffling algorithm.
\begin{enumerate}
\item Write $\w$ in two-line notation with the
numbers $1,2,\ldots,n$ on the first row and $w(1),w(2),\ldots,w(n)$ on
the second row, including orientations on fixed points.  
\item Highlight the columns $i_{1}<i_{2}<\cdots<i_k$ such
that $i_{j}>w(i_{j})$ or $\w(i_{j}) = \clockwise{i_j}$, corresponding with anti-exceedances. 
\item Keeping the columns intact, reorder the columns so that the
highlighted columns $i_{1}<i_{2}<\cdots<i_k$ come first followed by
the non-highlighted columns maintaining their relative order within the highlighted and non-highlighted blocks. Then, drop
any orientation arrows to obtain a $2\times n$ array of positive
integers $\uv$ with the one-line notation for $u$ determined by the lower row and the one-line notation for $v$ determined
by the upper row.
\end{enumerate}
From this shuffling algorithm, note that $v$ is a
$k$-Grassmannian permutation by construction and $v[k]=w^{-1}\cdot
I_{1}(\w)$.  Furthermore, $u\in S_{n}$ and its initial set $u[k] =
I_{1}(\w)$.  Also observe that $u_{i}\leq v_{i}$ for all $1\leq i \leq k$, and $u_{i}\geq
v_{i}$ for all $k+1\leq i \leq n$. Hence, $u\leq v$ in Bruhat order by
\Cref{Lemma: BS}.  Thus, $[u,v]$ is the Grassmann interval associated
to $\w$, since the shuffling algorithm above is equivalent
to the permutation multiplication $wv=u$.

Conversely, to go from $[u,v]$ to $\w=(w,\co)$, simply take the two
line array $\uv$, sort the columns by the top row to obtain $w$.
Observe that $j$ is a fixed point of $w$ if and only if $u^{-1}(j) =
v^{-1}(j)$. If $j$ is a fixed point of $w$, then $\co(j) = \cw$ 
if $j \in u[k]$, and $co(j) = \ccw$ otherwise.

\begin{Rem}
Observe that the chord diagram is just as easily obtained from
$[u,v]$ as it is from of $\w$.  The chord diagram arcs $\{(i
\mapsto w(i)) : i \in [n] \}$ are determined by the columns in the two
line notation of $\w$, so $\{(i \mapsto w(i)) : i \in [n] \} =
\{(v_{j} \mapsto u_{j}) : j \in [n] \}$ with fixed points oriented
appropriately.  Therefore, the graphical patterns determining
properties for decorated permutations also determine patterns for
Grassmann intervals.
\end{Rem}

\begin{Rem}\label{rem:apost.uv.bij} Note that the map above between
decorated permutations in $\SSnk$ and Grassmann intervals in
$\GIkn$ follows the work of \cite{KLS}. In
\cite[Sect. 20]{Postnikov.2006}, a slightly different map is given. In
particular, in Postnikov's work, the decorated permutation $\w$
obtained from a Grassmann interval $[u,v]$ is computed as $w = w_0
v u^{-1} w_0$, with a fixed point $i$ assigned a clockwise orientation
if $n+1-i$ is in $u[k]$.  Here $w_{0}=n\ldots 2 1 \in S_{n}$.  We will
return to this involution on decorated permutations and associated
objects in \Cref{rem:connectbacktoPost}.  Our reason for prioritizing
the Knutson-Lam-Speyer bijection is the direct connection to the
corresponding positroid given in \Cref{Th: positroid = initial sets}.
\end{Rem}

\begin{Ex}\label{ex:psi-map}

For the decorated permutation $\w =54127\clockwise{6}9\counterclockwise{8}3$, the anti-exceedance set
is $I_{1}(\w)=\{1,2,6,3 \}$.  These values occur in positions
$\{3,4,6,9 \}$ in $\w$.  The Grassmann necklace for $\w$ is 
\[
 (I_{1},\ldots, I_{9})=
(\{1236\},\{2356\},\{3456\},\{1456\},\{1256\},\{1267\},\{1267\}, 
\{1269\},\{1269\}).  \]
Write the two-line notation, highlighting the anti-exceedances and then 
shuffling the anti-exceedances to the front to identify the associated Grassmann interval $[u,v] = [126354798,
346912578]$,
\[
\w=\left[\begin{array}{cc>{\columncolor{yellow}}c>{\columncolor{yellow}}cc>{\columncolor{yellow}}ccc>{\columncolor{yellow}}c}
1&2&3&4&5&6&7&8&9\\
5&4&1&2&7&\clockwise{6}&9&\counterclockwise{8}&3
\end{array} \right]
\
\implies
\
\left[\begin{array}{c}
v
\\
u
\end{array} \right]
=
\left[\begin{array}{>{\columncolor{yellow}}c>{\columncolor{yellow}}c>{\columncolor{yellow}}c>{\columncolor{yellow}}cccccc}
3&4&6&9&1&2&5&7&8
\\
1&2&6&3&5&4&7&9&8 
\end{array} \right].
\]

Using Postnikov's map, the decorated permutation corresponding to the
interval $[u,v] = [126354798, 346912578]$ is
$3 \overleftarrow{2} 5 \overrightarrow{4} 98167$.  Postnikov's
inverse map would associate the Grassmann interval $[416732598,
478912356]$ to the original $\w = 54127\clockwise{6}9\counterclockwise{8}3$.
\end{Ex}

A key artifact of a decorated permutation in our work is its set of
alignments. To formally define alignments, we first establish the following
notation for a cyclic interval of elements in $[n]$ for a fixed
integer $n$.

\begin{Def} \label{Def:cyclic.interval}
Let $a,b \in [n]$. Then
\begin{align}
[a,b]^{cyc} := \begin{cases}
[a,b] & \text{ if } a \leq b \\
[a,n] \cup [1,b] & \text{ if } a > b
\end{cases}.
\end{align}
\end{Def}

\begin{Def} \label{Def: algnmt}
An \textit{alignment} of $\w =(w,co) \in \SSn$ is a pair of arcs $(p \mapsto w(p))$ and $(s \mapsto w(s))$ in $D(\w)$ which can be drawn as distinct noncrossing
arcs such that
\begin{enumerate}
\item both $w(p) \in [p, w(s)-1]^{cyc}$ and $w(s) \in [w(p) +1, s]^{cyc}$,
\item if $w(s) = s$, then $co(s) = \cw$, and
\item if $w(p) = p$, then $co(p) = \ccw$.
\end{enumerate}
In this case, we denote the alignment by $(p \mapsto w(p), s \mapsto w(s))$ or $A(p,s)$ if
$\w$ is understood from context. We say the arc $(p \mapsto w(p))$ is
the \textit{port side} and the arc $(s \mapsto w(s))$ is the
\textit{starboard side} of the alignment.  See \Cref{fig: port
starboard}. We use $\aligns(\w)$ to denote the set of all alignments
of $\w$.
\end{Def}

\begin{minipage}[t]{1.0\linewidth}
    \centering
   \vspace{-2ex}
   \includegraphics[scale=0.08]{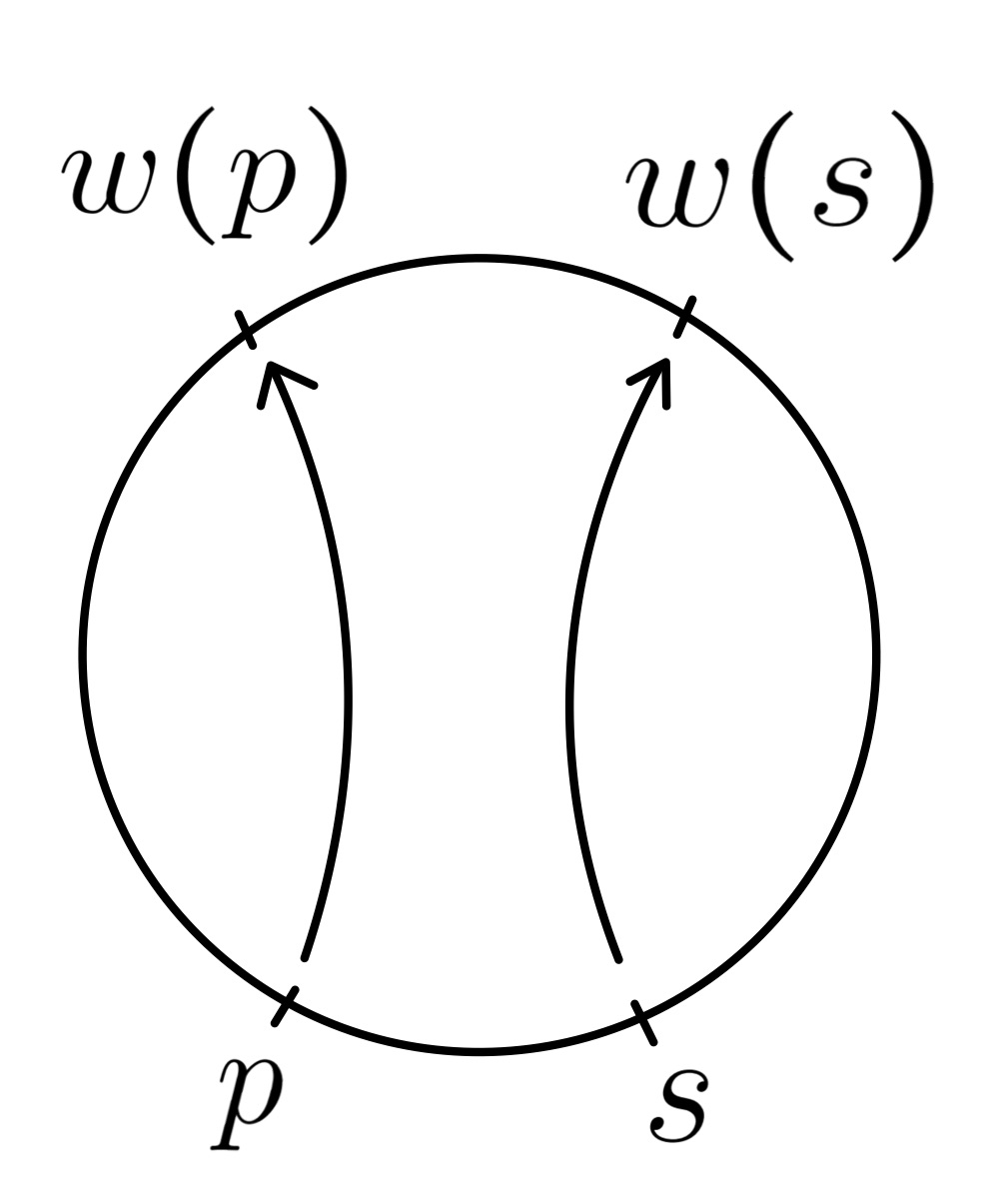}
   \captionof{figure}{Alignment with port side $(p \mapsto w(p))$ and starboard side $(s \mapsto w(s))$.}
   \label{fig: port starboard}
\end{minipage}

\begin{Ex}\label{ex:alignments}
Recall the chord diagram in \Cref{Example:intro} for $\w =
895\counterclockwise{4}7\clockwise{6}132$, the alignment $(9\mapsto 2,
8 \mapsto 3)$ highlighted in yellow has port side $(9\mapsto 2)$ and
starboard side $(8 \mapsto 3)$ as if they were two sides of a boat
with its bow pointing to the right.   Furthermore, $\w$ has 13 alignments,
\begin{align*}
\aligns(\w) = \{&A(3,1), A(3,2), A(3,6), A(4,1), A(4,2),A(4,3), A(4,6),\\
&A(5,1), A(5,2),A(7,6),A(8,6), A(9,6), A(9,8) \}.
\end{align*}

\end{Ex}


\begin{Def} \label{Def: crossed algnmt}
A \textit{crossed alignment} $A(p,s,x)$ consist of an alignment
$A(p,s)$ of $\w$ and an additional arc $(x \mapsto w(x))$ crossing both $(p
\mapsto w(p))$ and $(s \mapsto w(s))$.  We partition the set of crossed alignments
according to which side of the alignment the crossing arc intersects first as it passes from $x$ to $w(x)$.

\begin{enumerate}
\item If $x \in [w(s),s]^{cyc}$ and $w(x) \in [p,w(p)]^{cyc}$,  we say 
$A(p,s,x)$ is a \textit{starboard tacking crossed alignment}.
\item If $x \in [p,w(p)]^{cyc}$ and $w(x) \in [w(s),s]^{cyc}$, we say
$A(p,s,x)$ is a \textit{port tacking crossed alignment}.
\end{enumerate}
\end{Def}
Note that if $A(p,s,x)$ is a crossed alignment, then $p$ and $s$
cannot be fixed points. The highlighted crossed alignment in
\Cref{Example:intro} $A(9,8,7)$ is an example of a starboard
tacking crossed alignment, while $A(9,8,1)$ is a port tacking
crossed alignment.

The following family of (decorated) permutations are key to this work.
As a subset of $S_{n}$, these permutations have nice enumerative
properties \cite{Callan.2004}, \cite[A075834]{OEIS}.

\begin{Def}\label{def:SIFperms}\cite{Callan.2004}
A permutation $w \in S_{n}$ is \textit{stabilized-interval-free}
provided no proper nonempty interval $[a,b] \subset [n]$ exists such that
$w[a,b]=[a,b]$. 
\end{Def}

Note, the only stabilized-interval-free permutation with a fixed
point is the identity permutation in $S_{1}$.  Thus, the definition of
stabilized-interval-free permutations extends easily to
\textit{decorated stabilized-interval-free permutations}.  Both decorated
permutations in ${S_1^{\circ \bullet}}$ are stabilized-interval-free,
and for $n\geq 2$ the SIF permutations and the decorated SIF permutations are
the same.

\subsection{Positroids and Positroid Varieties}\label{sub:positroids}

Postnikov and Rietsch considered an important cell decomposition of
the totally nonnegative Grassmannian
\cite{Postnikov.2006,Rietsch.2006}.  The term \textit{positroid} does
not appear in either paper, but has become the name for the matroids
that index the nonempty matroid strata in that cell decomposition.
They also individually considered the closures of those cells, which
determines an analog of Bruhat order.  The cohomology classes for
these cell closures was investigated by Knutson, Lam, and Speyer
\cite{KLS,KLS2} and Pawlowski \cite{Paw}.  Further geometrical
properties of positroid varieties connected to Hodge structure and
cluster algebras can be found in \cite{Galashin-Lam.2021, Lam.2019}.

\begin{Def} \label{Def:tnn}
\normalfont A real valued $k\times n$ matrix $A$ is \textit{totally
nonnegative} (tnn) if each maximal minor $\Delta_I(A)$ satisfies
$\Delta_{I}(A)\geq 0$ for $I \in \binom{[n] }{k}$.  Let $\Gkn^{tnn}$
be the points in $\Gkn$ that can be represented by totally nonnegative
matrices.
\end{Def}

\begin{Def} \label{Def: positroid}
A \textit{positroid} of rank $k$ on the ground set $[n]$ is a matroid
of the form $\M_A$ for some matrix $A \in \Gkn^{tnn}$. More generally, a positroid can be defined on any ordered ground set $B = \{b_1 < \cdots < b_n\}$.
\end{Def}

Positroids are an especially nice class of realizable matroids. For
example, they are closed under the matroid operations of restriction,
contraction, and duality, as well as a cyclic shift of the ground set
\cite{Ardila-Rincon-Williams}.  The following theorem is nice
connection between positroids and non-crossing partitions as explored
by Ardila, Rinc\'on and Williams.  See also \cite{Ford.2015} for one
direction of this theorem.  

\begin{Th}\label{thm:ARW}\cite[Thm. 7.6]{Ardila-Rincon-Williams}
The connected components of a positroid on ground set $[n]$ give rise
to a non-crossing partition of $[n]$.  Conversely, each positroid $\M$
on $[n]$ can be uniquely constructed by choosing a non-crossing
partition $B_{1} \sqcup \cdots \sqcup B_{t}$ of $[n]$, and then putting the
structure of a connected positroid $\M_i$ on each block $B_{i}$, so
$\M= \M_{1} \oplus \cdots \oplus \M_{t}$.
\end{Th}

The connected components appearing in the decomposition of $\M$ above each correspond to decorated SIF permutations. This relationship is given in the following corollary.

\begin{Cor}\label{cor:SIFcor}
The decomposition of a positroid $\M= \M_{1} \oplus \cdots
\oplus\M_{t}$ into connected components gives rise to a
noncrossing partition $[n]= B_{1} \sqcup \cdots \sqcup B_{t}$ and a
decomposition of the associated decorated permutation $\w = w^{(1)}
\oplus \cdots \oplus w^{(t)}$, where each $w^{(i)}$ is a
decorated stabilized-interval-free permutation on $B_{i}$, and the
converse holds as well.
\end{Cor}

As mentioned in the introduction, positroids are in bijection with
decorated permutations.  The following bijection goes by way of the
associated Grassmann necklace.  The bijection depends on the
\textit{shifted Gale order} $\prec_r$ on $\binom{[n]}{k}$ using the
shifted linear order $<_{r}$ on $[n]$.  Specifically, if $I$ and $J$
are $k$-sets that can be written under $<_r$ as $I = \{i_1 <_r \cdots <_r i_k\}$
and $J = \{j_1 <_r \cdots <_r j_k\}$, then $I \preceq_r J$ if $i_h
\leq_r j_h$ for all $h \in [k]$.

\begin{Th}\label{th:decoratedpermtopositroid}\cite{Oh,Postnikov.2006}
For $\w \in \SSnk$, the set
\begin{align} \label{Def: positroid def from anti-ex sets}
\M(\w) \; := \; \Big\{I \in \binom{[n]}{k} \; : \; I_r(\w) \preceq_r I \text{ for all } r \in [n]\Big\}
\end{align}
is a positroid of rank $k$ on ground set $[n]$. Conversely, for every
positroid $\M$ of rank $k$ on ground set $[n]$, there exists a unique
decorated permutation $\w\in \SSnk$ such that the sequence of minimal
elements in the shifted Gale order on the subsets in $\M$ is the
Grassmann necklace of $\w$.
\end{Th}

\begin{Cor}\label{cor:allIrs are in positroid}
For $\w=(w,\co) \in \SSnk$ with associated positroid $\M(\w)$, every shifted
anti-exceedance set $I_r(\w)$ is in $\M(\w)$.   
\end{Cor}

\begin{ThmDef} \label{Def: positroid as projection}\cite[Thm 5.1]{KLS}
Given a decorated permutation $\w \in \SSnk$ along with its
associated Grassmann interval $[u,v]$ and positroid $\M = \M(\w) \subseteq
\binom{[n]}{k}$, the following are equivalent definitions of the
\textit{positroid variety} $\Pi_{\w} = \Pi_{[u,v]} = \Pi_{\M}$.  
\begin{enumerate}
\item The positroid variety $\Pi_{\M}$ is the homogeneous subvariety
of $Gr(k,n)$ whose vanishing ideal is generated by the Pl{\"u}cker coordinates indexed by the nonbases of $\M$, $\{\Delta_I \; : \; I \in \QQ(\M)\}$.
\item The positroid variety $\Pi_{[u,v]}$ is the projection of the
Richardson variety $X_u^v \subseteq \mathcal{F}\ell(n)$ to $Gr(k,n)$, so $\Pi_{[u,v]}=\pi_k(X_u^v)$.
\end{enumerate}
\end{ThmDef}


\begin{Ex}\label{ex:positroid}
In \Cref{ex:matroid}, the matrix $A$ has all nonnegative $2\times 2$
minors, so the associated matroid is a positroid.  The minimal elements
in shifted Gale order are $(\{2 4\}, \{2 4\}, \{3 4\}, \{4 6\}, \{5
6\},$ $\{2 6\})$, which is the Grassmann necklace for the decorated
permutation $\counterclockwise{1}36524$.  The associated Grassmann
interval is $[241365, 561234]$.  The set of nonbases of $\M$ is
\[
\QQ = \{\{1,2 \}, \{1,3 \}, \{1,4 \}, \{1,5 \}, \{1,6 \}, \{2,3 \},\{4,5 \} \}.
\]
Therefore, the points in the positroid variety $\Pi_{\w}$ are
represented by the full rank complex matrices of the form
\[
\begin{bmatrix}
0 & a_{12} & ca_{12} & a_{14} & da_{14} & a_{16} \\
0 & a_{22} & ca_{22} & a_{24} & da_{24} & a_{26}
\end{bmatrix}. 
\]
\end{Ex}

As mentioned in the introduction, there are many other objects in
bijection with positroids and decorated permutations. We refer the
reader to \cite{Ardila-Rincon-Williams} for a nice survey of many
other explicit bijections.

\section{Positroid Characterization Using Initial Sets}\label{sec:initial.sets}

The definition of positroid varieties in \Cref{Def: positroid as projection} is given in terms of positroids/positroid complements and
the Grassmann intervals.  In this section, we show how the initial
sets of permutations in Grassmann intervals determine the positroids,
proving \Cref{Th: positroid = initial sets}.  This theorem completes
the commutative diagram below.

\begin{center}
\begin{tikzcd}[sep=huge] 
\SSnk \arrow[r, "shuffle"] \arrow[d,  "necklace"]
& \GIkn \arrow["initial" red, d ] \\
\text{Grassmann necklaces in } \binom{[n] }{k}^{n} \arrow[r, "\hspace{.2in}
shifted Gale \hspace{.2in}" ]
& \text{Positroids in } \binom{[n] }{k}
\end{tikzcd}
\end{center}

\subsection{Canonical Representatives}
\label{Sec: Canonical reps}

For a $k$-set $I=\{i_1< i_{2}< \dots < i_{k}\}$ such that $i_{j}\leq
v_{j}$ for all $1\leq j\leq k$, we will show there exists a unique
maximal element below $v$ in each $S_{k} \times S_{n-k}$ coset with
initial set $I$. Let $u(I,v) \leq v$ be this maximal element when it
exists.  We think of $u(I,v)$ as the \textit{canonical representative}
of $I$ in $[id,v]$.

\begin{Lemma}\label{lem:maximal.reps} \normalfont Let $u\leq v$ be two
permutations in $S_{n}$, and assume $v$ is $k$-Grassmannian.  Then
there exists a unique permutation $u' \in [u,v]$ such that $u'$ is
maximal in Bruhat order among all permutations in the interval $[u,v]$
with $[u'_{1},\dotsc , u'_{k}]=[u_{1},\dotsc , u_{k}]$ fixed as an
ordered list.  Similarly, there exists a unique permutation $u' \in
[u,v]$ such that $u'$ is maximal in Bruhat order among all
permutations in the interval $[u,v]$ with $[u'_{k+1},\dotsc ,
u'_{n}]=[u_{k+1},\dotsc , u_{n}]$ fixed.
\end{Lemma}

\begin{proof}
Let $u'=u'_{1},\dotsc , u'_{n}$ be defined by $u'_{a}:=u_{a}$ for $a \leq k$ and for $b=n,n-1,\ldots,k+1$ by
\begin{equation}\label{eq:def.u.max1}
u'(b) = \min \; [v(b),n]\setminus(u'[k]\cup u'[b+1,n]).
\end{equation}
We claim that $u'$ is well defined. To see this, recall that $u \leq v$ implies $u(b) \geq v(b)$ for all $b>k$ by \Cref{Lemma: BS}, since $v$ is $k$-Grassmannian. Hence, for all $b'$ such that $k < b \leq b' \leq n$, we have $u(b') \geq v(b') \geq v(b)$, so $|[v(b),n] \cap u[k+1,n]| \geq n-b+1$. Then, at each step, the set $[v(b),n] \setminus (u'[k] \cup u'[b+1,n]) = ([v(b),n] \cap u[k+1,n]) \setminus u'[b+1,n]$ is nonempty.

By construction, $u'(a) =u(a) \leq v(a)$ for all $a \leq k$ and $u'(b)
\geq v(b)$ for all $b>k$, so again by \Cref{Lemma: BS}, we know
$u'\leq v$.  To see that $u\leq u'$, observe that we can obtain $u'$
from $u$ by applying a sequence of transpositions which increase
length at each step.  Since $u'(n)=u(b)$ for some $b\in [k+1,n]$, then
$u(b)\leq u(n)$ by \eqref{eq:def.u.max1}. Therefore, $u<ut_{bn}$ by
\Cref{th:tableau.criterion}.  Next, $u'(n-1)=ut_{bn}(c)$ for some
$c\in [k+1,n-1]$ such that $ut_{bn}(c)\leq ut_{bn}(n-1)$ by
\eqref{eq:def.u.max1}. Hence $u<ut_{bn}<ut_{bn}t_{c(n-1)}$, etc. for
$b=n-2,\ldots,k+1$.

Consider the permutation matrix $M_{u'}$ below row $k$.  For each
$k<b\leq n$, we see from the construction of $u'$ that the rectangle
northeast of $(b,v(b))$ and southwest of $(k+1,n)$ contains a
decreasing sequence of 1's ending at $u'(b)\geq v(b)$.  Therefore,
there is no $a<b$ such that $u'<u't_{ab}\leq v$ unless $a\leq k$ which
would change the first $k$ values of $u'$.  So, $u'$ is maximal among
all elements below $v$ with prefix $[u_{1},\dotsc , u_{k}]$.

It remains to show that $u'$ is the unique maximal element below $v$
with prefix $[u_{1},\dotsc , u_{k}]$.  Say $x\in S_{n}$ is another
element such that $u\leq x\leq v$ and $u_{i}=x_{i}$ for all $1 \leq
i\leq k$.  If $u' \neq x$, then there exists a maximal $j$ such that
$u'(j) \neq x(j)$. Since $x\leq v$, we know $x(j) \in [v(j),n]$ by
\Cref{Lemma: BS} again.  By construction $u'(j)$ is the minimal
element in
\[
[v(j),n]\setminus(u'[k] \cup u'[j+1,n] \}) = [v(j),n]\setminus(u'[k] \cup
x[j+1,n]),
\]
so $u'(j)<x(j)$.  Since $u,\ u'$, and $x$ are bijections, there exists
some $k<i<j$ such that $x(i)=u'(j)$.  Therefore, $x<xt_{ij}$, 
$x(i)\geq v(j)$, and $x(j)\geq v(i)$, so $xt_{ij}\leq v$ by \Cref{Lemma:
BS}.  Hence, $x$ is not a maximal element below $v$ with initial
values $[u(1),\dotsc , u(k)]$.

The proof for the second claim is very similar.  We leave the details
to the reader.
\end{proof}

\begin{Def} \label{def: uAv}  
Given a $k$-Grassmannian permutation $v \in S_{n}$ and a $k$-subset
$I=\{i_1< i_{2}< \dots < i_{k}\} \in \binom{[n] }{k}$ such that
$i_{j}\leq v_{j}$ for all $1\leq j\leq k$, let $u'=u(I,v)$ be the
permutation obtained as follows.

\begin{itemize}
\item Loop for $j$ from $1$ to $k$, assuming $u'(1),u'(2),\dotsc , u'(j-1)$ are defined. Let
\begin{equation}\label{eq:def.u.max2}
u'(j) = \max \; \{i \in I: i\leq v_{j} \}\setminus u'[j-1].
\end{equation}

\medskip

\item Loop for $j$ from $n$ down to $k+1$, assuming $u'(n),u'(n-1),\dotsc , u'(j+1)$ are defined. Let
\begin{equation}\label{eq:def.u.max3}
u'(j) = \min \; [v(j),n]\setminus(u'[k] \cup u'[j+1,n]).
\end{equation}
\end{itemize}

\end{Def}

\begin{Cor} \label{cor: uAv} 
Let $v$ be a $k$-Grassmannian permutation, and let $u\leq v$.  If
$I=\{u_{1},\ldots,u_k \}$ then $u \leq u(I,v) \leq v$, and $u(I,v)$
is the unique maximal element in the $S_{k}\times S_{n-k}$ right coset
containing $u$.   
\end{Cor}

\begin{proof}
Use \Cref{lem:maximal.reps} to find the maximal permutation $u\leq
u'\leq v$ with the prefix determined by $u$.  Then use
\Cref{lem:maximal.reps} again to find the maximal permutation below $v$
with suffix determined by $u'$ to obtain $u(I,v)$.  
\end{proof}

\begin{Rem}\label{rem:Oh-Richmond.2022}
In preparation of this manuscript, we have learned that the
restriction to $k$-Grassmannian permutations is unnecessary in
\Cref{cor: uAv}.  In fact, Oh and Richmond have proved a
substantially more general statement.  Their results imply that for
any elements $u \leq v$ in any Coxeter group $W$ and any parabolic
subgroup $W_{J}$, there exists a unique maximal element in $[u,v]\cap
uW_{J}$, see \cite[Thm. 2.1]{Oh.Richmond.2022}.
\end{Rem}

\subsection{Intervals and Positroids}
For $I \subseteq [n]$, define $w_0 \cdot I := \{n+1-i : i \in I\}$. For any integer $s$, we also establish the notation
\begin{align*}
I^{+s} := \{i+s \; : \; i \in I\}, \\
I^{-s} := \{i-s \; : \; i \in I\},
\end{align*}
where all values are taken mod $n$ in the range $[n]$.

\begin{Lemma} \label{Lemma: smallest/biggest sets} Let $\w = (w,co)$ be a
decorated permutation in $\SSnk$. Let $\M$ and $[u,v]$ be the
positroid and Grassmann interval corresponding to $\w$. Then every $I
\in \M$ satisfies
\[
I_1(\w) = u[k]  \preceq I  \preceq v[k] = w^{-1}(I_1(\w)).
\]
Furthermore, for any $r \in [n]$, $I$ satisfies
\[
I_r(\w) \preceq_r I \preceq_r w^{-1}(I_r(\w)).
\]
\end{Lemma}

\begin{proof} Recall the equalities $u[k] = I_1(\w)$ and $v[k] =
w^{-1}(I_1(\w))$ from the bijection between decorated permutations and
Grassmann intervals
given in \Cref{Sec: decperms}.  To see the lower bound, use the
characterization of $\M$ from \Cref{th:decoratedpermtopositroid}. Any
$I \in \M$ must satisfy $I_r(\w) \preceq_r I$ for every shifted
anti-exceedance set $I_r(\w)$. In particular, for $r=1$, $I$ must
satisfy $u[k]=I_1(\w) \preceq_1 I$.

To prove the upper bound $I \preceq v[k]$, let $I = \{i_1 < \cdots < i_k\}$, and
suppose that $I \not\preceq v[k] = \{v_1 < \cdots < v_k \}$. Then
there is some maximal index $r$ for which $i_r > v_r$. In particular,
we have $|v[k] \cap [i_r,n]| < k-r+1$. Consider the positroid variety
$\Pi_{\w} = \pi_k(X_u \cap X^v)$. Recall that the opposite Schubert
variety $X^v$ is defined by
\[
X^v = \{V_\bullet \in \mathcal{F}\ell(n) \; : \; \dim(\text{Proj}_{[n-j+1,n]}(V_i)) \leq |v[i] \cap [n-j+1,n]| \;\text{for all $i,j$}\}.
\]
Any point in $\Pi_{\w}$ is the projection of a point in $X^v$. Thus,
for any matrix $A$ representing a point in $\Pi_{\w}$, the rank of
columns $[i_r,n]$ of $A$ must be strictly less than $k-r+1$. In
particular, the rank of $A$ in the column set $\{i_r, \ldots, i_k\}$
is strictly less than $k-r+1$ so $\rank(A_I) < k$, which implies
$\Delta_I(M) = 0$. Since this argument holds for every matrix
representing a point in $\Pi_{\w}$, then by \Cref{Def: positroid as
projection} it follows that $I \notin \M$. Thus, $I \in \M$ must also satisfy $I \preceq v[k]$.

For general $r$, we have already seen above from the definition of $\M(\w)$ that $I$ must satisfy $I_r(\w) \preceq_r I$. Thus again, it remains to show that $I \preceq_r w^{-1}(I_r(\w))$. Consider the action of rotating $D(\w)$ in the counterclockwise direction by $r-1$ positions. As discussed further in \Cref{Sec: rigid transformations}, this transformation of $D(\w)$ has the action of cyclically shifting the anti-exceedance sets of $\w$ as well as the ground set of $\M(\w)$. Let $\z = (z,co')$ be the decorated permutation whose chord diagram is obtained by applying this rotation to $D(\w)$. Then $I_1(\z) = I_r(\w)^{-(r-1)}$, $z^{-1}(I_1(\z)) = w^{-1}(I_r(\w))^{-(r-1)}$, and $\M(\z) = \M(\w)^{-(r-1)}$. Applying the result above yields that $I^{-(r-1)} \preceq_r z^{-1}(I_1(\z)) = w^{-1}(I_r(\w))^{-(r-1)}$. Cyclically shifting the values back up by $r-1$ positions, it follows that $I = (I^{-(r-1)})^{+(r-1)} \preceq_r (w^{-1}(I_r(\w))^{-(r-1)})^{+(r-1)} = w^{-1}(I_r(\w))$, as desired.
\end{proof}

The result of \Cref{Lemma: smallest/biggest sets} gives a constraint on the elements in a positroid. From this constraint, we prove the characterization of a positroid as a collection of initial sets of the corresponding Grassmann interval.

\begin{proof}[Proof of \Cref{Th: positroid = initial sets}]
Let $\SSS
= \{y[k] \; : \; y \in [u,v]\}$ be the collection of initial sets for
$[u,v]$. For $y \in [u,v]$, the point represented by the permutation
matrix $M_y$ is in the Richardson variety $X_u^v$. Let $A_{y}=A_{y[k]}$ be the
submatrix of $M_{y}$ given by the top $k$ rows.  Then $A_{y}$
represents a point in the positroid variety $\pi_k(X_u^v) = \Pi_{[
u,v]}$, by Part (1) of \Cref{Def: positroid as projection}. The restriction of $A_y$ to the columns indexed by $I = y[k]$
is a $k \times k$ permutation matrix. Therefore, $\Delta_I$ does not
vanish at $A_y \in \Pi_{[u,v]}$. Since $\Delta_I$ is not in the vanishing
ideal of $\Pi_{[u,v]}$, we must have  $I \in \M$ by Part (2) of \Cref{Def: positroid as
projection}. Thus, we have the inclusion $\SSS \subseteq \M$.

For the reverse inclusion, let $I \in \M$.  By \Cref{Lemma:
smallest/biggest sets}, $I$ must satisfy $u[k] \preceq I \preceq
v[k]$. In particular, since $I \preceq v[k]$, we may use \Cref{def:
uAv} to define a permutation $y = u(I,v)$ with initial set $y[k] =
I$. By \Cref{cor: uAv}, $y \leq v$ under Bruhat order, so to prove $I
\in \SSS$ it remains only to show that this $y$ satisfies $u \leq y$.

By \cite[Ex. 8, Ch 2]{b-b}, $u \leq y$ if and only if $u[j]
\preceq y[j]$ for all $j \in [k]$ and $u[n-j+1,n] \succeq y[n-j+1,n]$
for all $j \in [n-k]$.  This exercise follows from
\Cref{th:tableau.criterion}.  We will prove that $u[j] \preceq y[j]$ for $j\in [k]$ by
induction on $j$.

Let $\w$ be the decorated permutation associated to $[u,v]$.  Since $I
\in \M$, we know by \Cref{th:decoratedpermtopositroid} that $I_r(\w)
\preceq_{r} I$ for all $r \in [n]$ where $I_{r}(\w)$ is the shifted
anti-exceedance set for $\w$. By the shuffling algorithm described in
\Cref{sec:background} mapping $\w$ to $[u,v]$, one can observe for all $1 \leq j \leq k$ that
$u_{j} \in I_r(\w)$ unless $u_{j} < r\leq v_{j}$.  In particular, $u[j] \subset I_{r}(\w)$ for all $r>v_{j}$
since $v_{1}<\cdots < v_{j}$ and $u \leq v$.  Therefore,
\begin{equation} \label{eq: induction 2}
|u[j] \cap [i,v_j]| \leq  |I_r(\w) \cap [i,v_j]| \; \forall \; r >v_{j}.
\end{equation}

For the case $j=1$, fix $r = v_1+1 $ (or if $v_{1}=n$ fix $r=1$) so
$v_1$ is maximal under the order $<_r$. As observed above, $u_{1} \in
I_r(\w)$.  Since $I_r(\w) \cap [u_1,v_1]$ is nonempty, $I_r(\w)
\preceq_{r} I$ implies $I \cap [u_1,v_1]$ must also be nonempty.  By
construction of $y=u(I,v)$ and \Cref{def: uAv}, $y_1$ is the maximal
element of $I \cap [u_1,v_1]$ .  Thus $u_1 \leq y_1$.

For the inductive step, fix $j \in [2,k-1]$, and assume that 
\begin{equation}\label{eq:as.bs}
u[j-1] =
\{a_{1}< \ldots < a_{j-1} \} \;\preceq  \; y[j-1] = \{b_{1}< \ldots <
b_{j-1} \}.
\end{equation}
We note that all $a_{i},b_{i} \in [v_{j}]$ since $v$ is
$k$-Grassmannian and $u,y\leq v$ in Bruhat order by \Cref{Lemma: BS}.
If $u_{j}\leq y_{j}$, then no matter where $u_{j}$ and $y_{j}$ fit in
among the increasing sequence of $a_{i}$'s and $b_{i}$'s, we will have
$u[j]\preceq y[j]$ by \eqref{eq:as.bs}.  Therefore, assume $1\leq
y_{j}<u_{j}\leq v_{j}$.  By the induction hypothesis, $1\leq a_{i}\leq
b_{i}\leq v_{i}$ for each $1\leq i <j$, so if $m$ is the largest index
such that $b_{m}<y_{j}$, then
\[
u[j] \cap [1,a_{m}] \preceq y[j] \cap [1,b_{m}].  
\]
So, to prove $u[j] \preceq v[j]$, it suffices to show that for all $i
\in [y_{j},v_{j}]$
\begin{equation} \label{eq: induction to show}
|u[j] \cap [i,v_j]| \leq |y[j] \cap [i,v_j]|.
\end{equation}

Fix $r = v_{j}+1$ modulo $n$ so $v_{j}$ is maximal in the shifted
$\leq_{r}$ order.  By definition of $y=u(I,v)$, we know $y_{j}$ is
maximal in $(I \setminus y[j-1]) \cap [v_{j}]$ under $\leq_{r}$ and
the usual order on the integers.  Therefore, for all $i \in [y_{j},
v_{j}]$ we have $ y[j] \cap [i,v_j] = I \cap [i,v_j]$.  Hence, since
$I_r(\w) \preceq_{r} I$, we have
\begin{equation} \label{eq: induction 1}
|I_r(\w) \cap [i,v_j]| \leq |I \cap [i,v_j]| = | y[j] \cap [i,v_j]|.
\end{equation}
By \eqref{eq: induction 2}, 
\begin{equation} \label{eq: induction 3}
|u[j] \cap [i,v_j]| \leq  |I_r(\w) \cap [i,v_j]|, 
\end{equation}
so \eqref{eq: induction 1} and \eqref{eq: induction 3} together imply
\eqref{eq: induction to show} for all $i \in [y_{j},v_{j}]$.  Hence,
$u[j] \preceq y[j]$ for all $j\in[k]$.

Observe that the map on permutations in $S_{n}$ sending $x \to
x'=w_{0}x w_{0}$ rotates the permutation matrices 180 degrees.  Using
\Cref{Lemma: BS} we have $[u,v] \in \GIkn$ if and only if
$[w_{0}uw_{0},w_{0}vw_{0}] \in Gi(n-k,n)$.  By the symmetry in
\Cref{def: uAv}, we observe that $u(y[k+1,n],
w_{0}vw_{0})=w_{0}yw_{0}$.  Therefore, by applying the argument above
to $y'\in [u',v']$ with $I'=w_{0} \cdot y[k+1,n]$, we have $u'[j] \preceq y'[j]$
for all $j \in [n-k]$ which implies $u[n-j+1,n] \succeq y[n-j+1,n]$
for $j \in [n-k]$ as needed to complete the proof.
\end{proof}

\begin{proof}[Proof of \Cref{The: Main theorem}, $ (2) \Leftrightarrow (3) \Leftrightarrow (4)$]
By \Cref{Th: positroid = initial sets}, we have the equality $\M = \{y[k] : y \in [u,v] \}$. By \Cref{thm:Alignments}, we have $\# \aligns(\w) = k(n-k) - [\ell(v) - \ell(u)]$. Hence, $\ell(v) - \ell(u) = k(n-k) - \# \aligns(\w)$. Together, these facts yield (2) $\Leftrightarrow$ (3).

The equality (3) $\Leftrightarrow$ (4) comes from the definition of the matroid Johnson graph $J(\M)$. For a set $J \in \M$, the vertices adjacent to $J$ in $J(\M)$ are exactly the sets $I \in \M$ such that $|I \cap J| = k-1$. So again, the equality $\ell(v) - \ell(u) = k(n-k) - \# \aligns(\w)$ from \Cref{thm:Alignments} implies (3) $\Leftrightarrow$ (4).
\end{proof}

\section{Reduction to $T$-fixed points}
\label{Sec: perm matrices suffice}

To determine whether a variety is singular, in general one would have
to test every point in the variety. In the results below, we show that
for a positroid variety $\Pi_{[ u,v ]}$, one need only check whether
$\Pi_{[u,v]}$ is singular at the $T$-fixed points $A_{y[k]}$ for $y \in [u,v]$. In particular, we will complete
the proofs of \Cref{thm:bounded.below} and
\Cref{thm:tangent.space.dim}. Together these two theorems along with \Cref{Th: positroid = initial sets} prove the
equality between (1) and (2) in \Cref{The: Main theorem} as mentioned
in the introduction.

First, we recall some specifics for the definition of the Jacobian
matrix of a positroid variety.  By \Cref{Def: positroid as projection}
Part (1), the polynomials generating the vanishing ideal of
$\Pi_{[u,v]}$ are exactly the determinants $\Delta_I$ for $I \in
\QQ_{[u,v]}$ where
\begin{equation}\label{eq:nolabel}
\QQ_{[u,v]} \; := \; \binom{[n]}{k} \setminus \M \; = \; \big\{I \in \binom{[n]}{k} : I \neq y[k] \; \forall y\in [u,v]
\big\}
\end{equation}
is the set of nonbases of the positroid $\M$ corresponding to $[u,v] \in
\GIkn$. The second equality in \eqref{eq:nolabel} holds by \Cref{Th: positroid = initial
sets}.  Each determinant $\Delta_I$ is a polynomial function using variables of
the form $x_{ij}$ indexed by row $i \in [k]$ and column $j\in [n]$.
Let $Jac_{[u,v]}$, or just $Jac$ if $[u,v]$ is understood, be the
Jacobian of $\Pi_{[u,v]}$.  Similarly, we will suppress the subscript
$\QQ=\QQ_{[u,v]}$. Then, the rows of $Jac$ are indexed by the sets $I
\in \QQ$, the columns are indexed by variables $x_{ij}$, and the $(I,
x_{ij})$ entry of $Jac$ is $\frac{\partial \Delta_I}{\partial
x_{ij}}.$ If $A$ is a $k\times n$ complex matrix representing a point
in $\Pi_{[u,v]} \subseteq \Gkn$, then $Jac|_{A}$ is the matrix with
entries in the complex numbers obtained from $Jac$ by evaluating each
entry at $A$,
\begin{equation}\label{eq:jac.A}
Jac|_{A}(I,x_{ij}) = \frac{\partial \Delta_I}{\partial x_{ij}}(A).
\end{equation}

By \Cref{Def: positroid as projection} Part (2), we can write $\Pi_{[ u,v ]}$ as the union
\begin{align} \label{eq: decomp of pos var by Rich}
\Pi_{[ u,v ]} \; = \; \pi_k (X_u^v) = \pi_k \Big(\bigsqcup_{u \leq y \leq v}(C_{y} \cap X^v) \Big) \;  = \bigcup_{u \leq y \leq v} \pi_k(C_{y} \cap X^v).
\end{align}
For any $y \in [u,v]$, recall $A_{y[k]}$ is the projection of the
permutation matrix $M_{y}$ to the top $k$ rows.

\begin{Lemma} \label{lem: Jac a partial perm matrix} For $[u,v] \in
\GIkn$ and $y \in [u,v]$, the matrix $Jac_{[u,v]}|_{A_{y[k]}}$ is a
partial permutation matrix, up to the signs of the
entries. Furthermore, the nonzero entries occur exactly in the entries
$(I,x_{st})$, where $I \in \QQ_{[u,v]}$ satisfies $|I \cap y[k]| =
k-1$, the value $s \in [k]$ is determined by $y[k] \setminus
I = \{y_{s}\}$, and $t \in [n]$ is
determined by $I \setminus y[k] = \{t\}$.
\end{Lemma}

\begin{proof} 

Write $A_{y[k]}|_I$ for the restriction of $A_{y[k]}$ to column set $I$. By expanding $\Delta_I$ along row $i \in [k]$ of the partial
permutation matrix $A_{y[k]}$, observe from \eqref{eq:jac.A} that
entry $(I, x_{ij})$ of $Jac|_{A_{y[k]}}$ is
\vspace*{-0.1in}
\begin{enumerate}
\item[(i)] ($j \in I$) up to sign, the $(k-1) \times (k-1)$ minor of $A_{y[k]}|_{I}$ in column set $I \setminus \{j\}$ and row set $[k] \setminus \{i\}$, or
\item[(ii)] ($j \notin I$) 0, since $\Delta_I$ does not contain $x_{ij}$.
\end{enumerate}
\vspace*{-0.1in}

Consider any $I \in \QQ$. By \Cref{Th: positroid = initial sets}, $y[k] \in \M$, so $I \neq y[k]$. Since $I$ and $y[k]$ are distinct $k$-sets, then $|I \cap y[k]|$ is at most $k-1$. If $|I \cap y[k]| \leq k-2$, then $A_{y[k]}$ has at most $k-2$ ones in $A_{y[k]}|_I$, and the remaining entries of $A_{y[k]}|_I$ are zeros. In that case, every $(k-1) \times (k-1)$ minor of $A_{y[k]}|_I$ is zero. Therefore, row $I$ of $Jac|_{A_{y[k]}}$ contains only zeros.  

Otherwise, $I \in \QQ$ satisfies $|I \cap y[k]| = k-1$.  Say $I \setminus
y[k] = \{t\}$ and $y[k] \setminus I = \{y_{s}\}$ for some $s\in
[k]$. Then column $t$ of $A_{y[k]}$ contains all zeros, row $s$ of
$A_{y[k]}|_I$ contains all zeros, and the submatrix of $A_{y[k]}|_I$
obtained by removing row $s$ and column $t$ is a permutation matrix.
Therefore, by cofactor expansion of the determinant 
\[
Jac\Big|_{A_{y[k]}}(I,x_{st})\; = \; \frac{\partial
\Delta_I}{\partial x_{st}} (A_{y[k]})\; = \; \pm 1 ,
\]
and $x_{st}$
is the unique variable such that $\frac{\partial \Delta_I}{\partial
x_{ij}} (A_{y[k]})$ is nonzero. 

Conversely, we claim column $x_{st}$ contains a unique nonzero entry
in row $I$.  If $I' \in \QQ$ is another set such that $|I' \cap y[k]|
= k-1$, let $I' \setminus y[k] = \{t'\}$ and $y[k] \setminus I' =
\{y_{s'}\}$.  Then, the only nonzero entry of $Jac|_{A_{y[k]}}$ in row
$I'$ occurs in column $x_{s't'}$ by the same reasoning as above. If $t \neq t'$,
then $x_{st} \neq x_{s't'}$. If $t=t'$, either $s=s'$ and $I'=I$ or
else $y_{s} \in I'$ so $s' \neq s$ since $y$ is a bijection, in which
case $x_{st} \neq x_{s't'}$.  Therefore, $Jac|_{A_{y[k]}}$ is a
partial permutation matrix up to signs as stated.
\end{proof}

\begin{proof}[Proof of \Cref{thm:tangent.space.dim}] From \Cref{lem:
Jac a partial perm matrix}, $Jac|_{A_{y[k]}}$ is a partial permutation
matrix, up to the signs of its entries. Thus, the rank of
$Jac|_{A_{y[k]}}$ is exactly the number of its nonzero entries, which
are exactly in the rows indexed by the sets $I \in \QQ$ satisfying $|I \cap y[k]|
= k-1$. Therefore, the rank of $Jac|_{A_{y[k]}}$ is equal to the
number of these sets,
\[
\rank (Jac|_{A_{y[k]}}) \; = \; \# \Big\{I \in \binom{[n]}{k} \setminus \M \; : \; |I \cap y[k]| = k-1
\Big\}.
\]
\end{proof}

Next, consider a flag $V_{\bullet} \in C_{y} \cap X^v$ and its
canonical matrix representative $A_{V_\bullet}$ as defined in
\Cref{sec:background}.  Let $A$ be the projection of $A_{V_\bullet}$ to $\Gkn$. We will see that the rank of the Jacobian matrix evaluated at $A$ is at least the rank of the Jacobian evaluated at $A_{y[k]}$. \\


\begin{proof}[Proof of \Cref{thm:bounded.below}] Without loss of
generality, we can assume $A \in \pi_k(C_{y} \cap X^v)$ is a canonical
matrix with leading ones in exactly the same entries as $A_{y[k]}$.  Therefore, 
permuting the columns of $A$ and $A_{y[k]}$ in the same way, we can
assume the first $k$ columns of both have the form of an upper
triangular matrix with ones on the diagonal. If we apply the same
permutation to the values in the $k$-sets in $\QQ$ and the variables $x_{ij}$, then such
a rearrangement of columns will not change the rank of
$Jac|_{A_{y[k]}}$ or $Jac|_A$.  After permuting, we can then
assume $y$ is the identity permutation, so $A_{y[k]} = A_{[k]}$.  Let
$N = \rank(Jac|_{A_{[k]}})$.  By \Cref{thm:tangent.space.dim}, we must
show that
\[
\rank(Jac|_A) \; \geq \; N = \; \# \{I \in \QQ \; : \; |I \cap [k]| = k-1 \}.
\]

Let $R = \{ I_1, \ldots, I_N\}$ and $C =\{(s_1,t_1), \ldots, (s_N,
t_N)\}$ be the unique row and column sets such that the $N \times N$
minor of $Jac|_{A_{[k]}}$ in rows $R$ and columns $C$ is nonzero, as
determined in \Cref{lem: Jac a partial perm matrix}.  Then, for each $
I_j \in R$, we know $I_{j} = ([k] \setminus \{s_{j}\}) \cup
\{t_{j}\}$, where $s \in [k]$ and $t_{j} > k$ after permuting values.
Since the first $k$ columns of $A$ form an upper triangular matrix
with ones along the diagonal, we observe that for all $j \in [N]$,
\begin{equation}\label{eq:jac.entries.R.C}
Jac|_A(I_{j}, x_{s_{j}t_{j}}) = \frac{\partial
\Delta_{I_{j}}}{\partial x_{s_{j}t_{j}}}(A) = \pm 1.
\end{equation}

Fix $j \in [N]$. Assume that the sets in $R$ are ordered lexicographically as sorted lists, and let $1 \leq h < j \leq N$. Either $t_h \neq t_j$ or $t_h=t_j$ and $s_h > s_j$ by lex order. If $t_{h} \neq t_{j}$, then
$I_{h}$ does not contain $t_j$, by construction of $I_{h}$. Hence, $\Delta_{I_{h}}$ does not depend on the variable $x_{s_{j}t_{j}}$, so $\frac{\partial
\Delta_{I_{h}}}{\partial x_{s_{j}t_{j}}}(A) = 0$.  Otherwise, $s_{h} > s_{j}$, and $\frac{\partial \Delta_{I_{h}}}{\partial
x_{s_{j}t_{j}}}(A)$ is the determinant of an upper triangular matrix
with a zero on its diagonal.  Hence, in either case
\begin{equation}\label{eq:jac.R.S.2}
Jac|_A(I_{h},x_{s_{j}t_{j}}) = \frac{\partial \Delta_{I_{h}}}{\partial
x_{s_{j}t_{j}}}(A) = 0.
\end{equation}

Arrange the rows and columns of $Jac$ so that $I_1, \ldots, I_N$ are
the top $N$ rows listed in order and $x_{s_{1}t_{1}},
x_{s_{2}t_{2}}, \ldots, x_{s_{N}kt_{N}}$ are the first $N$ columns
listed in order.  Then, the $N \times N$ upper left submatrix of
$Jac|_A$ is lower triangular with plus or minus one entries along the diagonal.  Thus
$Jac|_A$ contains a rank $N$ submatrix, so $\rank(Jac|_A) \geq N$.
\end{proof}

\begin{Ex} \label{ex: lower triangular}
Rearranging the columns of the matrix in \Cref{ex:matroid} so that $y[k] = [k]$, consider the matrix
\[
A = \begin{bmatrix}
1 & 2 & 3 & 0 & 4 & 0 \\
0 & 1 & 0 & 0 & 2 & 1
\end{bmatrix} .
\]
The nonbases of $\M_A$ are
\[
\QQ = \{ \{1,3\}, \{1,4\}, \{2,4\}, \{2,5\}, \{3,4\}, \{4,5\}, \{4,6\}\} \subseteq \binom{[6]}{2}.
\]
With $k = 2$, the sets $I \in \QQ$ satisfying $|[k] \cap I| = k-1$, ordered lexicographically, are $I_1 = \{1,3\}$, $I_2 = \{1,4\}$, $I_3 = \{2,4\}$, and $I_4 = \{2,5\}$. The corresponding $x_{st}$ so that $I_j = ([k] \setminus \{s_j\}) \cup \{t_j\}$ are $x_{23}$, $x_{24}$, $x_{14}$, and $x_{15}$, respectively. Arranging the rows and columns of $Jac|_A$ as described in the proof above, the upper left submatrix of $Jac|_A$ is a lower triangular matrix with plus or minus one entries on the diagonal, as shown below.

\begin{equation*}
  \begin{blockarray}{r *{4}{c}}
    \begin{block}{r *{4}{>{$\footnotesize}c<{$}} }
      & $x_{23}$ & $x_{24}$ & $x_{14}$ & $x_{15}$  \\
    \end{block}
    \begin{block}{>{$\footnotesize}r<{$}[*{4}{c}]}
   \bigstrut[t]  $\{1,3\}$ &  1 & 0 & 0 & 0   \\
     $\{1,4\}$ & 0 & 1 & 0 & 0  \\
     $\{2,4\}$ & 0 & 2 & -1 & 0  \\
     $\{2,5\}$ & 0 & 0 & 0 & 1  \\
    \end{block}
  \end{blockarray} 
\end{equation*}

\end{Ex}

\begin{Cor} \label{Cor: singular at piece of partition} Let $[u,v] \in \GIkn$ and $y \in [u,v]$. Then $\pi_k(C_y \cap X^v) \subseteq \Pi_{[u,v]}$ contains a singular point of $\Pi_{[u,v]}$ if and only if $\Pi_{[u,v]}$ is singular at $A_{y[k]}$. Furthermore, $A_{y[k]}$ is a singularity of $\Pi_{[u,v]}$ if and only if
\[
\# \{I \in \QQ_{[u,v]} \; ; \; |I \cap y[k]|=k-1\} \; < \; k(n-k) - [\ell(v)-\ell(u)].
\]
\end{Cor}

\begin{proof}
Recall from \eqref{eq: sing condition Jac algn} that $A \in \Pi_{[u,v]}$ is a singular point of $\Pi_{[u,v]}$ if and only if $\rank(Jac|_A) < \text{codim} \, \Pi_{[u,v]}$. By \Cref{thm:bounded.below}, for any $A \in \pi_k(C_y \cap X^v)$, $\rank(Jac|_{A_{y[k]}}) \leq \rank(Jac|_A)$. Therefore, if there exists some $A \in \pi_k(C_y \cap X^v)$ that is a singular point of $\Pi_{[u,v]}$, then
\[
\rank(Jac|_{A_{y[k]}}) \; \leq \; \rank(Jac|_A) \; < \; \text{codim} \, \Pi_{[u,v]}
\]
implies that $A_{y[k]}$ is also a singular point of $\Pi_{[u,v]}$.  Conversely, note that since the permutation matrix $M_y$ is in $C_y \cap X^v$, then $A_{y[k]} = \pi_k(M_y) \in \pi_k(C_y \cap X^v)$. Thus, if $A_{y[k]}$ is a singular point of $\Pi_{[u,v]}$, then $A_{y[k]}$ is already a singular point of $\Pi_{[u,v]}$ in $\pi_k(C_y \cap X^v)$. This proves the first statement.

For the second statement, recall from \Cref{thm:Alignments} that
\[
\text{codim} \ \Pi_{[u,v]} = k(n-k) - [\ell(v) - \ell(u)].
\]
By \Cref{thm:tangent.space.dim}, $\rank(Jac|_{A_{y[k]}}) = \# \{I \in \QQ \, : \, |I \cap y[k]|=k-1\}$. Therefore, $A_{y[k]}$ is a singular point of $\Pi_{[u,v]}$ if and only if
\[
\# \{I \in \QQ  :  |I \cap y[k]|=k-1\} = \rank(Jac|_{A_{y[k]}}) < \text{codim} \, \Pi_{[u,v]} = k(n-k) - [\ell(v)-\ell(u)].
\]
\end{proof}

\begin{proof}[Proof of \Cref{The: Main theorem}, $(1) \Leftrightarrow
(2)$] From the definition of a smooth variety and the decomposition 
\[
\Pi_{[ u,v ]} \; = \; \pi_k (X_u^v) = \pi_k \Big(\bigsqcup_{u \leq y
\leq v}(C_{y} \cap X^v) \Big) \;  = \bigcup_{u \leq y \leq v}
\pi_k(C_{y} \cap X^v), 
\]
it follows that $\Pi_{[u,v]}$ is smooth if and only if, for every $y \in [u,v]$, every point in $\pi_k(C_y \cap X^v)$ is a smooth point of $\Pi_{[u,v]}$. By \Cref{Cor: singular at piece of partition}, every point of $\pi_k(C_y \cap X^v)$ is a smooth point of $\Pi_{[u,v]}$ if and only if $A_{y[k]}$ is a smooth point of $\Pi_{[u,v]}$, which occurs if and only if
\begin{equation} \label{eq : smooth with Q}
\# \{I \in \QQ \, : \, |I \cap y[k]|=k-1\} \; \geq \; k(n-k) - [\ell(v)-\ell(u)].
\end{equation}
Since $\rank(Jac|_{A_{y[k]}})$ is the codimension of the tangent space to $\Pi_{[u,v]}$ at $A_{y[k]}$, then $\rank(Jac|_{A_{y[k]}})$ is bounded above by $\text{codim} \, \Pi_{[u,v]}$. Thus, the inequality in \eqref{eq : smooth with Q} can never be strict. Therefore, $A_{y[k]}$ is a smooth point of $\Pi_{[u,v]}$ if and only if
\begin{equation} \label{eq : smooth with Q 2}
 \# \{I \in \QQ \, : \, |I \cap y[k]|=k-1\} = \rank(Jac|_{A_{y[k]}})  =  \text{codim} \, \Pi_{[u,v]}  =  k(n-k) - [\ell(v)-\ell(u)].
\end{equation}

We can compute the left side of \eqref{eq : smooth with Q 2} as follows. For any
$j \in y[k]$ and any element $ i\in [n] \backslash y[k]$, observe that
there are $k(n-k)$ sets $I \in \binom{[n]}{k}$ of the form $I = (y[k]
\backslash \{j\}) \cup \{i\}$ so that $|I \cap y[k]| = k-1$. Since
$\QQ = \binom{[n]}{k} \setminus \M$, then
\begin{equation} \label{eq: Q + M}
\# \{I \in \M \, : \, |I \cap y[k]|=k-1\} \; + \; \# \{I \in \QQ \, : \, |I \cap y[k]|=k-1\} \; = \; k(n-k).
\end{equation}
Substituting the result of \eqref{eq: Q + M} into \eqref{eq : smooth with Q 2} yields that $A_{y[k]}$ is a smooth point of $\Pi_{[u,v]}$ if and only if
\begin{equation} \label{eq : smooth with M}
 \# \{I \in \M \, : \, |I \cap y[k]|=k-1\} \; =  \;\ell(v)-\ell(u).
\end{equation}
Thus, $\Pi_{[u,v]}$ contains only smooth points if and only if every $y \in [u,v]$ satisfies \eqref{eq : smooth with M}.
\end{proof}

Finally in this section, we relate the vertex degree conditions in the
Johnson graph to the sets of nonbases of positroids similar to \Cref{Cor:
singular at piece of partition}.    If a positroid $\M$ corresponds with
the decorated permutation $\w$, let $\aligns(\M)=\aligns(\w)$.

\begin{Cor}\label{cor:anti-exchangable.pairs.alignments}
Given a rank $k$ positroid $\M$ on ground set $[n]$, let $\QQ: = \binom{[n] }{k} \setminus \M$
be the corresponding set of nonbases.  For any $J \in \M$,
the codimension of the tangent space to $\Pi_{\M}$ at $A_J$ is
\begin{equation} \label{eq: num anti-exch pairs less equal codim}
\# \{I \in \QQ \; : \; |I \cap J| = k-1 \} \; \leq \;  \#
\aligns(\M) = \mathrm{codim}(\Pi_{\M}).
\end{equation}
Furthermore, $A_{J}$ is a singular point in $\Pi_{\M}$ if and only if
\begin{equation} \label{eq: num anti-exch pairs strictly less codim}
\# \{I \in \QQ \; : \; |I \cap J| = k-1 \} \; < \;  \# \aligns(\M).
\end{equation}
\end{Cor}

\begin{proof}
Equation \eqref{eq: num anti-exch pairs less equal codim} follows from
Equation \eqref{eq: sing condition Jac algn}, \Cref{Th: positroid =
initial sets}, and \Cref{thm:tangent.space.dim}.  The second claim now
follows by \Cref{Cor: singular at piece of partition} and \Cref{Th: positroid =
initial sets}.
\end{proof}

\section{Rigid Transformations} \label{Sec: rigid transformations} The
authors of \cite{Ardila-Rincon-Williams} show that the set of
positroids is closed under restriction, contraction, duality, and a
cyclic shift of the ground set. We add to this list reversal of the
ground set.  The fact that positroids are closed under duality, cyclic
shift of the ground set, and reversal of the ground set can be
obtained by considering rigid transformations of the chord diagram
of the associated decorated permutation. We consider
three types of rigid transformations on chord diagrams: arc reversal, reflection, and rotation. We associate these transformations with
the symmetric group operations of taking the inverse, conjugation by $w_0$, and
conjugation by a cycle. Because these rigid transformations are bijections on the set of chord diagrams, they generate a group of transformations on Grassmann intervals, Grassmann necklaces, and positroids. The results are collected in \Cref{prop: rigid
transformations}.

Throughout this section, fix $\w  = (w,co) \in \SSnk$, and let $\z = (z,co')$ be the decorated permutation whose chord diagram $D(\z)$ is obtained from $D(\w)$ by a rigid transformation. Denote by $F(\w)$ the set of fixed points of $\w$. Define the map $flip : \{\cw, \ccw\} \rightarrow \{\ccw, \cw\}$ to be the involution on $\{\cw, \ccw\}$ that reverses the orientation. Let $\chi \in S_n$ be the cycle with $\chi(i) = i+1$ mod $n$.

To give the maps on decorated permutations corresponding to the chord diagram transformations, we will consider the transformations on the arcs of the chord diagrams. The transformations of arcs then lead to transformations of the two-line notation for decorated permutations.

First, consider \textit{arc reversal}, where where $D(\z)$ is obtained from $D(\w)$ by reversing all the arcs in $D(\w)$. So, an arc $i \mapsto w(i)$ becomes $w(i) \mapsto i$. When $i \in F(\w)$, then $i$ becomes a fixed point in $\z$ with opposite orientation of $i$. Then $z = w^{-1}$, and $F(\z) = F(\w)$ with all fixed point orientations reversed. Therefore, $co' = flip \circ co$.  With these observations in mind, we define 
\begin{equation} \label{eq: arc reversal}
(\w)^{-1} := (w^{-1}, flip \circ co). 
\end{equation}
The two-line notation for $(\w)^{-1}$ is obtained from the two-line notation of $\w$ by swapping the rows, reversing arrows labeling fixed points, and reordering the columns so that the entries of the top line appear in increasing order.

Next, consider \textit{reflection}, where $D(\z)$ is obtained from $D(\w)$ by reflecting $D(\w)$ across the vertical axis. So, an arc $i \mapsto w(i)$ becomes $w_0(i) \mapsto w_0(w(i))$. When $i \in F(\w)$, then $w_0(i)$ becomes a fixed point of $\z$ with opposite orientation of $i$. Then $z = w_0 w w_0$, and $F(\z) = w_0 \cdot F(\w)$ with all fixed point orientations reversed after applying $w_0$ to the value of a fixed point in $F(\z)$. Therefore, $co' = flip \circ co \circ w_0$, and we define
\begin{equation} \label{eq: decperm conj w0}
w_0 \cdot \w := (w_0 w w_0, flip \circ co \circ w_0). 
\end{equation}
The two-line notation for $\z$ is obtained from the two-line notation of $\w$ by replacing $i$ with $w_0(i)$ in both lines, reversing all arrows labeling fixed points, and reversing the order of the columns.

Finally, consider \textit{rotation}, where $D(\z)$ is obtained from $D(\w)$ by rotating $D(\w)$ by $s$ units in the clockwise direction. So, an arc $i \mapsto w(i)$ becomes $i+s \mapsto w(i) + s$, taken modulo $n$. When $i \in F(\w)$, then $i+s$ becomes a fixed point of $\z$ with the same orientation as $i$. Then $z = \chi^s w \chi^{-s}$, and $F(\z) = F(\w)^{+s}$ with all fixed point labels preserved after applying $\chi^{-s}$ to the value of a fixed point in $F(\z)$. Therefore, $co' = co \circ \chi^{-s}$, and we define
\begin{equation} \label{eq: decperm rotation}
\chi^{s} \cdot \w := (\chi^s w \chi^{-s}, co \circ \chi^{-s}).
\end{equation}
The two-line notation for $\z$ is obtained from the two-line notation of $\w$ by replacing $i$ with $i+s$ mod $n$ in both lines and cyclically shifting all columns $s$ units to the right.

Associated with the maps listed above are maps on Grassmann intervals. Let $[u,v] \in \GIkn$ be the Grassmann interval corresponding to $\w$, and let $[u',v']$ be the Grassmann interval corresponding to $\z$. As we know from the shuffling algorithm described in \Cref{Sec: decperms}, $u'$ and $v'$ are easily extracted from the two-line notation for $\z$ by reordering the columns so that the highlighted columns, corresponding to the anti-exceedances of $\z$, appear on the left and the top line has the form of a $k(\z)$-Grassmannian permutation. This new array is the two-line array $\uprimevprime$. 

We have already described how the two-line notation for $\z$ is obtained from the two-line notation of $\w$ under the three rigid transformations of $D(\w)$. From these descriptions, we obtain the following maps $\uv \mapsto \uprimevprime$.
\begin{enumerate}
\item[(1)] For the arc reversal map, $\uprimevprime$ is obtained from $\uv$ by swapping the rows, swapping the left $k$-column block with the right $(n-k)$-column block, and reordering the columns within the two blocks so that the top row gives an $(n-k)$-Grassmannian permutation. Define $[u,v]^{-1} \in Gi(n-k,n)$ to be the Grassmann interval whose two-line array is obtained in this way from $\uv$.
\item[(2)] For the reflection map, $\uprimevprime$ is obtained from $\uv$ by replacing every $i$ with $w_0(i)$ in both rows and reversing all columns. It follows that $u' = w_0 u w_0$ and $v' = w_0 v w_0$. Define $w_0 \cdot [u,v] := [w_0uw_0, w_0vw_0] \in Gi(n-k,n)$. Note that this operation on $[u,v]$ of conjugation by $w_0$ was used at the end of the proof of \Cref{Th: positroid = initial sets}.
\item[(3)] For the rotation map, $\uprimevprime$ is obtained from $\uv$ by replacing every $i$ with $i+s$, highlighting all columns corresponding to anti-exceedances, then reordering the columns into a highlighted and an unhighlighted block so that the elements of the top row are increasing within each block. Define $\chi^{s} \cdot [u,v] \in \GIkn$ to be the Grassmann interval obtained from $[u,v]$ in this way.
\end{enumerate}

There are also naturally associated maps on the positroids corresponding to $\w$ and $\z$. These maps on positroids, listed in \Cref{prop: rigid transformations}, can be seen from the following facts.
\begin{enumerate}
\item[(1)] The dual map on positroids corresponds to arc reversal, as shown in the following lemma.
\item[(2)] The map corresponding to $\w \mapsto w_0\w w_0$ on Grassmann intervals is $[u,v] \mapsto [w_0uw_0, w_0vw_0] \in Gi(n-k,n)$, and $[w_0uw_0, w_0vw_0] = \{w_0yw_0 : y \in [u,v]\}$. In particular, the map $y \mapsto w_0 y w_0$ is an interval isomorphism between $[u,v]$ and $[w_0uw_0, w_0vw_0]$. The related map on positroids then follows from \Cref{Th: positroid = initial sets}.
\item[(3)] Every $r$-anti-exceedance of $\w$ translates to an $(r+s)$-anti-exceedance of $\z$ so that $I_r(\w)^{+s} = I_{r+s}(\z)$. Also, for sets $I,J \in \binom{[n]}{k}$, $I \preceq_r J$ if and only if $I^{+s} \preceq_{r+s} J^{+s}$.
\end{enumerate}

\begin{Lemma} \label{Lemma: inv = dual}
Let $\w \in \SSnk$ have associated positroid $\M(\w) \subseteq \binom{[n]}{k}$, and let $\z = (\w)^{-1} \in \SSnminusk$. The positroid associated with $\z$ is the dual of $\M(\w)$,
\begin{equation} \label{eq: inv = dual}
\M(\z) \; = \; \{[n] \setminus I \, : \, I \in \M(\w) \} \; \subseteq \; \binom{[n]}{n-k}.
\end{equation}
\end{Lemma}

\begin{proof}
Since arc reversal is an involution, to show that $\M(\z)$ is the dual of $\M(\w)$, as in \eqref{eq: inv = dual}, it suffices to show that $I \in \M(\w)$ implies that $[n] \setminus I \in \M(\z)$. In particular, we will show that $I \in \M(\w)$ implies that $I_r(\z) \preceq_r [n] \setminus I$ for all $r \in [n]$.

For $r \in [n]$, let $J_r(\w)$ be the set of $r$-exceedances of $\w$, $J_r(\w) = [n] \setminus I_r(\w)$. Fix $r \in [n]$ and $I \in \M(\w)$. By definition, $D(\z)$ is obtained from $D(\w)$ by reversing every arc. Then every arc $a \mapsto w(a)$ corresponding an $r$-anti-exceedance $w(a) \in I_r(\w)$ yields an arc $w(a) \mapsto a$ in $D(\z)$, which corresponds to an $r$-exceedance of $\z$, so that $a \in J_r(\z)$. Similarly, $r$-exceedances $w(a) \in J_r(\w)$ yield $r$-anti-exceedances $a \in I_r(\z)$. Therefore, we have
\begin{align*}
&J_r(\z) = w^{-1}(I_r(\w)) \\
&I_r(\z) = w^{-1}(J_r(\w)) = [n] \setminus w^{-1}(I_r(\w)).
\end{align*}


Since $I \in \M$, then by \Cref{Lemma: smallest/biggest sets}, $I$ must satisfy $I_r(\w) \preceq_r I \preceq_r w^{-1}(I_r(\w))$. In particular, $I \preceq_r w^{-1}(I_r(\w))$ implies that $[n] \setminus I \succeq_r [n] \setminus w^{-1}(I_r(\w)) = I_r(\z)$, as desired. Therefore, $[n] \setminus I$ is in $\M(\z)$ by \Cref{th:decoratedpermtopositroid}.
\end{proof}

The transformations of arc reversal, reflection, and rotation can now be applied to any of the bijectively equivalent objects. These transformations are summarized in the following proposition. 

\begin{Prop} \label{prop: rigid transformations}

Fix a decorated permutation $\w = (w,co) \in \SSnk$ with associated Grassmann interval $[u,v]$, Grassmann necklace $(I_1(\w), \ldots, I_n(\w))$, and positroid $\M(\w)$. Let $\z$ be a decorated permutation whose chord diagram is obtained from $D(\w)$ by (1) arc reversal, (2) reflection, or (3) rotation. Let $[u',v']$ be the Grassmann interval associated with $\z$. The table below describes $\z$ and its associated objects and values.


\begin{table}[ht]
\centering
\begin{tabular}{c||c|c|c|c|c} 
Trans $\backslash$ Obj & $\z$ & $[u',v']$ & $I_r(\z)$ & $\M(\w) \mapsto \M(\z)$ & $k(\z)$ \\
\hline
\hline
arc reversal& $(\w)^{-1}$ & $[u,v]^{-1}$ & $[n] \setminus w^{-1}(I_r(\w))$ & $I \mapsto [n] \setminus I$ & $n-k$ \\
\hline
reflection & $ w_0 \cdot \w $ & $w_0 \cdot [u,v]$ & $w_0 \cdot ([n] \setminus I_{w_0(r)}(\w))$ & $I \mapsto w_0 \cdot ([n] \setminus I)$ & $n-k$ \\
\hline
rotation & $\chi^{s} \cdot \w $ & $\chi^{s} \cdot [u,v]$ & $I_{r-s}(\w)^{+s}$ & $I \mapsto I^{+s}$ & $k$ 
\end{tabular}
\end{table}

\end{Prop}

\begin{Cor} \label{cor: ground set reversal}
The set of positroids is closed under reversal of the ground set $[n]$.
\end{Cor}

\begin{proof}
Let $\M \subseteq \binom{[n]}{k}$ be a positroid, and let $\w$ be the decorated permutation associated with $\M$ so that $\M = \M(\w)$. Let $w_0 \cdot \M = \{w_0 \cdot I \, : \, I \in \M\}$ be the matroid obtained from $\M$ by reversal of the ground set $[n]$. We show that $w_0 \cdot \M$ is a positroid.

Let $\z$ be the decorated permutation whose chord diagram is obtained
from $D(\w)$ by arc reversal followed by a reflection across the
vertical axis. By \Cref{prop: rigid transformations}, $\z = w_0  \cdot
(\w)^{-1}$, and in particular, $\M(\z)$ is obtained from $\M(\w)$ via
the sequence of maps corresponding with arc reversal followed by
reflection on the positroid,
\[
I \quad \mapsto \quad [n] \setminus I \quad \mapsto \quad w_0 \cdot ([n] \setminus ([n] \setminus I)) \, = \, w_0 \cdot I.
\]
Therefore, $\M(\z)$ is exactly $w_0 \cdot \M(\w)$, so $w_0 \cdot \M(\w) = w_0 \cdot \M$ is a positroid by \Cref{th:decoratedpermtopositroid}.
\end{proof}

\begin{Rem}
The closure of the set of positroids under reversal of the ground set can also be obtained using the fact that positroids are the matroids of totally nonnegative matrices. In particular, suppose $A$ is a totally nonnegative $k \times n$ matrix with positroid $\M(A)$. Let $A'$ be obtained by reversing the columns of $A$ and multiplying the bottom row by $(-1)^{\binom{k}{2}}$. Then $A'$ is totally nonnegative, and $\M(A') = \{w_0 \cdot I : I \in \M(A)\}$.
\end{Rem}

\begin{Rem}\label{rem:connectbacktoPost}
Using the notation of rigid transformations, Postnikov's bijection
from decorated permutations to Grassmann intervals mentioned in
\Cref{rem:apost.uv.bij} maps $\w$ to $w_0 \cdot [u,v]^{-1}$, which
preserves the size of the anti-exceedance set.  This map is similar
to the reversal of ground set involution on positroids. 
\end{Rem}

Observe that each of the three transformations in \Cref{prop: rigid transformations} induces a bijection from $\aligns(\w)$ to $\aligns(\z)$. Then $\# \aligns(\w) = \# \aligns(\z)$, so it follows from \Cref{thm:Alignments} that $\text{codim} \, \Pi_{\w} = \text{codim} \, \Pi_{\z}$. Furthermore, as we have seen from the results of \Cref{Sec: perm matrices suffice}, one may determine whether a positroid variety $\Pi_{\M}$, corresponding to a positroid $\M$, is smooth or singular by performing certain computations involving the sets in $\M$. This fact and the maps given in \Cref{prop: rigid transformations} imply the following relationship between $\Pi_{\w}$ and $\Pi_{\z}$.

\begin{Lemma} \label{Lemma: invariances} Let $w \in \SSnk$, and let $\z$ be a decorated permutation whose chord diagram is obtained from $D(\w)$ by (1) arc reversal, (2) reflection, or (3) rotation. Then, in any of these three cases, $\Pi_{\w}$ is smooth if and only if $\Pi_{\z}$ is smooth.
\end{Lemma}

\begin{proof} 
We know from the equivalence of Parts (1) and (3) of \Cref{The: Main theorem} that $\Pi_{\w} \subseteq \Gkn$ is smooth if and only if, for every $J \in \M(\w)$, $J$ satisfies 
\begin{equation} \label{eq: inv smooth condition w}
\# \{I \in \M(\w) \, : \, |I \cap J| = k-1\} = k(n-k) - \# \aligns(\w).
\end{equation}
As noted above, in all three cases, $\# \aligns(\z) = \# \aligns(\w)$. Observe that $\Pi_{\z} \subseteq Gr(k(\z),n)$, where $k(\z) = n-k$ in the cases of arc reversal and reflection, and $k(\z) = k$ in the case of rotation. Therefore, $\Pi_{\z}$ is smooth if and only if every $J \in \M(\z)$ satisfies 
\begin{equation} \label{eq: inv smooth condition z}
\# \{I \in \M(\z) \, : \, |I \cap J| = k(\z)-1\} = k(n-k) - \# \aligns(\w).
\end{equation}

For the cases of arc reversal and reflection, recall from \Cref{prop: rigid transformations} that the positroids $\M(\z)$ and their nonbases in these cases are obtained by the maps $I \mapsto [n] \setminus I$ and $I \mapsto w_0 \cdot ([n] \setminus I)$, respectively. Let $I, J \in \binom{[n]}{k}$. Then $|I \cap J| = k-1$ if and only if $|([n] \setminus I) \cap ([n] \setminus J)| = (n-k)-1$ if and only if $|w_0 \cdot ([n] \setminus I) \cap w_0 \cdot ([n] \setminus J)| = (n-k)-1$. Hence, $I \in \M(\w)$ contributes to the set in \eqref{eq: inv smooth condition w} if and only if $[n] \setminus I \in \M((\w)^{-1})$ contributes to the set in \eqref{eq: inv smooth condition z} for the case of arc reversal if and only if $w_0 \cdot ([n] \setminus I) \in \M(w_0 \cdot \w)$ contributes to the set in \eqref{eq: inv smooth condition z} for the case of reflection. Therefore, the result follows in these cases.

For case of rotation, recall from \Cref{prop: rigid transformations} that $\M(\z) = \M(w)^{+s}$. For $I, J \in \binom{[n]}{k}$, $|I \cap J| = k-1$ if and only if $|I^{+s} \cap J^{+s}| = k-1$. Therefore $I \in \M(\w)$ contributes to the set in \eqref{eq: inv smooth condition w} if and only if $I^{+s} \in \M(\z)$ contributes to the set in \eqref{eq: inv smooth condition z}, so the result also follows for this case.
\end{proof}

\begin{Rem}\label{rem:circularBruhat}
In \cite{Postnikov.2006}, Postnikov defined a partial order on
decorated permutations called circular Bruhat order.  This order
determines the containment relations on positroid varieties just as
Bruhat order determines the containment relation on Schubert
varieties.    The covering relations in circular Bruhat order on $\SSnk$
are determined by exchanging a simple crossing with a simple alignment
by \cite[Thm. 17.8]{Postnikov.2006}.  Since the rigid
transformations preserve simple crossings and simple alignments, we
see these operations are order preserving under circular Bruhat order.
\end{Rem}

\section{Reduction to Connected Positroids}\label{sec:connected}

Recall that by \Cref{thm:ARW}, each positroid $\M$ on $[n]$ can be
uniquely constructed by choosing a non-crossing partition $B_{1}\sqcup
\cdots \sqcup B_{t}$ of $[n]$, and then putting the structure of a
connected positroid $\M_i$ on each block $B_{i}$, so $\M= \M_{1}
\oplus \cdots \oplus \M_{t}$.  The noncrossing partition also
determines a decomposition of the chord diagram of the associated
decorated permutation into connected components as a union of directed
arcs inscribed in the plane.  We will show that a positroid variety is smooth if
and only if the positroid varieties corresponding with each connected
component are smooth. By direct analysis of the Jacobian matrix at a $T$-fixed point, we show that this matrix can be decomposed into a block diagonal matrix corresponding with the connected components of the associated positroid. We will need a slight refinement of \Cref{cor:SIFcor}.

\begin{Lemma} \label{Lemma: concat matroid}
Let $\M$ be a positroid on ground set $[n]$.  If $\M$ is not
connected, then up to a possible cyclic shift, it has a decomposition
of the form
\[
\M = \M_1 \oplus \M_2^{+n_1} = \{I \cup J^{+n_{1}}: I \in M_{1}, J\in M_{2} \}.
\]
where $n=n_{1}+ n_{2}$, $\M_{1}$ and $\M_{2}$ are positroids on ground
sets $[n_{1}]$ and $[n_{2}]$ respectively. 
\end{Lemma}

\begin{proof} 
By \cite[Prop 7.4]{Ardila-Rincon-Williams}, if $\M$ is not connected
then we can assume it is the direct sum of two positroids $\M_1$ and
$\M_2$ on disjoint cyclic intervals.  Using \Cref{Lemma: invariances}
and the fact that positroids are closed under the transformation of
rotating the set $[n]$, we can assume that $\M_1$ has ground set
$[1,n_{1}]$ and $\M_2$ has ground set $[n_{1}+1,n]$.
\end{proof}

As discussed in \Cref{Sec: perm matrices suffice}, for a set $I \in \M$, one may classify the point $A_I \in \Pi_\M$ as a smooth or singular point by computing the rank of the Jacobian matrix, $Jac(\M)$, for $\Pi_{\M}$ evaluated at $A_I$. In the regime where $\M = \M_1 \oplus \M_2^{+n_1}$, as in \Cref{Lemma: concat matroid}, $Jac(\M)|_{A_I}$ can be written as a block diagonal matrix as follows.

\begin{Lemma} \label{lem: Jacob.decomp}
Let $\M_i \subseteq \binom{[n_i]}{k_i}$ be positroids for $i \in \{1,2\}$, and let $\M = \M_1 \oplus \M_2^{+n_1}$. Let $I_i \in \M_i$ for $i \in \{1,2\}$, and let $I = I_1 \sqcup I_2^{+n_1} \in \M$. Then the Jacobian matrix $Jac(\M)|_{A_I}$ can be written as a block diagonal matrix whose first and second blocks are, up to signs of entries, the matrices $Jac(\M_i)|_{A_{I_i}}$ for $i \in \{1,2\}$, and whose third block has rank $k_1(n_2-k_2) + k_2(n_1-k_1)$.
\end{Lemma}

\begin{proof}
Recall that for a positroid $\M \subset \binom{[n]}{k}$, the rows of $Jac(\M)$ are indexed by the nonbases in $\QQ(\M) = \binom{[n]}{k} \setminus \M$, and the columns of $Jac(\M)$ are indexed by variables $x_{ij}$, where $i \in [k]$ and $j \in [n]$. Consider $I \in \M$. By \Cref{Th: positroid = initial sets}, if $[u,v] \in \GIkn$ is the Grassmann interval corresponding to $\M$, then there is some $y \in [u,v]$ such that $I=y[k]$. By \Cref{lem: Jac a partial perm matrix}, $Jac(\M)|_{A_I}$ is, up to the signs of the entries, a partial permutation matrix whose nonzero entries occur exactly in the cells $(J,x_{st})$, where $J \in \QQ(\M)$ satisfies $|I \cap J|=k-1$,  $I \setminus J = \{y_s\}$, and $J \setminus I = \{t\}$. 

Set $n=n_1+n_2$, $k=k_1+k_2$, $Jac=Jac(\M)$, $Jac_i=Jac(\M_i)|_{A_{I_i}}$, and $\QQ_i=\QQ(\M_i) \subset \binom{[n_i]}{k_i}$ for $i \in \{1,2\}$. Let $[u^{(i)},v^{(i)}]$ be the Grassmann interval corresponding to $\M_i$, and let $y^{(i)}=u(I_i,v) \in [u^{(i)},v^{(i)}]$, as in \Cref{def: uAv}. By considering the decorated permutations corresponding to $\M$ and the $\M_i$ and the associated decomposition, as in \Cref{cor:SIFcor},  $v_j=v^{(1)}_j$ for $j \in [k_1]$ and $v_{j+k_1}=v^{(2)}_j+n_1$ for $j \in [k_2]$. Since every element of $I_2^{+n_1}$ is greater than every element of $v[k_1]$, then the construction of $y=u(I,v)$ replicates the constructions of the $y^{(i)}$ so that  $y_j=y^{(1)}_j$ for $j \in [k_1]$ and $y_{j+k_1}=y^{(2)}_j+n_1$ for $j \in [k_2]$. 


By definition, $\M = \M_1 \oplus \M_2^{+n_1}$. Therefore, the collection of nonbases can be partitioned as $\QQ = B \sqcup C \sqcup D \sqcup E \sqcup F$, where the sets in the partitioned are defined as
\begin{align*} 
B &= \QQ_1 \oplus \{I_2^{+n_1}\} , \hspace{.2in} C = \{I_1\} \oplus \QQ_2^{+n_1}, \\
D &= \Big\{J \in \binom{[n]}{k} \; :
\; |J \cap [n_1]|> k_1 \text{ or } |J \cap [n_1+1,n]| > k_2 \Big\}, \\
E &=  \QQ_1 \oplus \Big(\binom{[n_1+1,n]}{k_2} \setminus \{I_2^{+n_1}\} \Big), \text{ and}\hspace{.2in}
F=  (\M_1 \setminus \{I_1\}) \oplus \QQ_2^{+n_1}.
\end{align*}
We also partition the set $[k] \times [n]$ =  \textbf{(I)} $\sqcup$ \textbf{(II)} $\sqcup$ \textbf{(III)} in the following way.
\begin{itemize}
\item \textbf{(I)} = $[k_1] \times [n_1]$
\item \textbf{(II)} = $[k_1+1,k] \times [n_1+1,n]$
\item \textbf{(III)} = $([k_1+k,k] \times [n_1]) \cup ([k_1] \times [n_1+1,n])$
\end{itemize}

The partitions above yield partitions of the row and column sets of $Jac|_{A_I}$. From these partitions, $Jac|_{A_I}$ can be decomposed as a block diagonal matrix in the following way.

\noindent \textbf{\underline{Block 1}:} ($J \in B$) The rows of the first block will be indexed by the sets $J \in B$, and the columns of this block will be indexed by $x_{rc}$ with pairs $(r,c)$ in \textbf{(I)}. The maps $J_1 \mapsto J_1 \sqcup I_2^{+n_1}$ for $J_1 \in \QQ_1$ and $(r,c) \mapsto (r,c)$ for $(r,c) \in [k_1] \times [n_1]$ together give a bijection between the cells in $Jac_1$ and the cells in this upper left block of $Jac|_{A_I}$.

Since $I = I_1 \sqcup I_2^{+n_1}$, then a set $J = J_1 \sqcup I_2^{+n_1} \in \QQ_1 \oplus \{I_2^{+n_1}\}$ satisfies $|I \cap J| = k-1$ if and only if $|I_1 \cap J_1| = k_1-1$. By \Cref{lem: Jac a partial perm matrix} and \Cref{Th: positroid = initial sets}, for any such $J_1 \in \QQ_1$, the unique nonzero entry of $Jac_1$ in row $J_1$ is in column $x_{st}$, where $I_1 \setminus J_1 = \{y^{(1)}_s\}$ and $J_1 \setminus I_1 = \{t\}$. Since $I = y[k] = y^{(1)}[k_1] \cup y^{(2)} [k_2]^{+n_1}$, then $I \setminus J = \{y_s\}$ and $J \setminus I = \{t\}$. Hence, for this same pair $(s,t)$, entry $(J, x_{st})$ of $Jac|_{A_I}$ will be nonzero. Note also that this unique pair is in \textbf{(I)}.

Therefore, this upper left block with rows indexed by $B$ and columns indexed by \textbf{(I)} looks, up to sign, like $Jac_1$. Furthermore, all entries in the rows indexed by $B$, but outside of the columns indexed by \textbf{(I)}, will be zeros.

\noindent \textbf{\underline{Block 2}:} ($J \in C$) This case is similar to the previous case. The rows of the second block will be indexed by the sets $J \in C$, and the columns of this block will be indexed by $x_{rc}$ with $(r,c)$ in \textbf{(II)}. The maps $J_2 \mapsto I_1 \sqcup J_2^{+n_1}$ for $J_2 \in \QQ_2$ and $(r,c) \mapsto (r+k_1,c+n_1)$ for $(r,c) \in [k_2] \times [n_2]$ together give a bijection between the cells in $Jac_2$ and the cells in the second block of $Jac|_{A_I}$. A similar argument to that of the Block 1 case with elements and indices shifted by these maps shows that this second block looks, up to sign, like $Jac_2$ and that all entries in the rows indexed by $C$, but outside of the columns indexed by \textbf{(II)}, will be zeros.


\noindent \textbf{\underline{Block 3}:} ($J \in D \sqcup E \sqcup F$) The rows of the third block will be indexed by the sets $J \in D \sqcup E \sqcup F$, and the columns of this block will be indexed by pairs $(r,c)$ in set \textbf{(III)}. First, consider $J = J_1 \cup J_2 \in D$, where $J_1 = J \cap [n_1]$ and $J_2 = J \cap [n_1+1,n]$. Then $J$ satisfies $|I \cap J|=k-1$ if and only if either 
\begin{enumerate}
\item[(i)] $J_1 = I_1 \setminus \{y_s\}$ for some $y_s \in I_1$ and $J_2 = I_2^{+n_1} \cup \{t\} \in \binom{[n_1+1,n]}{k_2+1}$, or
\item[(ii)] $J_1 = I_1 \cup \{t\} \in \binom{[n_1]}{k_1+1}$ and $J_2 = I_2^{+n_1} \backslash \{y_s\}$ for some $y_s \in I_2^{+n_1}$. 
\end{enumerate}

In case (i), the column pair $(s,t)$ satisfying \Cref{lem: Jac a partial perm matrix} has $s \in [k_1]$ since $y_s \in [n_1]$ and $t \in [n_1+1,n]$. There are $k_1$ choices for $s$ and $n_2-k_2$ choices for $t$. In case (ii), the pair $(s,t)$ satisfying \Cref{lem: Jac a partial perm matrix} has $s \in  [k_1+1,k]$ since $y_s \in [n_1+1,n]$ and $t \in [n_1]$. There are $n_1-k_1$ choices for $t$ and $k_2$ choices for $s$. Thus, the rank of this block restricted to the rows in $D$ is $k_1(n_2-k_2)+k_2(n_1-k_1)$. Furthermore, all nonzero entries in the rows indexed by $D$ occur in the columns in \textbf{(III)} by arguments similar to the proof of \Cref{lem: Jac a partial perm matrix}.

We claim that for $J \in E \sqcup F$, $|I \cap J| < k-1$. Hence, $Jac|_{A_I}$ has all zeros in row $J$ by \Cref{lem: Jac a partial perm matrix}. To prove the claim, consider $J = J_1 \sqcup J_2 \in E$, where $J_1 \in \QQ_1$ and $J_2 \in \binom{[n_1+1,n]}{k_2} \setminus \{I_2^{+n_1}\}$. Since $J_1$ is in $\QQ_1$, so $J_1 \neq I_1$, then $|I_1 \cap J_1| \leq k_1-1$. Similarly, since $J_2 \neq I_2^{+n_1}$, then $|I_2^{+n_1} \cap J_2| \leq k_2-1$. Therefore, $|I \cap J| \leq (k_1-1)+(k_2-1) = k-2$. Next, consider $J = J_1 \sqcup J_2^{+n_1} \in F$, where $J_1 \in \M_1 \setminus \{I_1\}$ and $J_2 \in \QQ_2$. Since $J_1 \neq I_1$, then $|I_1 \cap J_1| \leq k_1-1$. Similarly, since $J_2$ is in $\QQ_2$, and therefore is not $I_2$, then $|I_2 \cap J_2| = |I_2^{+n_1} \cap J_2^{+n_1}| \leq k_2-1$. Therefore, $|I \cap J| \leq (k_1-1)+(k_2-1) = k-2$.
\end{proof}

\begin{Lemma} \label{Lemma: concat alignments}
Let $\M_{1}$ and $\M_{2}$ be positroids on grounds sets $[n_{1}]$ and
$[n_{2}]$. Let $\M = \M_1\oplus \M_{2}^{+n_{1}}$ and let $\w=w^{(1)}
\oplus w^{(2)} $ be the associated decorated permutation.  Let $k_{i}$
be the number of anti-exceedances of $w^{(i)}$.  Then, 
\[
\# \aligns(\w) \; = \; \# \aligns(w^{(1)}) + \# \aligns(w^{(2)}) + k_1(n_2-k_2) + k_2(n_1-k_1).
\]
\end{Lemma}

\begin{figure}
\begin{center}
   \includegraphics[scale=0.12]{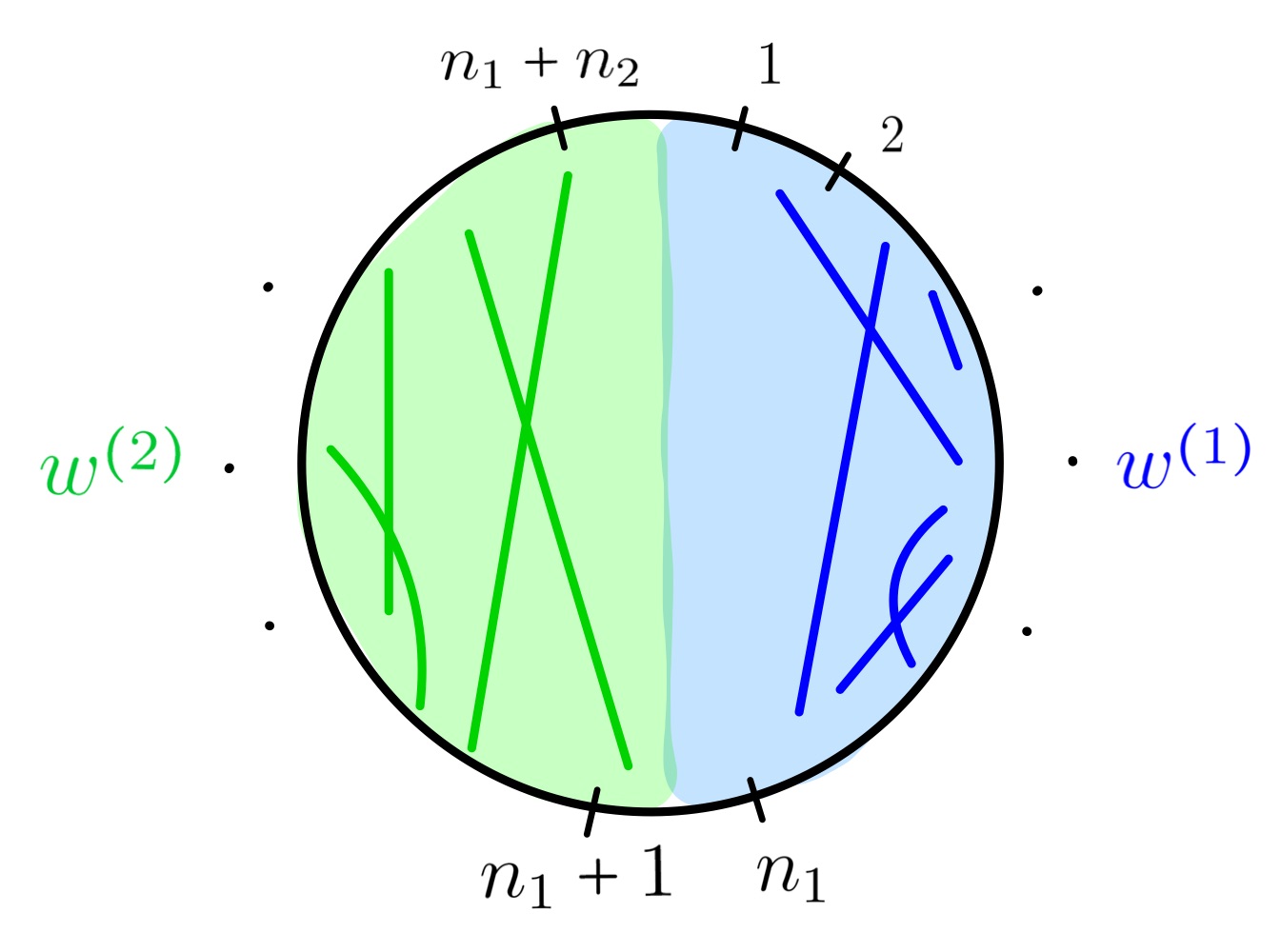}
\end{center}
\caption{Chord diagram for $w^{(1)} \oplus w^{(2)}$.}
   \label{Fig: concat}
\end{figure}

\begin{proof}
Consider the chord diagram for $\w$ partitioned by a line separating the connected components of the noncrossing partition corresponding with $w^{(1)}$ and $w^{(2)}$.
See \Cref{Fig: concat}.  Every alignment of $w^{(i)}$ remains an alignment
of $\w$. In addition, every anti-exceedance arc of $w^{(i)}$ forms an alignment with
every exceedance arc of $w^{(i')}$, where $i' =
\frac{3+(-1)^{i-1}}{2}$. Since $w^{(i)}$ has $k_i$ anti-exceedances
and $n_i-k_i$ exceedances, the result follows.
\end{proof}

\begin{Lemma} \label{Lemma: concat sm-sing}
Let $\M_{1}$ and $\M_{2}$ be positroids on grounds sets $[n_{1}]$ and
$[n_{2}]$ of rank $k_{1}$ and $k_{2}$, respectively. Let $\M =
\M_1\oplus \M_{2}^{+n_{1}}$.  For any $I = I_1 \sqcup I_2^{+n_1} \in \M$,
$A_I$ is a smooth point of $\Pi_{\M}$ if and only if the
$A_{I_1}$ and $A_{I_{2}}$ are both smooth points of $\Pi_{\M_1}$ and
$\Pi_{\M_2}$, respectively.
\end{Lemma}

\begin{proof}
By \Cref{thm:ARW} and the construction, $\M$ is a positroid of rank $k=k_{1}+k_{2}$ on ground set $[n]$ for $n=n_1+n_2$. Let $\w = w^{(1)} \oplus w^{(2)}$ be the corresponding decomposition of the associated decorated permutation, as in \Cref{Lemma: concat alignments}. Let $Jac = Jac(\M)|_{A_I}$ and $Jac_i = Jac(\M_i)|_{A_{I_i}}$ for $i \in \{1,2\}$. Since $\text{rank}(Jac)$ is bounded above by $\#\aligns(\w)$, then by \eqref{eq: sing condition Jac algn}, $A_I$ is a smooth point of $\Pi_\M$ if and only if $\text{rank}(Jac) = \#\aligns(\w)$. Similarly, $A_{I_i}$ is a smooth point of $\Pi_{\M_i}$ if and only if $\text{rank}(Jac_i) = \#\aligns(w^{(i)})$.

By \Cref{lem: Jacob.decomp}, $\text{rank}(Jac) = \text{rank}(Jac_1) + \text{rank}(Jac_2) + k_1(n_2-k_2) + k_2(n_1-k_1)$. It follows from \Cref{Lemma: concat alignments} that $A_I$ is a smooth point of $\Pi_{\M}$ if and only if
\begin{equation} \label{eq:Jac.rank.sum}
\text{rank}(Jac_1) + \text{rank}(Jac_2) = \#\aligns(w^{(1)}) + \#\aligns(w^{(2)}).
\end{equation}
Once again, since $\text{rank}(Jac_i) \leq \#\aligns(w^{(i)})$, then \eqref{eq:Jac.rank.sum} holds if and only if $\text{rank}(Jac_i) = \#\aligns(w^{(i)})$ for both $i=1,2$, which holds if and only if $A_{I_i}$ is a smooth point of $\Pi_{\M_i}$ for both $i=1,2$. 
\end{proof}

\begin{Cor} \label{Cor: concat smooth/sing}
Let $\M= \M_{1} \oplus \cdots \oplus\M_{t}$ be a positroid decomposed
into its connected components.  Then $\Pi_{\M}$ is smooth if and only
if $\Pi_{\M_{i}}$ is smooth for each $i \in [t]$.
\end{Cor}

\begin{proof}
If $\M$ is connected the statement holds, so assume $\M$ is not
connected.  We know that $\Pi_{\M}$ is smooth if and only if all
cyclic rotations of $\M$ correspond with smooth positroid varieties by
\Cref{Lemma: invariances}.  Therefore, by \Cref{Lemma: concat
matroid}, we can assume it has a decomposition of the form $\M = \M_1
\oplus \M_2^{+n_1}$ where $n=n_{1}+ n_{2}$ for some $1\leq
n_{1},n_{2}<n$, and $\M_{1}$ and $\M_{2}$ are positroids on ground sets
$[n_{1}]$ and $[n_{2}]$, respectively.  By \Cref {Lemma: concat
sm-sing}, $\Pi_{\M}$ has a singular $T$-fixed point $A_{I}$ for $I\in
\M$ if and only if either $\Pi_{\M_1}$ or $\Pi_{\M_2}$ has a singular
$T$-fixed point.  Therefore, the result holds by \Cref{Th:
positroid = initial sets}, \Cref{Cor: singular at piece of
partition}, and induction on $t$.   
\end{proof}




\begin{Rem}\label{remark:reducution.to.SIF}
By \Cref{cor:SIFcor}, the connected components of positroids are
associated with decorated SIF permutations.  In order to complete the
proof of \Cref{The: Main theorem}, it remains to characterize the
decorated SIF permutations indexing smooth positroid varieties by
\Cref{Cor: concat smooth/sing}.  We will show that these are in bijection
with the spirograph permutations, discussed further in the next section.  
\end{Rem}

\section{Crossed Alignments and Spirograph Permutations}
\label{Sec: spirograph direction}

Recall the class of spirograph permutations introduced in \Cref{sec:
intro}.  It is straightforward to show that spirograph permutations
always index smooth positroid varieties.  We will also show in this
section that direct sums of spirograph permutations on a noncrossing
partition are exactly the same as decorated permutations that have no
crossed alignments.  In particular, we prove the equivalence of Parts
(5), (6), and (7) in \Cref{The: Main theorem}.

By \Cref{def:spirograph}, the spirograph permutations are a subset of the
decorated SIF permutations.  The chord diagram
and the positroid associated to a spirograph permutation is always
connected. The two decorated permutations in $S_1^{\circ \bullet}$ are
spirographs, and for $n>1$, there are $n-1$ distinct spirograph
permutations corresponding with $m\in[n-1]$.  Therefore, the generating
function for spirograph permutations in $S_n$ for $n \geq 1$ is
\begin{equation}\label{eq:spirgraph.enumeration}
S(x) = 2x+x^{2}+2x^{3}+3x^{4}+\ldots = \frac{2x - 3 x^{2} + 2 x^3}{(1 - x)^2}.
\end{equation}

\begin{Rem}\label{rem:spiro.circular}
Each $\SSnk$ has a unique spirograph permutation, denoted
$\pi_{n,k}$. This spirograph permutation is defined so that
$\pi_{n,k}(i) = i+k$ for all $i$ and is the unique maximal element in
circular Bruhat order in $\SSnk$ for $n>1$.  For $k=n=1$, $\pi_{1,1}$
is the unique decorated permutation in $\SSnk$, which consists of a
single clockwise fixed point.  For $k=0, n=1$, $\pi_{1,0}$ is the
unique decorated permutation in $\SSnk$, which has a counterclockwise
fixed point.
\end{Rem}

\begin{Lemma} \label{Lemma: star perm} \normalfont If $\w$ is the unique the spirograph permutation in $\SSnk$, then $\Pi_{\w}$ is $\Gkn$, which is a smooth
variety.
\end{Lemma}
\begin{proof}
Observe that the chord diagram of a spirograph permutation $\w$ has no
alignments, so the codimension of $\Pi_{\w}$ as a subvariety of $\Gkn$
is zero by \Cref{thm:Alignments}.  Hence, $\Pi_{\w} = \Gkn$.
\end{proof}

\begin{proof}[Proof of \Cref{The: Main theorem}, $(6) \Rightarrow
(1)$] Let $D(\w)$ be a disjoint union of spirographs corresponding to
the decomposition $\w = w^{(1)} \oplus \cdots \oplus w^{(t)}$ of $\w$
into decorated SIF permutations with $\M(\w) = \M(w^{(1)}) \oplus
\cdots \oplus \M(w^{(t)})$ using \Cref{cor:SIFcor}.  By \Cref{Lemma: star
perm}, the positroid varieties $\Pi_{w^{(i)}}$ are all smooth. Hence,
$\Pi_{\w} = \Pi_{\M(\w)}$ is smooth by \Cref{Cor: concat smooth/sing}.
\end{proof}

\begin{proof}[Proof of \Cref{The: Main theorem}, $(5) \Leftrightarrow (6)$]
By definition of a chord diagram $D(\w)$, arcs from distinct components of the associated noncrossing partition must be drawn so that they do not intersect. It follows that, if $\w$ contains a crossed alignment, then all arcs involved in the crossed alignment must be contained in the same connected component of $D(\w)$.

A spirograph permutation has no crossed alignments since it has no alignments. By the observation above, if every connected component of $D(\w)$ is a spirograph, then $\w$ has no crossed alignments.

For the converse, observe that the property of
containing a crossed alignment is invariant under rotation of the chord diagram.  Hence, by
\Cref{cor:SIFcor} and the fact that a crossed alignment is contained a single connected component of $D(\w)$, the argument may be completed by assuming that $\w$ is a decorated SIF
permutation that is not a spirograph and showing that $\w$ has a crossed
alignment.  For $n < 4$, $\w$ is either a spirograph permutation or has
more than one connected component.  Therefore, we must have $n \geq 4$, in which case we will denote $\w$
simply by $w$ since it has no fixed points.


Since $w$ has no fixed points, $m_i:= w(i)-i\
(\mathrm{mod}\ n) \in [n-1]$ for each $i\in [n]$. Since $w$ is not a spirograph permutation, the $m_{i}$ are not all equal.
Therefore, taking indices mod $n$, there must exist some $i
\in [n]$ such that $m_{i}<m_{i+1}$.  Observe from the chord diagram that $m_{i}<m_{i+1}$ implies the arcs $(i \mapsto w(i))$
and $(i+1 \mapsto w(i+1))$ form a crossing. Use the crossing arcs to create three disjoint cyclic intervals.
Let $A = [w(i)+1,w(i+1)-1]^{cyc}$. The fact that $m_i < m_{i+1}$ implies that $A$ is nonempty.  Let $B = [w(i+1),i-1]^{cyc} \setminus [i,w(i+1)-1]^{cyc}$ and $C =[i+2,w(i)]^{cyc} \setminus [w(i)+1,i+1]^{cyc}$. See \Cref{Fig: nonstar crossing}.

\begin{figure}
\begin{center}
   \includegraphics[scale=0.18]{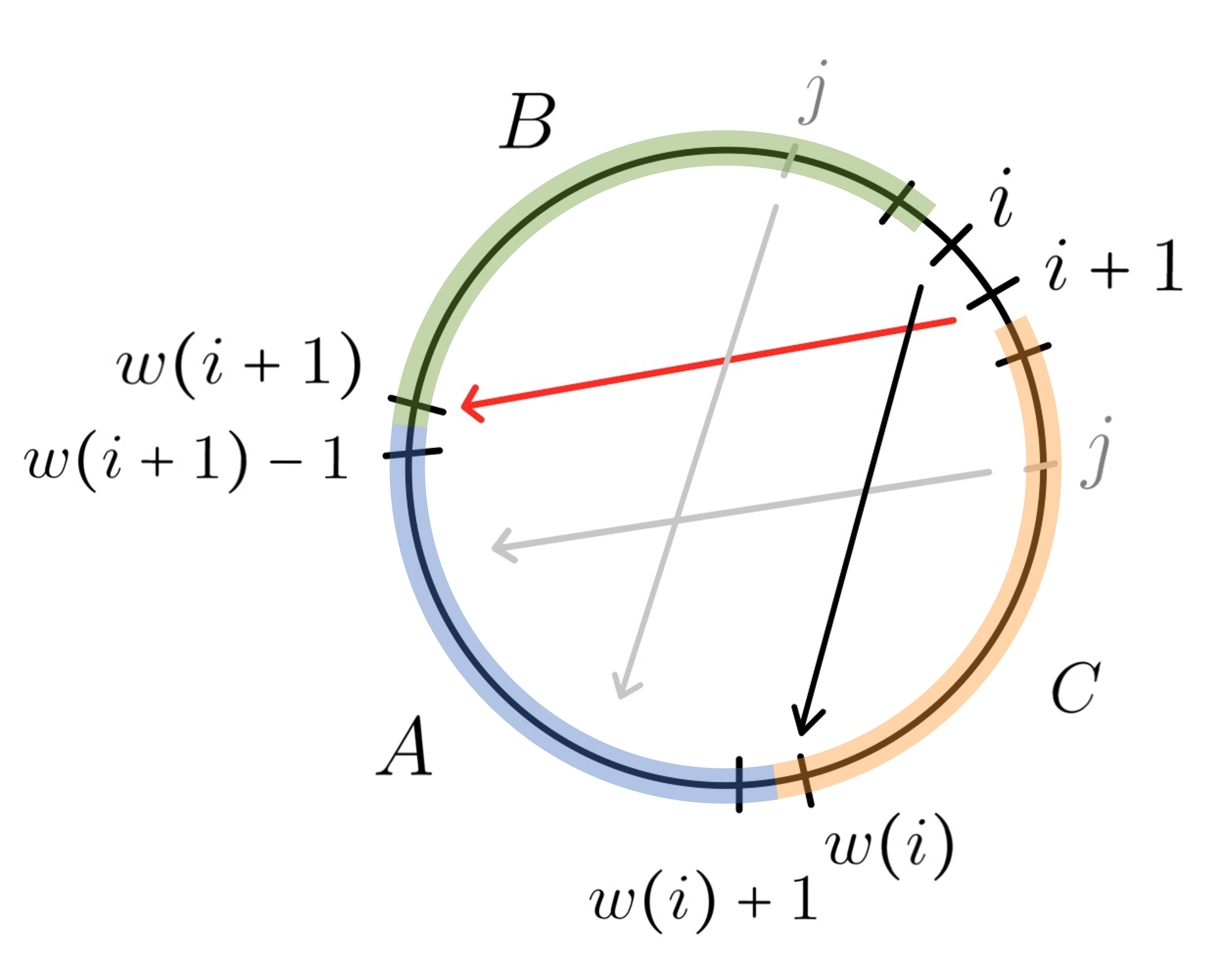}
\hspace{.15in}
   \includegraphics[scale=0.18]{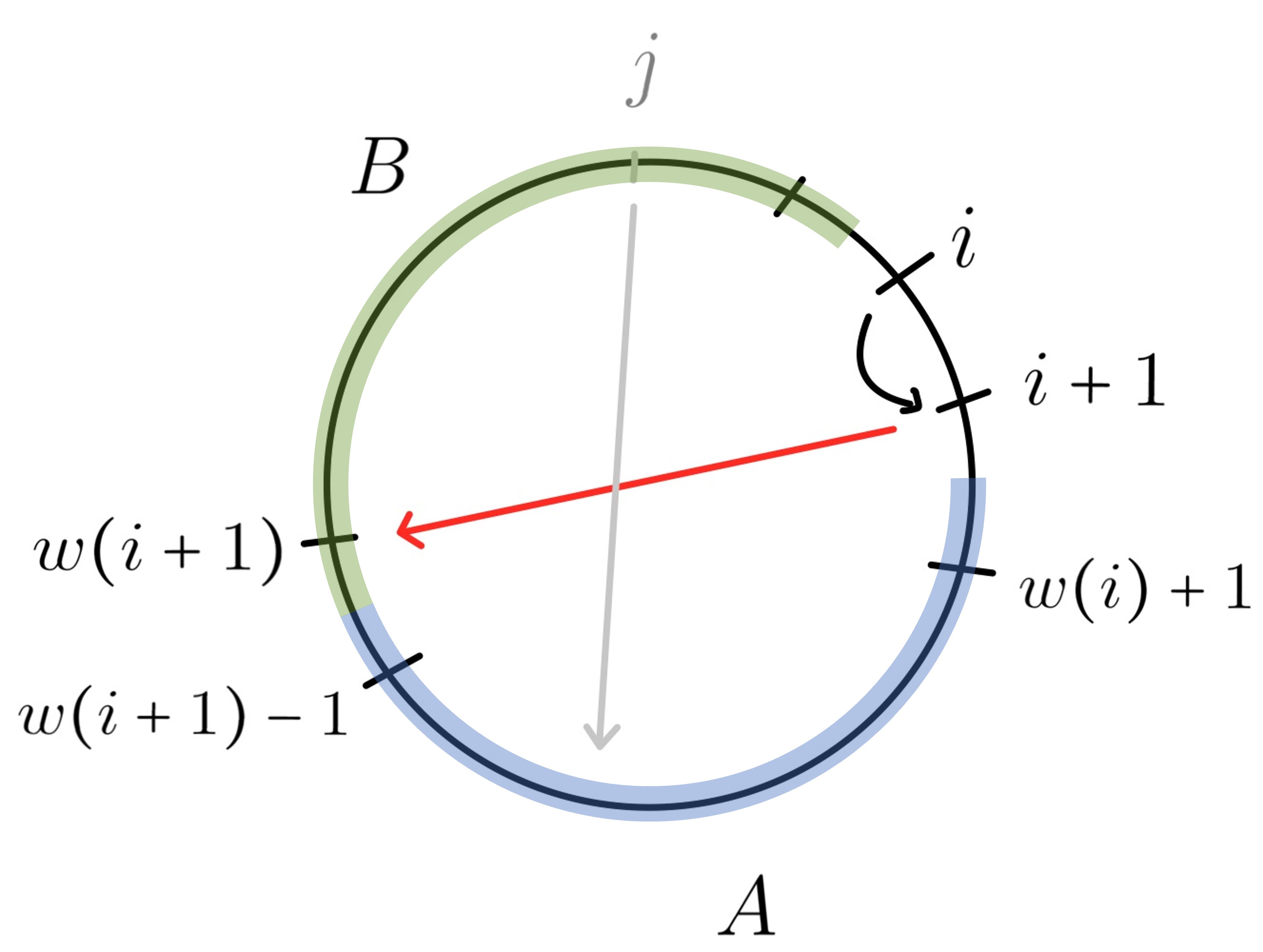}
\hspace{.15in}
   \includegraphics[scale=0.18]{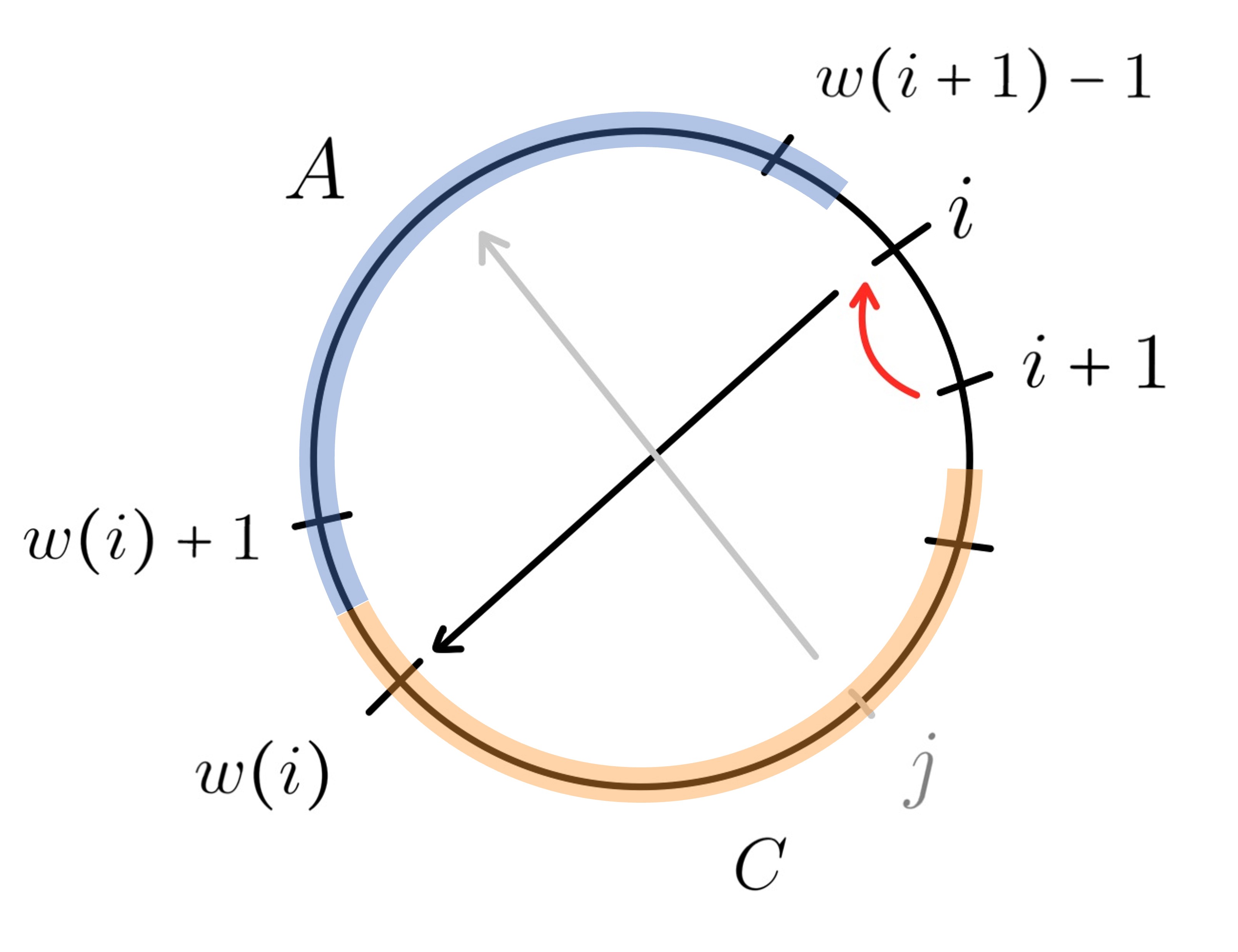}
\end{center}
\caption{Chord diagrams with $(i \mapsto w(i))$ and $(i+1 \mapsto w(i+1))$ crossing.}
   \label{Fig: nonstar crossing}
\end{figure}

Since $w$ is a SIF permutation, the arcs incident to elements of $A$ cannot
form an isolated connected component. Thus, there must be an arc
mapping some $j \in [n] \setminus (A \cup \{i,i+1\})=B \cup C $ to $w(j)
\in A$.  If $j \in B$, then $(i \mapsto w(i), j \mapsto w(j))$ is an
alignment crossed by $(i+1 \mapsto w(i+1))$. If $j \in C$, then $(j
\mapsto w(j), i+1 \mapsto w(i+1))$ is an alignment crossed by $(i
\mapsto w(i))$. Therefore, $w$ contains a crossed alignment.
\end{proof}

\begin{proof}[Proof of \ref{The: Main theorem}, $(6) \Leftrightarrow (7)$]
The unique spirograph permutation $\pi_{n,k}$ in $\SSnk$, mentioned in \Cref{rem:spiro.circular}, has shifted anti-exceedance sets $I_r(\pi_{n,k}) = [k]^{+(r-1)}$ for all $r \in [n]$. Since $[k]^{+(r-1)}$ is the minimal element of $\binom{[n]}{k}$ under $\prec_r$, it follows from \Cref{th:decoratedpermtopositroid} that $\M(\pi_{n,k}) = \binom{[n]}{k}$. Hence, spirograph permutations correspond exactly to uniform matroids. By the correspondence in \Cref{cor:SIFcor}, a positroid $\M = \M(\w)$ is a direct sum of uniform matroids if and only if $\w$ is a direct sum of spirograph permutations.
\end{proof}

\section{Anti-Exchange Pairs and Crossed Alignments}
\label{Sec: anti-exchange pairs and crossed aligns}

In this section, we relate decorated permutations with crossed
alignments to singular positroids varieties as determined by the
vertex degrees of the matroid Johnson graphs in Parts (1)-(4) in
\Cref{The: Main theorem}.  By the end of \Cref{Sec: spirograph
direction}, we proved Parts (5) $\Leftrightarrow$ (6) $\Leftrightarrow$ (7) and (6) $\Rightarrow$ (1) in
\Cref{The: Main theorem}.  To finish the proof of \Cref{The:
Main theorem}, it suffices to prove (1) $\Rightarrow$ (5), namely if
$\Pi_{\w}$ is smooth, then $\w$ has no crossed alignments. This final step will be proved using the characterization of
smoothness from
\Cref{cor:anti-exchangable.pairs.alignments} utilizing the nonbases in $\QQ(\w)$.  Specifically, for any
$\w \in \SSnk$, $\Pi_{\w}$ is singular if and only if there exists a
$J \in \M(\w)$ such that
\begin{equation} \label{eq: num anti-exch pairs less codim}
\# \{I \in \QQ(\w) \; : \; |I \cap J| = k-1 \} \; < \;  \# \aligns(w).
\end{equation}
Thus, we now carefully study the sets $I$ that occur on the left side
of this equation. 
The
following vocabulary refers back to the Basis Exchange Property of
matroids from \Cref{sec:background}.

\begin{Def} \label{Def: anti-exchange pair} For a set $J$ in a
positroid $\M \subseteq \binom{[n]}{k}$, let $a \in J$ and $b \in [n] \setminus J$. If $(J \setminus
\{a\}) \cup \{b\} \in \M$, we say the pair $(a,b)$ is an
\textit{exchange pair for} $J$ and that the values $a$ and $b$ are
\textit{exchangeable}.  Otherwise, $(J \setminus \{a\}) \cup \{b\}
\in \QQ = \binom{[n]}{k} \setminus \M$, in which case we say $(a,b)$ is an \textit{anti-exchange pair
for} $J$ and that the values $a$ and $b$ are \textit{not
exchangeable}.
\end{Def}

For a set $J \in \M \subseteq \binom{[n]}{k}$, pairs $(a,b)$ with $a
\in J$ and $b \notin J$ are in bijection with the vertices in
the full Johnson graph $J(k,n)$ adjacent to $J$. Exchange pairs
for $J$ are in bijection with the vertices adjacent to $J$ in the
matroid Johnson graph $J(\M)$. Anti-exchange pairs for $J$ are in
bijection with the vertices adjacent to $J$ in $J(k,n)$ that do not
appear in $J(\M)$.

\subsection{Characterizing Exchange Pairs in the Johnson Graph}
\label{Sec: anti-exchange pairs}
For a positroid $\M(\w)$, we will focus on the set $J = I_1(\w)$, which is in $\M(\w)$ by \Cref{cor:allIrs are in positroid}. The following lemma characterizes exchange pairs for $I_{1}(\w)$.

\begin{Lemma} \normalfont \label{Lemma: anti-exchange pair condition} For
$\w=(w,\co) \in \SSnk$, suppose that $a \in
I_1(\w)$ and $b \in [n] \setminus I_1(\w)$. Then $I = (I_1(\w) \setminus \{a\}) \cup \{b\}$
is in $\M(\w)$ if and only if $a < b$ and for every $r \in [a+1,b]$,
both of the following conditions hold:
\begin{enumerate}
\item[$(1)$] there exists $x \in [a,r-1]$ such that $w^{-1}(x) \in [r,n]$, and
\item[$(2)$] there exists $y \in [r,b]$ such that $w^{-1}(y) \in [1,r-1]$.
\end{enumerate}
\end{Lemma}

\begin{proof}
Let $I_r = I_r(\w)$.  Recall from \Cref{th:decoratedpermtopositroid} that $I$ is in $\M = \M(\w)$ if and only if $I_{r}\preceq_{r} I$ for all
$r\in[n]$. We use this characterization of $\M$ in terms of the elements of the Grassmann necklace to derive the conditions of the lemma.

Consider the case when $a>b$. Then $I$ is obtained by replacing an element of $I_1$ with a smaller element. By definition of the Gale order, $I$ cannot satisfy $I_1 \preceq I$, so $I \not \in \M$.  

For the remainder of the proof, assume $a<b$.  For $r \in [b+1,a]^{cyc}$, the inequality $a <_r b$ holds, which
implies that $I$ is obtained by replacing an element of $I_1$ with an element that is larger under $<_r$. Thus, $I_1 \prec_r I$. By \Cref{cor:allIrs are in positroid}, $I_1$ is in $\M$, so by \Cref{th:decoratedpermtopositroid}, $I_1$ must satisfy $I_r \preceq_r I_1$. Hence, the sequence of inequalities  $I_r \preceq_r I_1 \prec_r I$ holds for every $r \in [b+1,a]^{cyc}$. Therefore, $I \in
\M$ if and only if $I_r \preceq_r I$ for every $r \in [a+1,b]$.

Fix some $r
\in [a+1,b]$, so $a>_{r}b$. An arc $(i \mapsto j)$ with $i \in [r,n]$ and $j \in [a,r-1]$ is an anti-exceedance arc, but not an $r$-anti-exceedance arc, so $j \in I_1 \setminus I_r$. Similarly, an arc $(i \mapsto j)$ with $i \in [1,r-1]$ and $j \in [r,b]$ is an $r$-anti-exceedance arc, but not an anti-exceedance arc, so $j \in I_r \setminus I_1$. Therefore, the following statements hold.
\begin{equation}
\text{For } x \in [a,r-1], \; w^{-1}(x) \in [r,n] \Leftrightarrow x \in I_1 \setminus I_r.
\end{equation}
\begin{equation}
\text{For } y \in [r,b], \; w^{-1}(y) \in [1,r-1] \Leftrightarrow y \in I_r \setminus I_1.
\end{equation}
See \Cref{fig: exchange lemma conditions} for these two cases. From these observations, for the fixed $r \in [a+1,b]$, Conditions $(1)$ and $(2)$ in the statement of the lemma are equivalent to the following:
\begin{enumerate}
\item[$(1')$] $(I_1 \setminus I_r) \cap [a,r-1] \neq \emptyset$, and
\item[$(2')$] $(I_r \setminus I_1) \cap [r,b] \neq \emptyset$.
\end{enumerate}

\begin{figure*}[h!]
    \centering
  {\includegraphics[width=1.7in]{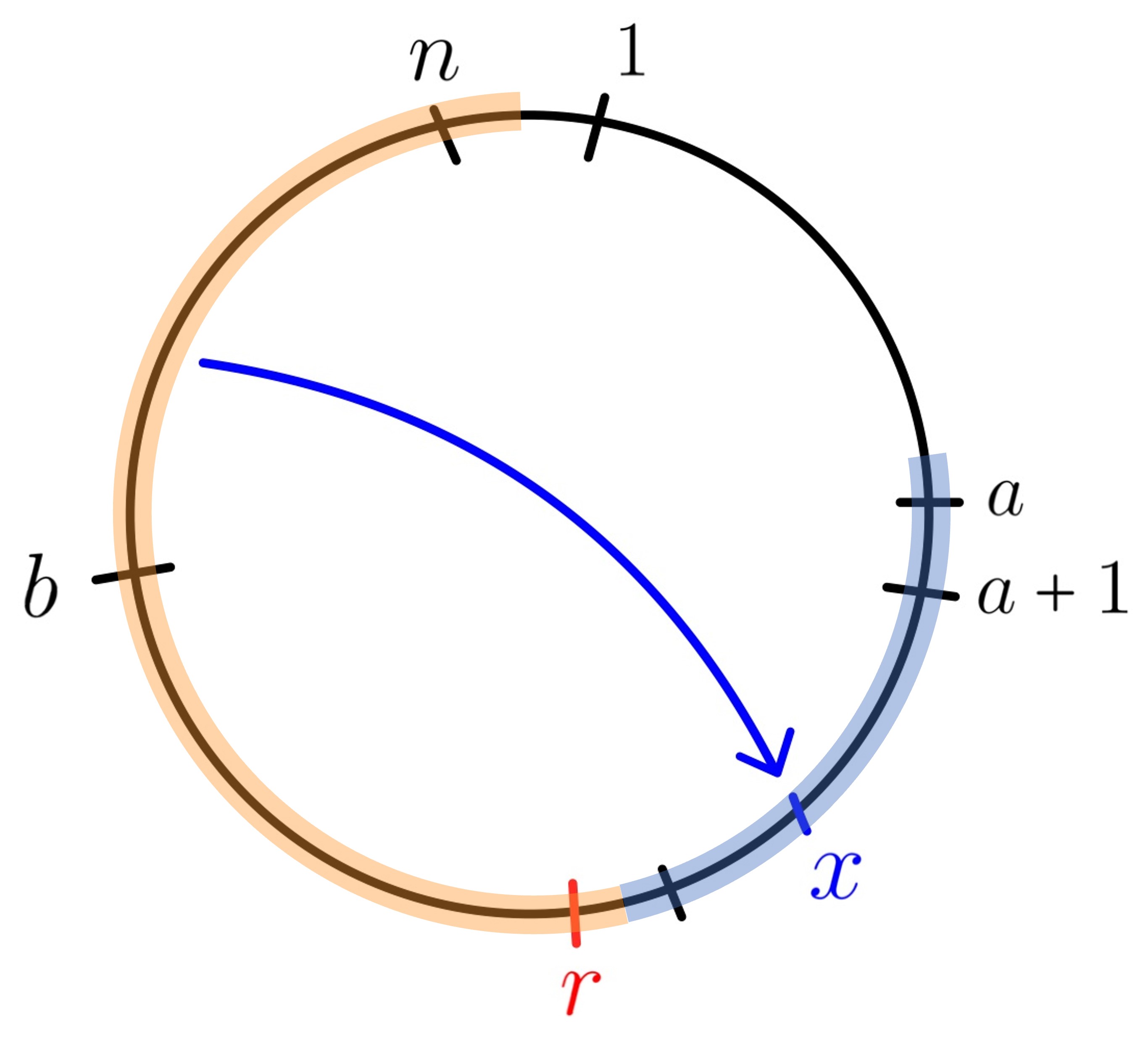}}\hspace{5em}%
  {\includegraphics[width=1.7in]{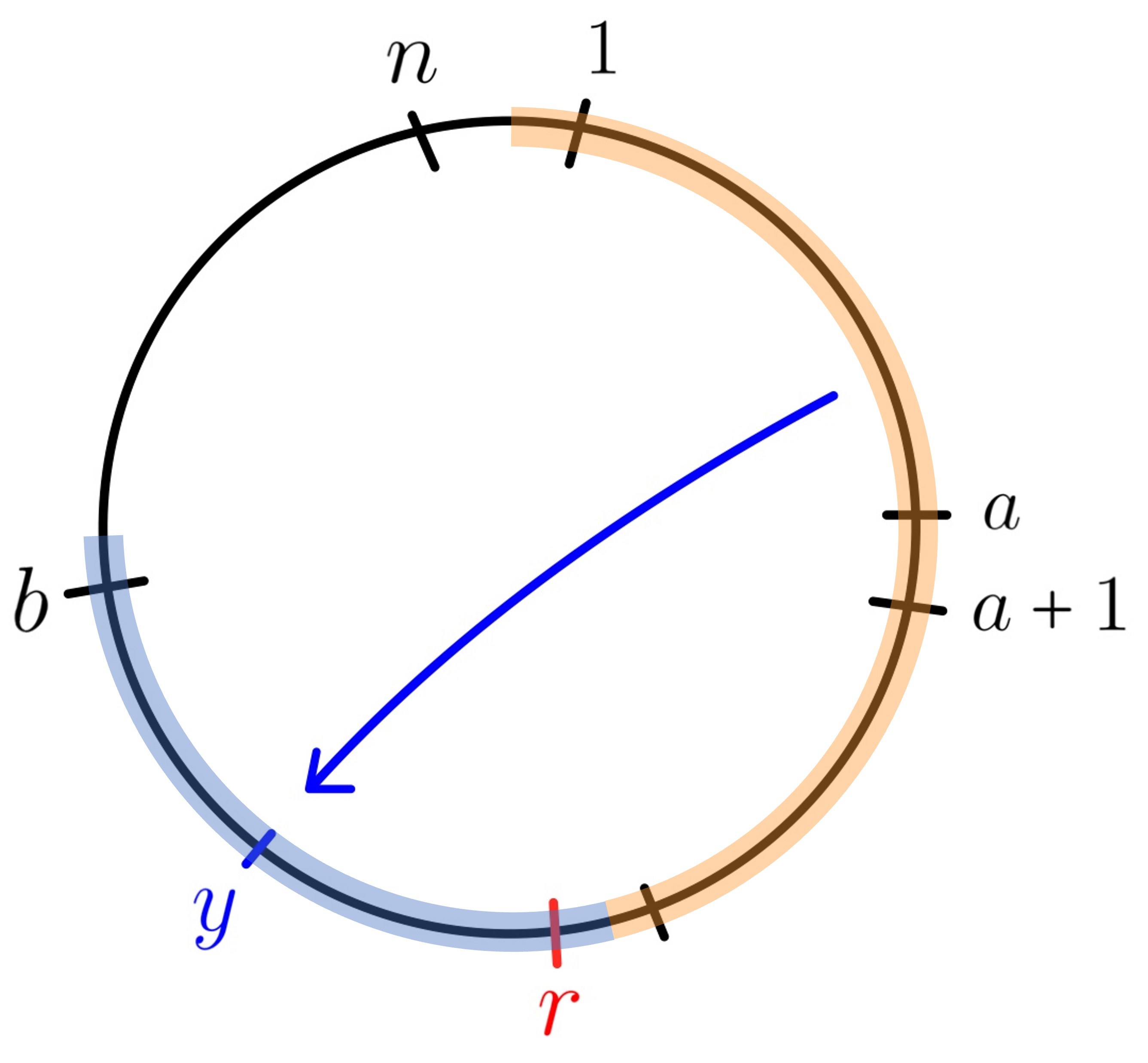}}
    \caption{Two conditions of \Cref{Lemma: anti-exchange pair condition}}
\label{fig: exchange lemma conditions}
\end{figure*}


We must now compare the sets $I_r$ and $I$ under the shifted Gale order $\prec_r$ for $r \in [a+1,b]$. Recall that $I$ is obtained from $I_1$ by exchange of the elements $a$ and $b$, where $a \in I_1$ and $b \notin I_1$. Thus we begin by writing $I_1 = \{j_1 <_r \cdots <_r j_k\}$ increasing under the $<_r$ order. Setting $f = |I_1 \cap [r,b]| +1$, $g = |I \cap [r,n]| +1$, and $h = |I \cap [r,a]^{cyc}|$, we refine the expression of $I_1$ in the following way,
\begin{align*} 
I_1 \, = \, \{\underbrace{j_1 \, <_r \, \cdots \, <_r \, j_{f-1}}_{\in [r,b]} \, <_r \, \underbrace{j_f \, <_r \, \cdots \, <_r \, j_{g-1}}_{\in [b+1,n]} \, <_r \, \underbrace{j_g \, <_r \, \cdots \, <_r \, j_{h-1} <_r \, a}_{\in [1,a]} \, <_r \, \underbrace{j_{h+1} \, <_r \, \cdots \, <_r \, j_k}_{\in [a+1,r-1]}\}.
\end{align*}
Here, $a = j_h$. Note that since $b \notin I_1$, then $j_{f-1} <_r b$. Write $I_r$ and $I$ under $<_r$ as
\begin{align*} 
I_r \, &= \, \{i_1 \, <_r \, \cdots \, <_r \, i_{f-1} \, <_r \, i_f \, <_r \, i_{f+1} \, <_r \, \cdots \, <_r \, \; \, i_{g} \;  \; <_r \; i_{g+1} \, <_r \, \cdots \;  <_r \, \; \; i_h \,  \; <_r \; i_{h+1} \, <_r \, \cdots \, <_r \, i_k\} \\ 
I \, &= \, \{j_1 \, <_r \, \cdots \, <_r \, j_{f-1} \, <_r  \, b \, \; <_r  \; \; j_f \; \;  <_r \, \cdots \, <_r \, j_{g-1} \, <_r \,  \; j_g \; \; \, <_r \, \cdots \, <_r \, j_{h-1} \, <_r \, j_{h+1} \, <_r \, \cdots \, <_r \, j_k\}.
\end{align*}
Then $I_r \preceq_r I$ if and only if all three of the following hold:
\begin{enumerate}
\item[(i)] $i_\ell \leq_r j_\ell$ for all $\ell \in [1,f-1] \cup [h+1,k]$;
\item[(ii)] $i_f \leq_r b$ and $i_{\ell+1} \leq_r j_{\ell}$ for all $\ell \in [f,g-1]$; and
\item[(iii)] $i_{\ell+1} \leq_r j_{\ell}$ for all $\ell \in [g,h-1]$.
\end{enumerate}
For $\ell \in [1,f-1] \cup [h+1,k]$, $I_r \preceq_r I_1$ already implies that $i_\ell \leq_r j_\ell$, so (i) is always satisfied. Thus, $I_r \preceq_r I$ if and only if (ii) and (iii) are satisfied.

Consider Condition (ii) above. We will show that (ii) is satisfied if and only if $(I_r \setminus I_1)\cap [r,b] $ is nonempty, as in Condition $(2')$.

We require the following observations about the elements $j_1, \ldots, j_{g-1}$. By construction, $\{j_1, \ldots, j_{g-1}\} = I_1 \cap [r,n]$. Each of these $j_\ell$ has $w^{-1}(j_\ell) \in [j_\ell,n]$, as in \Cref{fig: j_ell}. But then $j_\ell$ and $w^{-1}(j_\ell)$ also satisfy $r \leq_r j_\ell \leq_r w^{-1}(j_\ell) \leq_r n$.  In particular, this implies that every such $j_\ell$ is also an element of $I_r$ so that $\{j_1, \ldots, j_{g-1}\} \subseteq I_r$. 
Furthermore, $j_1 <_r \cdots <_r j_{g-1}$ appear in $I_r$ in the same relative order as they do in $I_1$, but there may be additional elements of $I_r$ that are interspersed among them. Therefore, for every $\ell \in [1,g-1]$, there is some $\ell' \geq \ell$ such that $i_{\ell'} = j_{\ell}$.

For the only if direction, (ii) implies $(2')$, assume $(I_r \setminus I_1) \cap [r,b]$ is empty. By construction, $I_1 \cap [r,b] = \{j_1, \ldots, j_{f-1}\}$. Since $\{j_1, \ldots, j_{f-1}\} \subset I_r$ by the previous paragraph and $(I_r \setminus I_1) \cap [r,b] = \emptyset$, then $I_r \cap [r,b] = \{j_1, \ldots, j_{f-1}\} = \{i_1, \ldots, i_{f-1}\}$. In particular, $i_f \notin [r,b]$, so $b <_r i_f$, violating Condition (ii).

Conversely, suppose $(I_r \setminus I_1) \cap [r,b]$ is nonempty. Then $|I_r \cap [r,b]| \geq |I_1 \cap [r,b]| +1 = f$. In particular, $i_f$ must be in $[r,b]$, so $i_f \leq_r b$. Furthermore, this extra element in $I_r \cap [r,b]$ shifts the elements $j_f, \ldots, j_{g-1}$ to the right in $I_r$. Specifically, for $\ell \in [f,g-1]$, $i_{\ell'} = j_{\ell}$ for some $\ell'>\ell$. Therefore, for $\ell \in [f,g-1]$, we have the inequalities $i_{\ell+1} \leq_r  i_{\ell'} = j_\ell$. Thus, $(I_r \setminus I_1) \cap [r,b] \neq \emptyset$ implies that Condition (ii) is satisfied.

Now, consider Condition (iii) above. We will show that (iii) is satisfied if and only if Condition $(1')$ holds. The argument is symmetric to the argument above for Condition (ii).

Consider any $i_\ell \in I_r \cap [1,r-1]$, as in \Cref{fig: i_ell}. Since $i_\ell \in I_r$, then  $w^{-1}(i_\ell)$ must be in $[i_\ell, r-1]$ so that $1 \leq_r i_\ell \leq_r w^{-1}(i_\ell) \leq_{r} r-1$. Hence, $i_\ell$ is also in $I_1$, so $I_r \cap [1,r-1] \subseteq I_1 \cap [1,r-1]$. Thus, these $i_\ell \in I_r \cap [1,r-1]$ appear among the elements $j_g, \ldots, j_k$, with possibly some additional elements. So, each of these $i_{\ell}$ has some $j_{\ell'}$ with $\ell' \leq \ell$ for which $i_{\ell} = j_{\ell'}$. After possibly deleting $a$, all of these elements are in $I$.

For the only if direction in this case, suppose $[a,r-1] \cap (I_1 \setminus I_r)$ is empty. By construction, $I_1 \cap[a,r-1] = \{j_h, \ldots, j_k\}$. Then $I_1 \cap [a,r-1] = I_r \cap [a,r-1]$, so $i_\ell = j_\ell$ for every $\ell \in [h,k]$. In particular, $i_h = j_h = a$. But then we have $i_h = a >_r j_{h-1}$, which violates Condition (iii) for $\ell = h-1$.

Conversely, suppose $(I_1 \setminus I_r) \cap [a,r-1] $ is nonempty. Then $|I_1 \cap [a,r-1]| > |I_r \cap [a,r-1]|$.  Since $I_1 \cap [a,r-1] = \{j_h, \ldots, j_k\}$, then $i_h$ cannot be in $[a,r-1]$, so $i_h <_r a$, and thus $i_\ell <_r a$ for all $\ell \in [g+1,h]$. Furthermore, the existence of at least one extra element in $I_1 \cap [a,r-1]$ implies that the elements of $I_r \cap [1,a-1]$ all appear shifted to the left in $I_1$. Specifically, for $i_{\ell} \in \{i_{g+1}, \ldots, i_h\} \cap [1,a-1]$, there is some $\ell'<\ell$ for which $i_{\ell} = j_{\ell'}$. Then, for these $i_{\ell}$, we have $i_{\ell} = j_{\ell'} \leq_r j_{\ell-1}$. The remaining $i_{\ell} \in \{i_{g+1}, \ldots, i_h\}$ are in $[r,n]$, so $i_{\ell} <_r j_{\ell-1} \in [1,a-1]$. Hence, $i_{\ell} \leq_r j_{\ell-1}$ for all $\ell \in [g+1,h]$, which is equivalent to Condition (iii).


\begin{figure*}[h!]
    \centering
  \subcaptionbox{$j_\ell \in I_1 \cap [r,n] \Rightarrow j_\ell \in I_r$. \label{fig: j_ell}}{\includegraphics[width=1.7in]{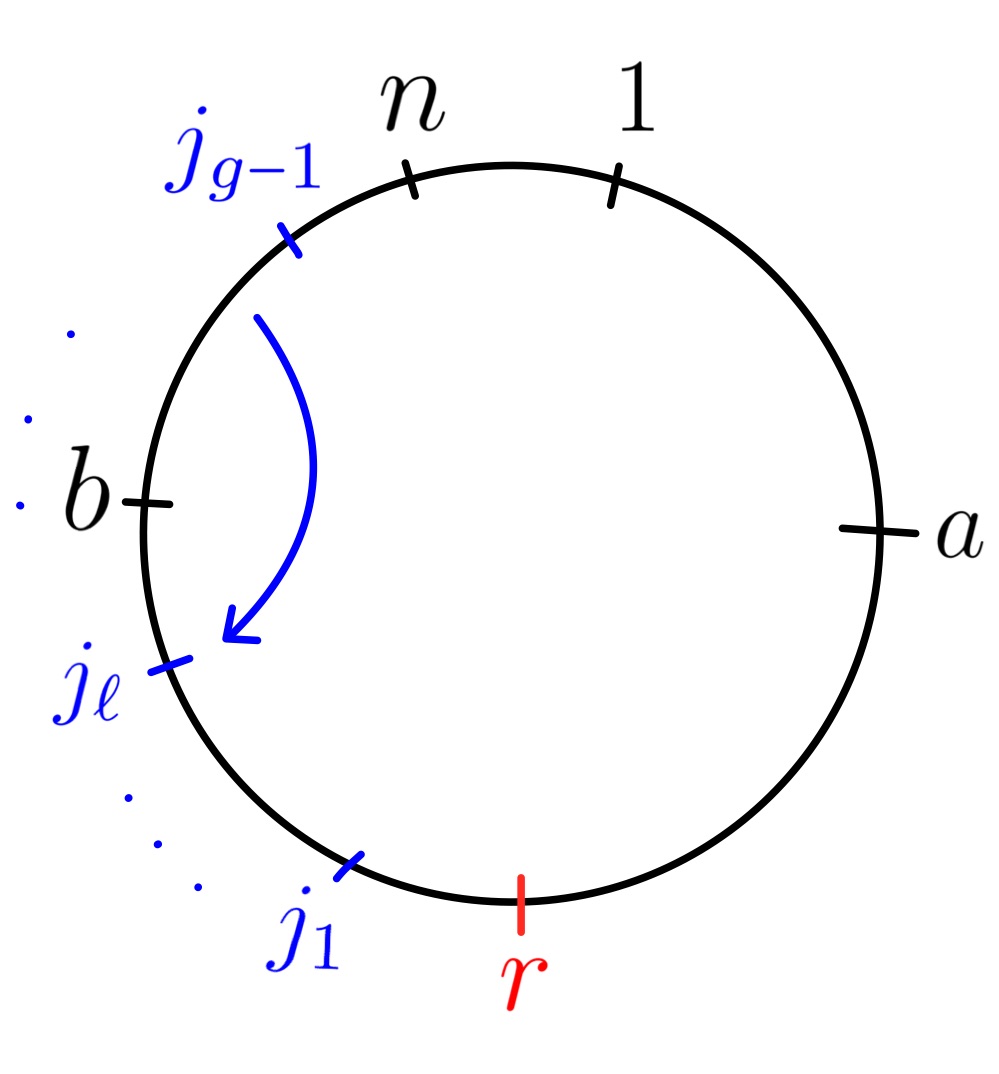}}\hspace{4em}%
  \subcaptionbox{$i_\ell \in I_r \cap [1,r-1] \Rightarrow i_\ell \in I_1$.\label{fig: i_ell}}{\includegraphics[width=1.8in]{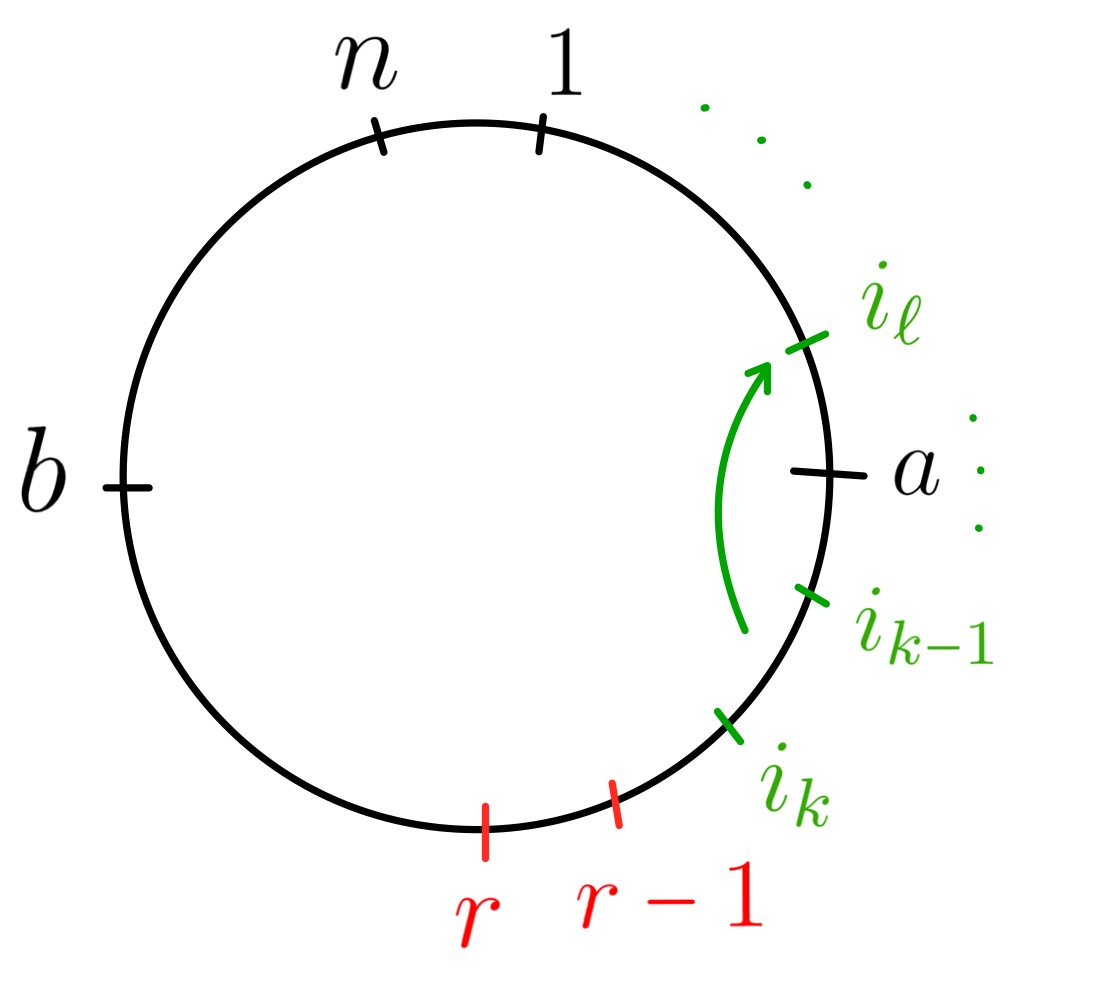}}
    \caption{Elements in both $I_1$ and $I_r$.}
\label{fig: elts in both I_1 and I_r}
\end{figure*}
\end{proof}

\subsection{Mapping Anti-Exchange Pairs to Alignments}\label{sec:strdpairs.to.alignments}
For a fixed $\w \in \SSn$ whose chord diagram contains a crossed alignment, we will show that \eqref{eq: num anti-exch pairs less codim} holds for $J = I_1(\w)$. This will be achieved using an injective map $\Psi_{\w}$, which maps anti-exchange pairs for $I_1(\w)$ to alignments in $D(\w)$ and is defined below for any $\w \in \SSn$ using \Cref{Lemma: anti-exchange pair condition}.


Given $\w \in \SSn$, let $\antiexch(\w)$ be the set of all anti-exchange pairs $(a,b)$ for $I_1(\w)$, so $a \in I_{1}(\w)$, $b \not \in I_{1}(\w)$, and $(I_1(\w) \setminus \{a\}) \cup \{b\} \notin \M(\w)$. For each $(a,b) \in \antiexch(\w)$, either $a>b$ or $a<b$ and there must exist some $r \in [a+1,b]$ for which Condition (1)
or (2) of \Cref{Lemma: anti-exchange pair condition} fails.  In the
$a<b$ case, we say $r$ is a \textit{witness} for $(a,b)$ to be in
$\antiexch(\w)$. Let $\antiexch_{>}(\w)$ be all the
anti-exchange pairs with $a>b$. Let $\antiexch_{1}(\w)$ be all the anti-exchange pairs
with $a<b$ such that Condition (1) fails for some $r \in [a+1,b]$.
Let $\antiexch_{2}(\w)$ be the anti-exchange pairs with $a<b$ such that
Condition (1) holds for all $r \in [a+1,b]$. \Cref{Lemma:
anti-exchange pair condition} then implies for $(a,b) \in \antiexch_{2}(\w)$ that Condition (2) fails
for some $r \in [a+1,b]$. These three sets form a
partition of the anti-exchange pairs for $I_{1}(\w)$,
\begin{equation}\label{eq:anti-exchange.partition}
\antiexch(\w)= \antiexch_{>}(\w) \sqcup \antiexch_{1}(\w) \sqcup
\antiexch_{2}(\w).
\end{equation}

\begin{Lemma} \label{lem: b > bar r} Fix $\w = (w,co) \in \SSn$ and $a \in
I_{1}(\w)$. Let $\{ b_1 < \cdots < b_s \}$ be the set of all elements
in $[a+1,n] \setminus I_1(\w)$ such that for each $1\leq i\leq s$,
there exists some minimal $r_i \in [a+1,b_i]$ for which Condition (1)
of \Cref{Lemma: anti-exchange pair condition} fails, so $(a,b_{i}) \in
\antiexch_{1}(\w)$. If $\{ b_1 < \cdots < b_s\}$ is nonempty, then
defining $\bar r(a,\w) :=r_{1}$ we have
\[
\{b_1 < \cdots < b_s\} \; = \; [\bar r(a,\w), n] \setminus I_1, 
\]
and $\bar r(a,\w) \in [a+1,b_i]$ is a witness for Condition (1) failing for
$(a,b_{i})$ for each $i$.
\end{Lemma}

\begin{proof}
Assume $\{b_1 < \cdots < b_s\}$ is nonempty.  By construction, $r_{1}
\in [a+1,b_1] \subseteq [a+1,b_i]$ for all $i$.  Furthermore, $r_{1}$ is chosen so that
Condition (1) of \Cref{Lemma: anti-exchange pair condition} fails for the pair $(a,b_1)$,
and hence there is no $x \in [a, r_{1}-1]$ such that $w^{-1}(x)\in [r_{1},n]$.
This last condition only depends on $r_{1}$, so Condition (1) still fails
for $r_{1}$ when determining whether or not $a$ is exchangable with
any given $b \in [r_{1},n] \setminus I_1$, according to \Cref{Lemma: anti-exchange pair condition}.  Hence, $\{b_1 < \cdots <
b_s\}= [r_{1},n] \setminus I_1.$
\end{proof}

\begin{Cor} \label{cor: type 2 < bar r} If $\{b_1 < \cdots < b_s\}$
from \Cref{lem: b > bar r} is nonempty, then $b < \bar r(a,\w)$ for
any $(a,b) \in \antiexch_{2}(\w)$.
\end{Cor}

\begin{proof}
For $(a,b) \in \antiexch_{2}(\w)$, Condition (1) of
\Cref{Lemma: anti-exchange pair condition} is satisfied for every $r
\in [a+1,b]$. In particular, $b$ is not one of the $b_i$s as in
\Cref{lem: b > bar r}, so $b$ is not in $[\bar r(a,\w),n] \setminus
I_1$. Since $b \notin I_1$ by definition of an anti-exchange pair, $b$
must be in $[1, \bar r(a,\w)-1]$.
\end{proof}

\begin{Lemma} \label{lem: a < underline r}
Fix  $\w = (w,co) \in \SSn$  and $b \not \in
I_{1}(\w)$.  Let $\{a_1 < \cdots < a_s \}$ be the elements of $[b] \cap
I_1(\w)$ such that there exists some maximal $r_i \in [a_i+1,b]$ for which
Condition (2) of \Cref{Lemma: anti-exchange pair condition} fails, so
$(a_{i},b)\in \antiexch(\w)$. If $\{a_1 < \cdots < a_s \}$ is not
empty, then defining $\underline r(b,\w) :=r_{s}$ we have 
\[
\{a_1 < \cdots < a_s \} \; = \; [1,\underline r(b,\w)-1] \cap I_1(\w), 
\]
and $\underline r(b,\w) \in [a_i+1,b]$ is a witness for Condition (2) failing for
$(a,b_{i})$ for each $i$.
\end{Lemma}

\begin{proof}
The proof is similar to the proof of \Cref{lem: b > bar r} by reflecting $D(\w)$ across the vertical axis. 
\end{proof}

For each $\w \in \SSn$, we construct a map
$\Psi_{\w}$ from the anti-exchange pairs for $I_1(\w)$ paritioned according
to \eqref{eq:anti-exchange.partition} to ordered pairs of
the form $(p \mapsto w(p), s \mapsto w(s))$ for $p \neq s$.  We will
show in the lemma that follows that the range of $\Psi_{\w}$ is
contained in the set of alignments for $\w$, as defined in \Cref{Def: algnmt}.


\begin{Def}\label{alg: psi algorithm}  For  $\w = (w,co) \in \SSn$, let $\Psi_{\w}$ be
defined by the following algorithm.

\begin{algorithm}[H] 
    \SetKwInOut{Input}{Input}
    \SetKwInOut{Output}{Output}
    

    \Input{$(a,b) \in \antiexch(\w)= \antiexch_{>}(\w) \sqcup\antiexch_{1}(\w)\sqcup
    \antiexch_{2}(\w)$}
    \Output{ $(p \mapsto w(p), s \mapsto w(s)) \in
    \aligns(\w)$}
  
\medskip

    \textbf{set}  $p \gets w^{-1}(b)$  and  $s \gets w^{-1}(a)$.

    \If{$(a,b) \in \antiexch_{1}(w)$}
    	{
        \While{$p \not \in [\bar r(a,\w),n]$,}{
         	update $p \leftarrow w^{-1}(p)$.}
	}

    \If{$(a,b) \in \antiexch_{2}(w)$}
    	{
        \While{$s \not \in [1, \underline r(b,\w)-1]$,}{
         	update $s \leftarrow w^{-1}(s)$.
		}
	}

    \Return{$(p \mapsto w(p), s \mapsto w(s))$.}  
\end{algorithm}
\end{Def}


  




      

\bigskip

\begin{Lemma} \label{lem: psi well defined}
For  $\w \in \SSn$ and $(a,b) \in \antiexch(\w)$, the algorithm defined in \Cref{alg: psi algorithm} terminates in finitely many steps, and the image $\Psi_{\w}(a,b)$ is in $\aligns(\w)$.
\end{Lemma}

\begin{proof}

First, consider the case when $(a,b) \in \antiexch_{>}(\w)$, so $a > b$. Since $a \in I_1(\w)$ and $b \notin I_1(\w)$, we have the inequalities $w^{-1}(b) \leq b < a \leq w^{-1}(a)$, as indicated in \Cref{Fig: stradd pair 1}. Thus, $\Psi_{\w}(a,b) = (w^{-1}(b) \mapsto b, w^{-1}(a) \mapsto a)$ is an alignment. Note that $a$ and $b$ could be fixed points with $co(a) = \cw$ and $co(b) = \ccw$, according to the definition of anti-exceedances from \Cref{Sec: decperms}.

\begin{minipage}[t]{1.0\linewidth}
    \centering
   \vspace{-2ex}
   \hspace*{-2cm}
   \includegraphics[scale=0.10]{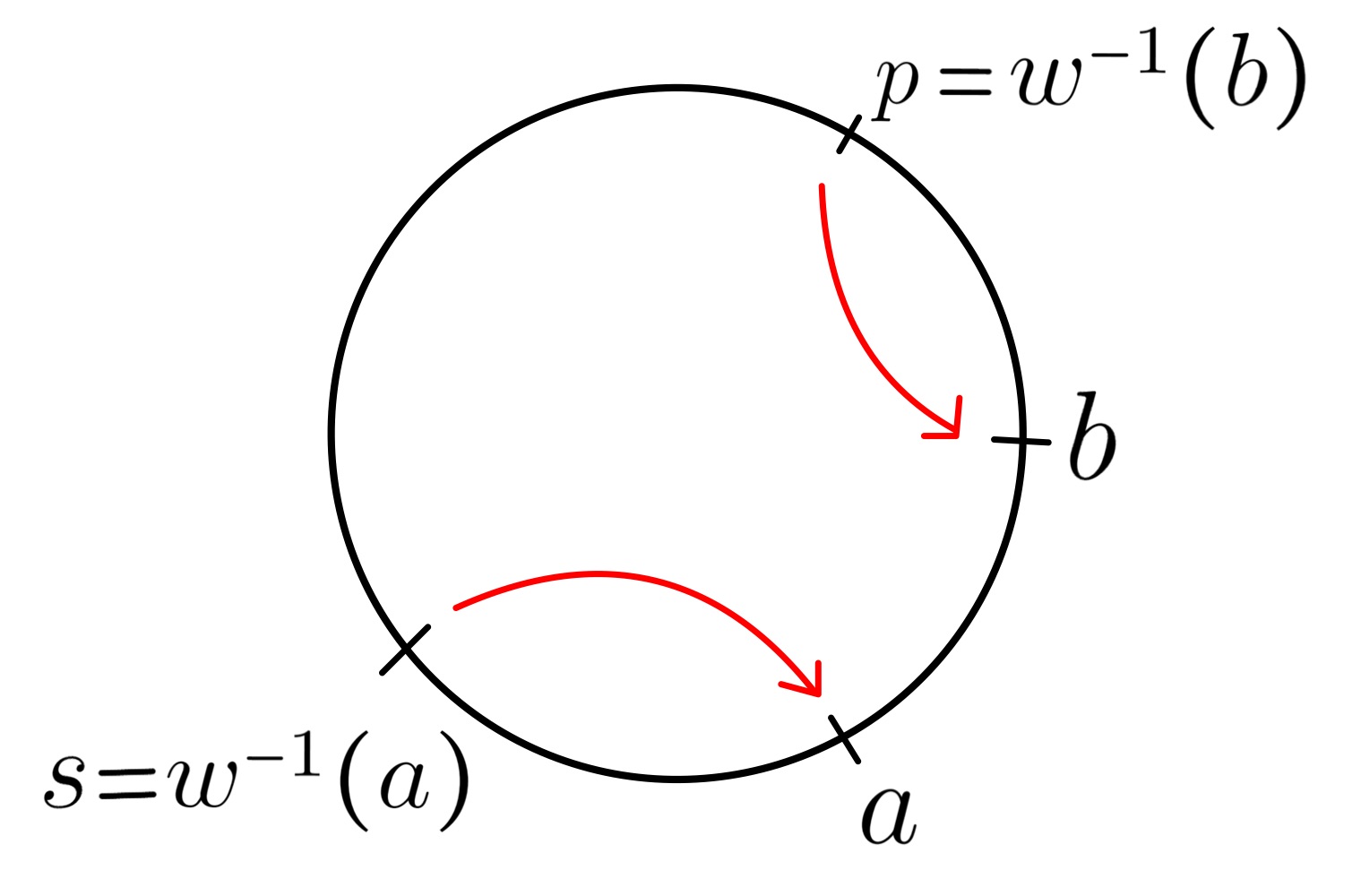}   \captionof{figure}{The image of $(a,b)$ when $a>b$ is $(w^{-1}(a) \mapsto a, w^{-1}(b) \mapsto b)$.}
   \label{Fig: stradd pair 1}
\end{minipage}

Next, consider the case when $(a,b) \in \antiexch_{1}(\w)$, and let $\bar
r = \bar r(a,\w)$. Then $\Psi_{\w}(a,b) = (p \mapsto w(p), w^{-1}(a) \mapsto a)$, where $p$ is in $[\bar r,n]$. By definition of $\antiexch_{1}(\w)$, Condition (1) of \Cref{Lemma: anti-exchange pair
condition} fails for some $r \in [a+1,b]$. It follows from
\Cref{lem: b > bar r} that $b \geq \bar r$, and furthermore, that Condition (1) fails for this $\bar r$. In particular, there is no $x \in [a,\bar r
-1]$ with preimage in $[\bar r,n]$.  Since $a \in I_1(\w)$, then
$w^{-1}(a) \in [a,n]$. But $a \in [a, \bar r-1]$ implies that
$w^{-1}(a) \notin [\bar r,n]$. Hence, $w^{-1}(a) \in [a, \bar r -1]$, as drawn in \Cref{fig: case1}.

The termination of the while loop in Line 3 of \Cref{alg: psi algorithm} comes from the fact that $p$ is obtained by tracing in
reverse the cycle containing $b$. By \Cref{lem: b > bar r}, $b \in
[\bar r,n]$. If $w^{-1}(b) \in [\bar r,n]$, then $p = w^{-1}(b)$, and the loop
ends immediately. In this case, since $b$ is an exceedance by the
assumption that $(a,b) \in \antiexch_{1}(\w)$, then $w^{-1}(b) \in [\bar r, b]$. See \Cref{fig:
case1a}. However, even if $w^{-1}(b)$ is not in $[\bar r,n]$ so that
the reverse cycle immediately leaves the interval $[\bar r,n]$, it
must eventually return to $[\bar r,n]$, since $b \in [\bar r,n]$. In
this case, the first element at which the reverse cycle returns to
$[\bar r,n]$ is $p$, which implies that $w(p) \in [1, \bar r
-1]$. However, the fact that Condition (1) is not satisfied for $\bar
r$ then implies that $w(p) \notin [a, \bar r -1]$, and thus $w(p) \in
[1,a-1]$. See \Cref{fig: case1b}. In both cases, $(p \mapsto
w(p), w^{-1}(a) \mapsto a)$ is an alignment with port side $(p \mapsto w(p))$ and starboard side $(w^{-1}(a) \mapsto a)$.

\begin{figure*}[h!]
    \centering
  \subcaptionbox{$w^{-1}(b) \geq \bar r$\label{fig: case1a}}{\includegraphics[width=1.8in]{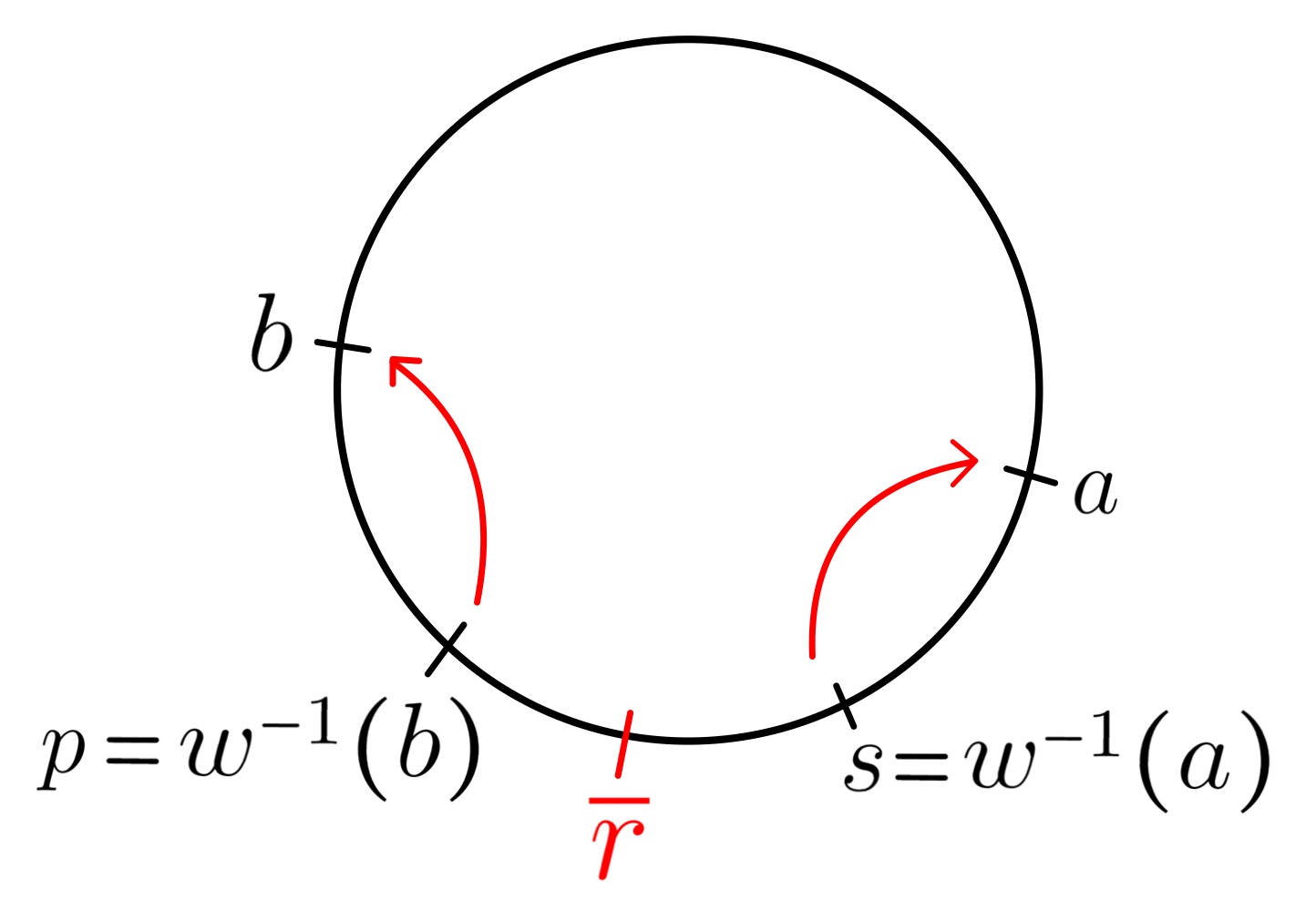}}\hspace{3em}%
  \subcaptionbox{$w^{-1}(b) < \bar r$\label{fig: case1b}}{\includegraphics[width=1.6in]{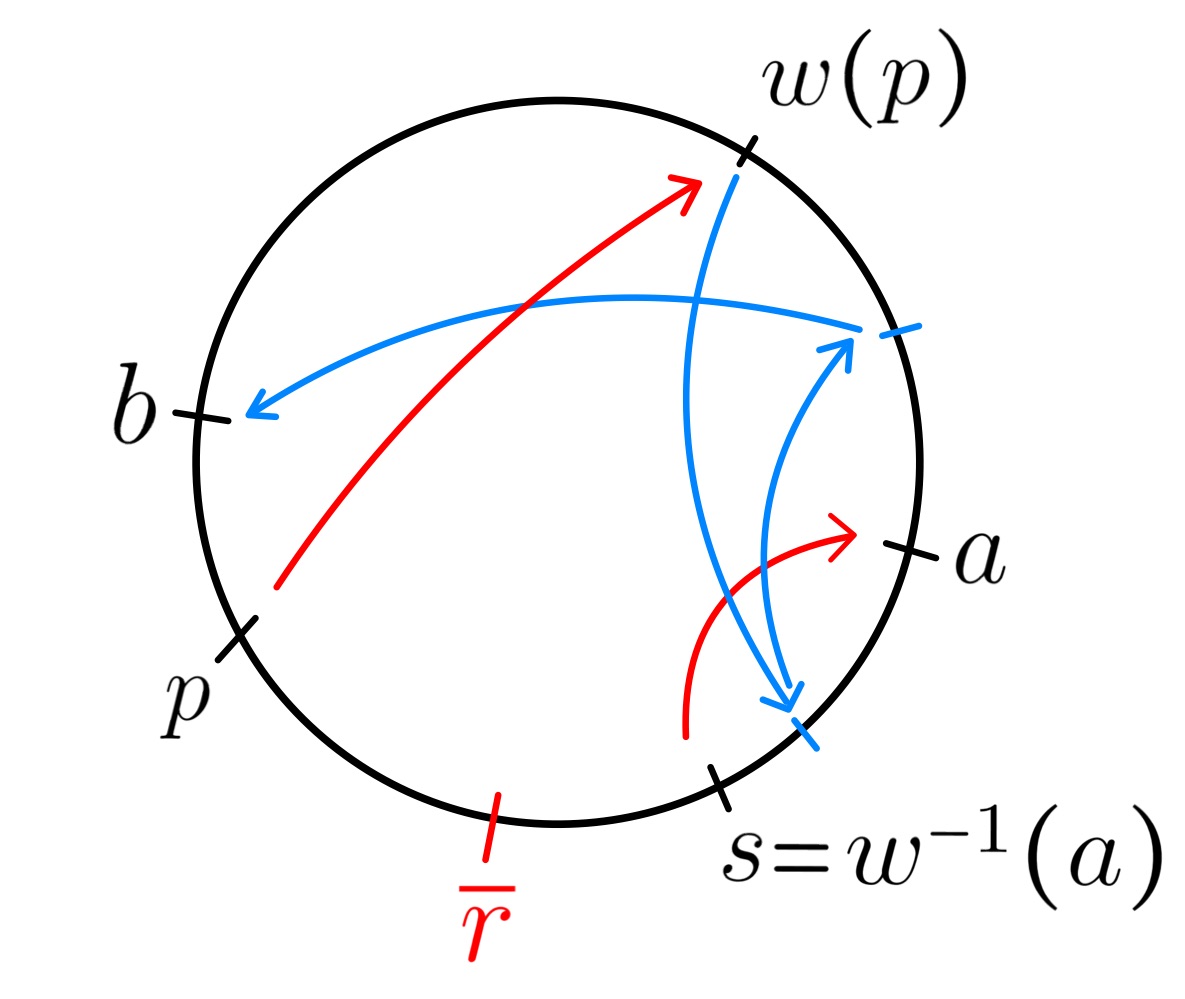}}
    \caption{The red alignment is the image of $(a,b) \in \antiexch_{1}$.}
\label{fig: case1}
\end{figure*}

Finally, consider the case where $(a,b) \in \antiexch_{2}(\w)$, and let $\underline r = \underline r(b,\w)$. Again, this case is symmetric to the $\antiexch_{1}(\w)$ case above. From the algorithm, $\Psi_{\w}(a,b) = (w^{-1}(b) \mapsto b, s \mapsto w(s))$, where $s \in [1, \underline r -1]$. By \Cref{lem: a < underline r}, $a < \underline r$, and Condition (2) is not satisfied for $\underline r$. Thus, no $y \in [\underline r,b]$ has preimage in $[1, \underline r-1]$. Since $b \notin I_1(\w)$, we must have $w^{-1}(b) \in [\underline r, b]$, as drawn in \Cref{fig: case2}.

The termination of the while loop in Line 8 of \Cref{alg: psi
algorithm} is again seen to terminate by tracing the cycle containing $a$ in reverse.  Since $a \in I_1(\w)$, then $a \leq w^{-1}(a)$. If $w^{-1}(a) \in [a, \underline r-1]$, then $s = w^{-1}(a)$, as in \Cref{fig: case2a}. If $w^{-1}(a) \in [\underline r, n]$, then $s$ is the first element at which the reverse cycle returns to the interval $[1, \underline r-1]$, as in \Cref{fig: case2b}. Thus $w(s) \in [\underline r,n]$, and the fact that Condition (2) of \Cref{Lemma: anti-exchange pair condition} fails for $\underline r$ implies that $w(s) \in [b+1,n]$. So, in both cases, $(w^{-1}(b) \mapsto b, s \mapsto w(s))$ is an alignment with port side $(w^{-1}(b) \mapsto b)$ and starboard side $(s \mapsto w(s))$.

\begin{figure*}[h!]
    \centering
  \subcaptionbox{$w^{-1}(a) < \underline r$\label{fig: case2a}}{\includegraphics[width=1.8in]{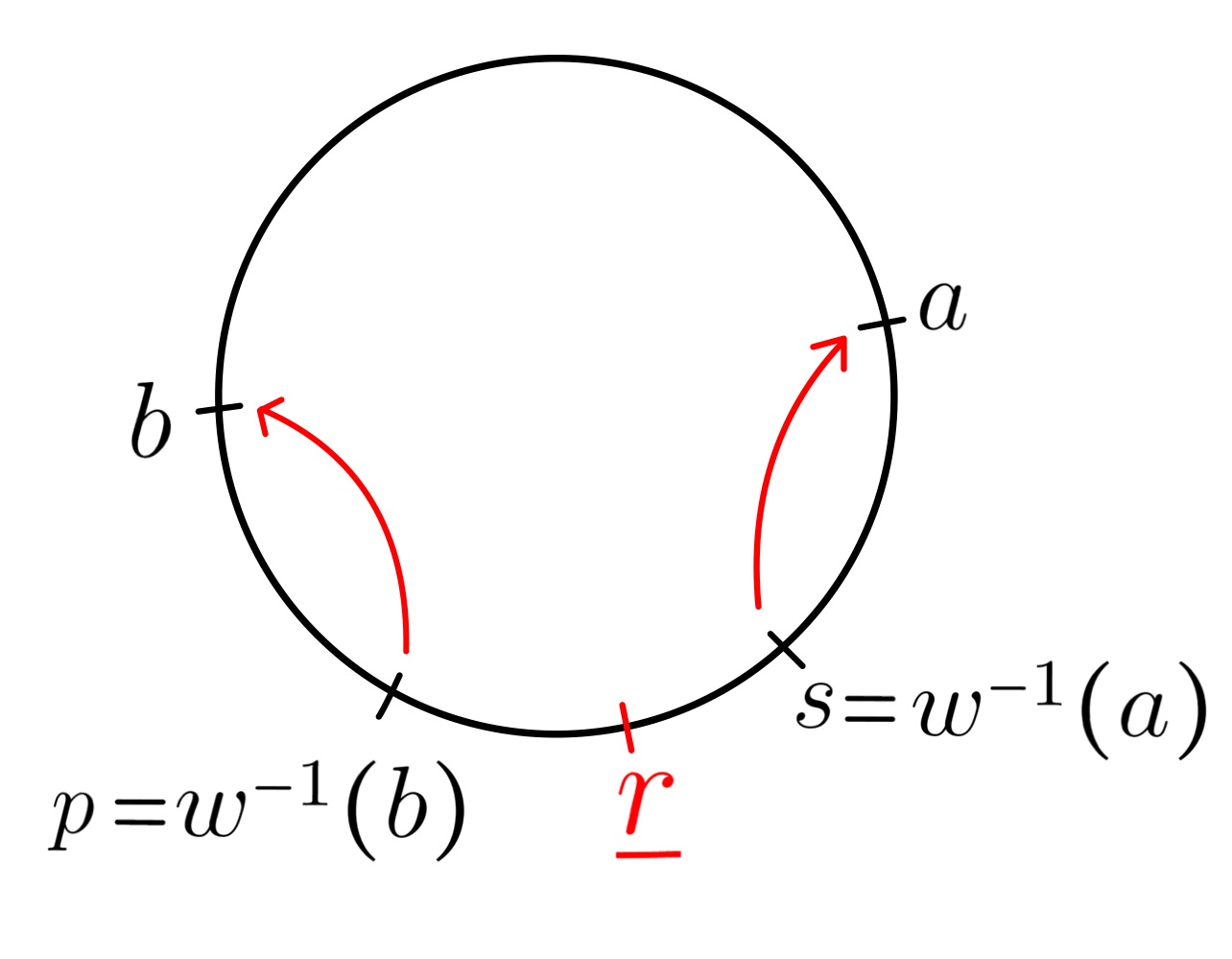}}\hspace{3em}%
  \subcaptionbox{$w^{-1}(a) \geq \underline r$\label{fig: case2b}}{\includegraphics[width=1.6in]{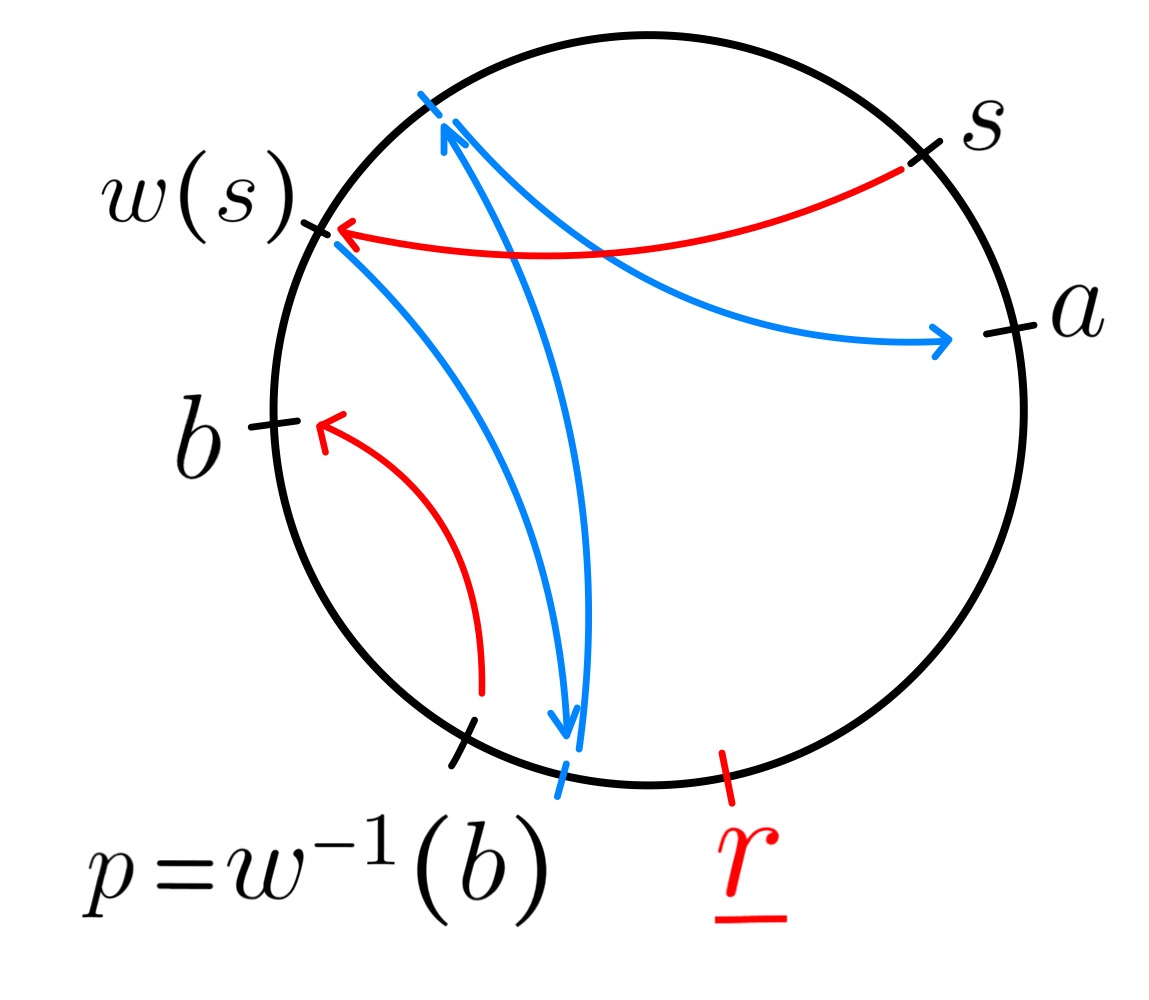}}
    \caption{The red alignment is the image of $(a,b) \in \antiexch_{2}$.}
\label{fig: case2}
\end{figure*}
\end{proof}


\begin{Lemma} \label{Cor: all distinct alignments}
Fix $\w = (w,co) \in \SSn$. The map $\Psi_{\w}: \antiexch(\w) \longrightarrow \aligns(\w)$ is injective.
\end{Lemma}

\begin{proof}
Let $(a_1,b_1)$ and $(a_2,b_2)$ be two distinct anti-exchange pairs
for $I_1(\w)$. We need to show $\Psi_{\w}(a_1,b_1) \neq \Psi_{\w}(a_2,b_2)$.
We again utilize the partition of $\antiexch(\w)$ from
\eqref{eq:anti-exchange.partition} and consider several cases.

\noindent \textbf{Case 1}: Assume at least one of $(a_1,b_1)$
and $(a_2,b_2)$ is in $\antiexch_{>}(\w)$. Observe from \Cref{alg:
psi algorithm} that anti-exchange pairs $(a,b) \in \antiexch_{>}(\w)$
are the only pairs assigned alignments $\Psi_{\w}(a,b) = (p
\mapsto w(p),s \mapsto w(s))$ with $p \leq w(p) < s \leq w(s)$. Thus, $\Psi_{\w}(a_1,b_1) \neq \Psi_{\w}(a_2,b_2)$ unless
both $(a_1,b_1)$ and $(a_2,b_2)$ are in $\antiexch_{>}(\w)$.  In this case, $\Psi_{\w}(a_i,b_i) = (w^{-1}(b_i) \mapsto b_i,
w^{-1}(a_i) \mapsto a_i)$ for both $i$ by definition.  Therefore, since $(a_1,b_1)$ and $(a_2,b_2)$ are distinct, we must have that
$\Psi_{\w}(a_1,b_1) \neq \Psi_{\w}(a_2,b_2)$.

\noindent \textbf{Case 2}: Assume both $(a_1,b_1)$ and $(a_2,b_2)$ are in $\antiexch_{1}(\w)$ and $\Psi_{\w}(a_i,b_i) = (p_i \mapsto
w(p_i), w^{-1}(a_i) \mapsto a_i)$. If $a_{1}\neq a_{2}$, then the two
alignments have distinct starboard sides. If  $a_{1}=a_{2}=a$, 
then $b_1 \neq b_2$.  Since the $p_i$ are determined uniquely by tracing the cycle containing $b_i$ backwards and finding the first element in the interval $[\bar
r(a,w),n]$, then $b_1 \neq b_2$ implies that $p_1 \neq p_2$, so the alignments have distinct port sides.  Either
way, we have $\Psi_{\w}(a_1,b_1) \neq \Psi_{\w}(a_2,b_2)$.

 \noindent \textbf{Case 3}: Assume both $(a_1,b_1)$ and $(a_2,b_2)$ are in $\antiexch_{2}(\w)$ and $\Psi_{\w}(a_i,b_i) = (w^{-1}(b_i) \mapsto
b_i, s_i \mapsto w(s_i))$.  If $b_{1}\neq b_{2}$, then the two
alignments have distinct port sides. If $b_{1}=b_{2}=b$, then $a_1
\neq a_2$.  Since the $s_i$ are determined uniquely by tracing the
cycle containing $a_i$ backwards and finding the first element in the
interval $[1, \underline r(b,\w) -1]$,
then $a_1 \neq a_2$ implies $s_1 \neq s_2$.  Either way,
$\Psi_{\w}(a_1,b_1) \neq \Psi_{\w}(a_2,b_2)$.

 \noindent \textbf{Case 4}: Assume one of the $(a_i,b_i)$
is in $\antiexch_{1}(\w)$ and the other anti-exchange pair is in $\antiexch_{2}(\w)$.
Without loss of generality, we may assume that $(a_1,b_1) \in
\antiexch_{1}(\w)$ with $\bar r(a_{1},\w) \in [a_{1}+1,b_{1}]$ and
$(a_2,b_2) \in \antiexch_{2}(\w)$ with $\underline r(b_2,\w) \in
[a_2+1,b_2]$.  Then $\Psi_{\w}(a_{1},b_1) = (p \mapsto w(p),
w^{-1}(a_{1}) \mapsto a_{1})$ for some $p \geq \bar
r(a_{1},\w)$, and $\Psi_{\w}(a_{2},b_2) = (w^{-1}(b_2) \mapsto b_2,
s \mapsto w(s))$ for some $s < \underline
r(b_{2},\w)$.

If $a_{1}= a_{2}=a$, then \Cref{cor: type 2 < bar r} shows that $b_2 <
\bar r(a,\w)$ since $(a,b_2) \in \antiexch_{2}(\w)$ and
$(a,b_1) \in \antiexch_{1}(\w)$. Since $b_2 \notin I_1(\w)$, $w^{-1}(b_2) \leq b_2 < \bar r(a,\w) \leq
p$. Therefore, once again, the two alignments cannot have the same
port sides, so $\Psi_{\w}(a_{1},b_1) \neq \Psi_{\w}(a_{2},b_2)$.

If $a_{1}\neq a_{2}$, observe from \Cref{alg: psi algorithm} that $(s
\mapsto w(s))$ is either the arc $(w^{-1}(a_2) \mapsto a_2)$, or it is
an arc with $s < \underline r(b_2,w)$ and $w(s) \geq
\underline r(b_2,\w)$,
in which case it is an exceedance arc. See \Cref{fig: case2} for these two cases. In either case, $(s \mapsto w(s))$ is
not the anti-exceedance arc $(w^{-1}(a_1) \mapsto a_1)$. Hence, $\Psi_{\w}(a_1,b_1) \neq
\Psi_{\w}(a_2,b_2)$.
\end{proof}

\subsection{The Case of Crossed Alignments}

\begin{figure*}[h!]
    \centering
  \subcaptionbox{$d \neq 1$\label{fig: crossed d neq 1}}{\includegraphics[width=1.7in]{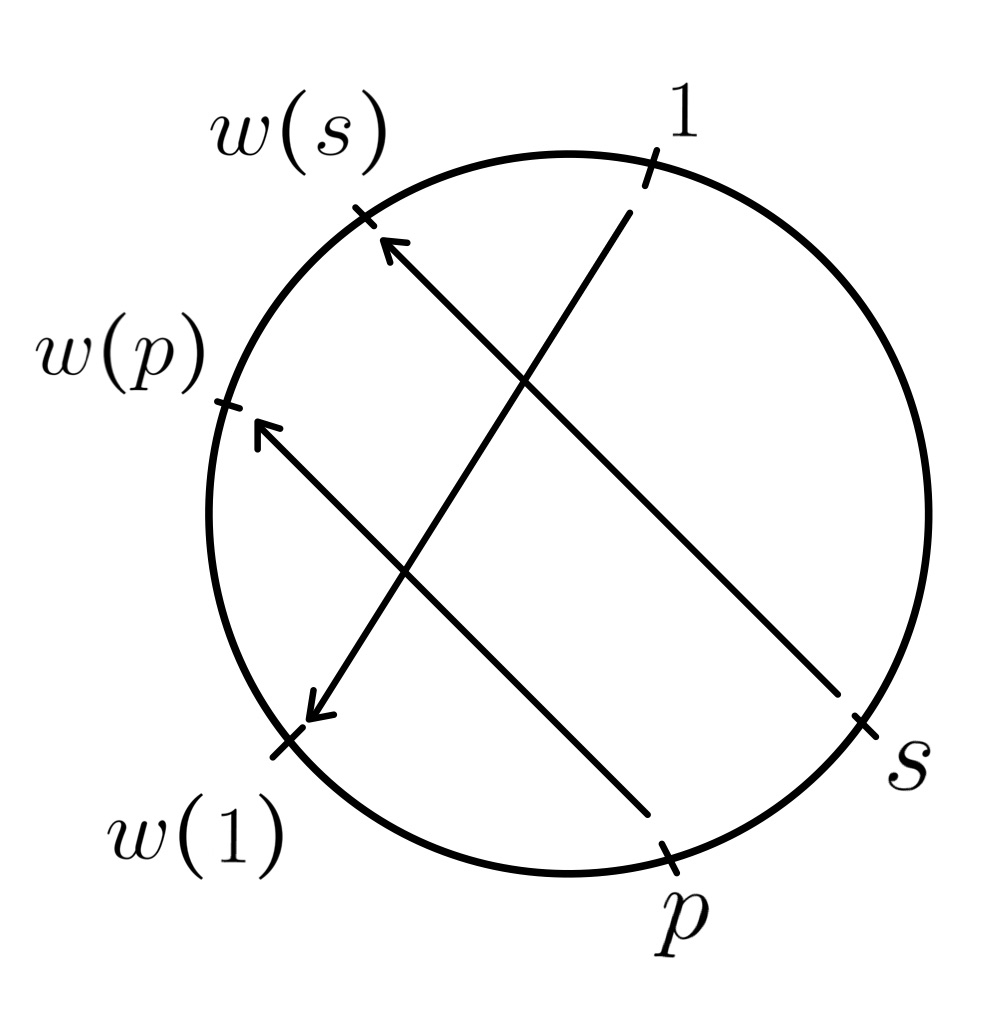}}\hspace{3em}%
  \subcaptionbox{$d = 1$\label{fig: crossed d=1}}{\includegraphics[width=1.8in]{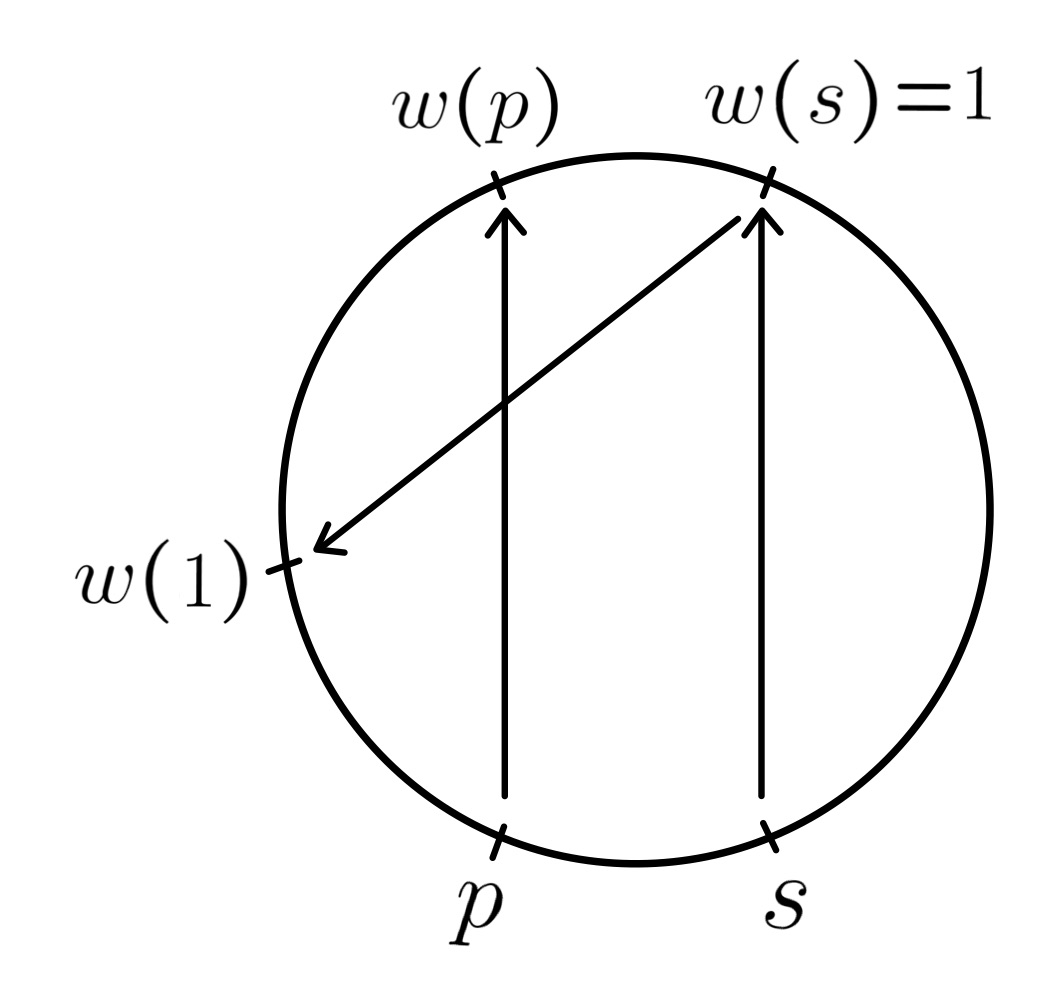}}
    \caption{Two cases for starboard tacking crossing arc with tail at 1.}
\label{fig: crossing arc}
\end{figure*}

Recall the definition of a starboard tacking crossed alignment from \Cref{Def: crossed algnmt}. \Cref{fig: crossing arc} depicts two examples of a starboard tacking crossed alignment where the tail of the crossing arc is 1.

\begin{Lemma} \label{Lemma: strict for crossed almt}
If $\w  = (w,co) \in \SSn$ has a starboard tacking crossed alignment $(p \mapsto
w(p), s \mapsto w(s))$, where 1 is the tail of the crossing arc, then there
is no element of $\antiexch(\w)$ mapping to $(p \mapsto w(p), s \mapsto
w(s))$ under $\Psi_{\w}$.  Hence, $\Psi_{\w}$ is not surjective.
\end{Lemma}


\begin{proof}
Assume by way of contradiction that there is some anti-exchange pair $(a,b) \in \antiexch(\w)$ such that $\Psi_{\w}(a,b) = (p \mapsto w(p), s \mapsto w(s))$. Since fixed point loops cannot cross other arcs, then none of the arcs involved in a crossed alignment are loops. Therefore, none of $1$, $p$, or $s$ is a fixed point. There are two cases to consider, depending on whether $w(s) = 1$. 

If $w(s) \neq 1$, as in \Cref{fig: crossed d neq 1}, then $w(p)$ and $w(s)$ are both exceedances. By considering the various cases in the algorithm defined in \Cref{alg: psi algorithm}, $w(p)$ and $w(s)$ can both be exceedances only if $a < b$ and $(a,b) \in \antiexch_{2}(\w)$, as in \Cref{fig: case2b}. In this case, we must have that $w(p) = b$. Hence, the arc $(1 \mapsto w(1))$ satisfies Condition (2) of \Cref{Lemma: anti-exchange pair condition} for all $r \in [2,w(1)]$, and the arc $(p \mapsto w(p) = b)$ satisfies Condition (2) for $r \in [p+1,w(p)]$. Since the alignment is crossed by $(1 \mapsto w(1))$, then $w(1) \geq p$. It follows that Condition (2) is satisfied for all $r \in [2,w(p)]$ and thus for all $r$ in $[a+1,w(p)]$. This contradicts that $(a,b) \in \antiexch_{2}(\w)$.


If $w(s)=1$, as in \Cref{fig: crossed d=1}, then the arrangement of
the arcs in $(p \mapsto w(p), s \mapsto w(s))$ implies that $a <
b$. Thus, $(a,b) \in \antiexch_{1}(\w) \sqcup \antiexch_{2}(\w)$. If
$(a,b) \in \antiexch_{2}(\w)$, then we must have $w(p)=b$. Again,
looking at the arcs $(1 \mapsto w(1))$ and $(p \mapsto w(p))$
contradicts the assumption that Condition (2) fails for some $r \in
[a+1,w(p)]$. Otherwise, $(a,b) \in \antiexch_{1}(\w)$, and the
starboard side being $(s \mapsto w(s))$ implies that $a=1$. The port
side being an exceedance implies that $w(p)=b$, as in \Cref{fig:
case1a}. Thus, Condition (1) fails for some $r \in [2,w(p)]$. If $r
\in [2, w(1)]$, let $x=1$, and if $r \in [w(1)+1,w(p)]$, let $x=p$. In
either case, the arc $(x \mapsto w(x))$ maps from $[1,r-1]$ to
$[r,n]$. The cycle containing the arc $(x \mapsto w(x))$ must then
also contain an arc $(y \mapsto w(y))$ from $[r,n]$ to $[1,r-1]$. The
existance of the arc $(y \mapsto w(y))$ contradicts the fact that Condition (1) fails for $r$. Therefore, $(a,b) \notin \antiexch_{1}(\w)$. 

For these two cases of $w(s) \neq 1$ and $w(s) = 1$, we have obtained contradictions that $(a,b)$ is an anti-exchange pair. Therefore, there can be no such anti-exchange pair mapping to $(p \mapsto w(p), s \mapsto w(s))$ under $\Psi_{\w}$. 
\end{proof}


\begin{proof}[Proof of \Cref{The: Main theorem}, $(1) \Rightarrow (5)$]
We prove the implication $(1) \Rightarrow (5)$ using the contrapositive. Thus, assume that $\w$ is a decorated permutation whose chord diagram contains a crossed alignment. The chord diagram $D(\w)$ may be rotated and reflected as necessary to produce a new chord diagram with a starboard tacking crossed alignment whose crossing arc has its tail at 1. By
\Cref{Lemma: invariances}, any such operations on $D(\w)$ preserve the
property of being smooth or singular. Therefore, we may assume that $D(\w)$ has a starboard tacking crossed alignment whose crossing arc has its tail at 1. It follows from \Cref{Lemma: strict for
crossed almt} that $\# \antiexch(\w) < \# \aligns(\w)$. Since $I_1(w) \in
\M$, then \Cref{cor:anti-exchangable.pairs.alignments} and
\Cref{Def: anti-exchange pair} of anti-exchange pairs imply that $\Pi_{\w}$ is
singular.
\end{proof}

\section{Enumeration of Smooth Positroids}\label{sec:enumeration}

Let $\w$ be a decorated permutation in $\SSn$ with associated
positroid $\M$ of rank $k$ on ground set $[n]$.  We say $\w$ is a
\textit{smooth decorated permutation} and $\M$ is a \textit{smooth
positroid} if $\Pi_{\w}=\Pi_{\M}$ is a smooth positroid variety.  
\begin{Def} \label{def: num smooth}
Let $s(n)$ be the number of smooth positroids on ground set $[n]$, and let $s_c(n)$ be the number of connected smooth positroids on ground set $[n]$.
\end{Def}
By \Cref{The: Main theorem}, each connected smooth
positroid can be bijectively associated with a spirograph
permutation. Thus, by \eqref{eq:spirgraph.enumeration}, we have 
\begin{equation} \label{eq:num.spirographs}
s_c(n) = \begin{cases}
2& n=1  \\
n-1& n>1
\end{cases}.
\end{equation}

In
analogy with the enumeration of smooth Schubert varieties studied by
Haiman, Bona, Bousquet-M\'elou and Butler
\cite{bona98,BousquetMelou-Butler,Haiman}, the enumeration of smooth
positroids on $[n]$ gives rise to several interesting sequences. See
\cite[A349413, A349456, A349457, A349458, A353131, A353132]{OEIS}.  The following formula
for enumerating smooth positroids is very similar to the results in
\cite[Thm. 10.2]{Ardila-Rincon-Williams}.  In particular, we use the
results of Beissinger \cite{Beissinger1985} and Speicher
\cite{Speicher1994} for counting structures induced on the blocks of
non-crossing partitions and the Lagrange inversion formula for formal
power series \cite[Sect. 5.4]{ec2}.  The sequence begins $2, 5, 16,
61, 256, 1132, 5174, 24229, 115654$ for $n=1,\ldots,10$
\cite[A349458]{OEIS}. If $G$ is a polynomial or power series in $x$, then $[x^n]G(x)$ denotes the coefficient of $x^n$ in $G(x)$.


\begin{Th}\label{lem:smooth.enum.coeff}
The number of smooth positroids on ground set $[n]$ is the coefficient
\begin{equation}\label{eq:smooth.en}
s(n) = [x^n]\frac{1}{n+1}\left(1+2x+\sum_{i=2}^{\infty}(i-1)x^{i}\right)^{n+1} = [x^n]\frac{1}{n+1}\left(1+2x+\sum_{i=2}^{n}(i-1)x^{i}\right)^{n+1}.
\end{equation}
\end{Th}

\begin{proof}
Let $\M= \M_{1} \oplus \cdots \oplus\M_{t}$ be a positroid decomposed
into its connected components on a non-crossing partition.  Then,
$\Pi_{\M}$ is smooth if and only if $\Pi_{\M_{i}}$ is smooth for each
$i \in [t]$ by \Cref{Cor: concat smooth/sing}.  
Each connected smooth
positroid can be bijectively associated with a spirograph permutation
by \Cref{The: Main theorem}.  Therefore, every smooth
positroid on $[n]$ can be uniquely determined by a non-crossing
partition of $[n]$ along with a spirograph permutation on each block.
Hence, the number of smooth positroids on ground set $[n]$ is given by
\begin{equation}\label{eq:sum.over.nc}
s(n)= \sum s_c(\#B_1)\,s_c(\#B_2)\, \cdots \, s_c(\#B_t)
\end{equation}
where the sum is over all non-crossing partitions $B_{1} \sqcup \cdots
\sqcup B_{t}$ of $[n]$. Equation \eqref{eq:smooth.en} now follows directly from
\cite[Corollary 0]{Speicher1994} and \eqref{eq:num.spirographs}.
\end{proof}

The \textit{partial Bell polynomial}, $B_{n,k}(x_1, \ldots, x_{n-k+1})$ introduced in \cite{Bell}, is defined as
\begin{equation} \label{eq:bell.poly}
B_{n,k}(x_1, \ldots, x_{n-k+1}) = \sum_{B_1 \sqcup \cdots \sqcup B_k} \prod_{i=1}^k x_{|B_i|},
\end{equation}
where the sum is taken over all set partitions of $[n]$ into $k$
blocks. The following formula of Fa{\'a} di Bruno expresses the $n$th
derivative of a composition of functions in terms of the partial Bell
polynomials \cite{Johnson.on.Faa} 
\begin{equation} \label{eq:fa.di.bruno}
\frac{d^n}{d x^n} f(g(x)) = \sum_{k=1}^n f^{(k)}(g(x)) \cdot B_{n,k}(g'(x), g''(x), \ldots, g^{(n+1-k)}(x)).
\end{equation}
Define the triangle of numbers $b_{n,k}$ using
the Bell polynomials evaluated at $x_{i}=s_{c}(i)\cdot i!$ for all
$i$, so 
\begin{equation} \label{eq:bnk}
b_{n,k}  = B_{n,k}(2\cdot 1!, 1 \cdot 2!, 2 \cdot 3!, \ldots,
(n-k)\cdot (n-k+1)!) = \sum_{i=1}^{n-k+1}\binom{n-i}{i-1} \cdot s_{c}(i)\cdot i! \cdot b_{n-i,k-1}
\end{equation}
for $1\leq k\leq n$ along with initial condition $b_{0,0}=1$, and
$b_{0,k}=b_{n,0}=0$ if $n>0$ or $k>0$, see \cite[A353131]{OEIS}.
Also, let $(n)_k$ denote the falling factorial,
\[
(n)_k \; := \; n \, (n-1)\,  \cdots \, (n-k+1).
\]


\begin{Cor} \label{cor:smooth.enum.bell}
The number of smooth positroids on ground set $[n]$ is
\begin{equation} \label{eq:smooth.enum.bell}
s(n) = \frac{1}{(n+1)!}\sum_{k=1}^n (n+1)_k \cdot b_{n,k} = \sum_{k=1}^n \frac{b_{n,k}}{(n-k+1)!}.
\end{equation}
\end{Cor}

\begin{proof}
Set $f(x) = x^{n+1}$ and $g(x) = 1 + 2x + \sum_{i=2}^n (i-1)x^i$. By \Cref{lem:smooth.enum.coeff}, $s(n) = \frac{1}{n+1} \cdot [x^n] f(g(x))$. Furthermore, the coefficient of $x^n$ in the composition $f(g(x))$ can be computed as $[x^n] f(g(x)) = \frac{1}{n!} \frac{d^n}{dx^n} f(g(x))\Big|_{x=0}$. Then, by Fa{\'a} di Bruno's formula \eqref{eq:fa.di.bruno},
\begin{align} 
s(n) &= \frac{1}{n+1} \cdot \frac{1}{n!} \cdot \frac{d^n}{dx^n} f(g(x))\Big|_{x=0} \\
  &= \frac{1}{(n+1)!} \sum_{k=1}^n f^{(k)}(g(x)) \cdot B_{n,k}(g'(x), g''(x), \ldots, g^{(n+1-k)}(x)) \Big|_{x=0}. \label{eq:bell}
\end{align}

The $k$th derivative of $f(x) = x^{n+1}$ is $f^{(k)}(x) = (n+1)_k \cdot x^{n-k+1}$, so
\[
f^{(k)}(g(x)) \Big|_{x=0} = (n+1)_k \Big(1+2x + \sum_{i = 1}^n (i-1) x^i \Big)^{n-k+1} \Big|_{x=0} = (n+1)_k.
\]
The derivatives of $g$ are $g'(x) = 2 + \sum_{i=1}^{n-1}i (i+1) x^i$
and $g^{(j)}(x) = \sum_{i=0}^{n-j} (i+j-1) \cdot (i+j)_j x^i$ for
$j>1$. Hence, $g^{(j)}(x) |_{x=0} = s_c(j) \cdot j!$ for all $j \geq
1$. Therefore, the formula in \eqref{eq:smooth.enum.bell} follows from
\eqref{eq:bnk} and \eqref{eq:bell}.
\end{proof}

Note that the value $k$ in the derivation above is different from the value $k = k(\w)$ used elsewhere in the paper. These distinct values may be used to refine the enumeration of $s(n)$. To this end, we define
\begin{align} \label{def:refined.enum}
s_1(n,k) &:= \# \{\text{smooth positroid varieties in } \Gkn\} \hspace{0.05in} \text{ for $0 \leq k \leq n$}\\
    &= \# \{\text{smooth decorated permutations in $\SSnk$}\}, \\
s_2(n,k) &:= \# \left\{
     \begin{aligned}
     &\text{smooth decorated permutations in $\SSn$ with} \\
     &\text{exactly $k$ components in its SIF decomposition}
     \end{aligned}
  \right\} \hspace{0.0in} \text{for $1 \leq k \leq n$}, \text{ and}\\
s_3(n,k) &:= \frac{b_{n,k}}{(n-k+1)!}  \hspace{0.05in} \text{ for $1 \leq k \leq n$.}
\end{align}
The terms $s_i(n,k)$ are displayed in \Cref{Fig:s1} and \Cref{Fig:s2} for $1 \leq n \leq 10$.

For each $i \in \{1,2,3\}$, we may define a $q$-analog of $s(n)$ as
\begin{equation}
\label{eq:q.analogs}
s^{(i)}(n;q) = \sum_{k \leq n} s_i(n,k) q^k.
\end{equation}
In fact, the $q$-analogs for $i=2$ and $i=3$ coincide, leading to the
following theorem. This theorem was recently proved by the second
author.  The proof will appear in her PhD thesis.   

\begin{Th} \cite{Weaver-PhD-2022}\label{conj:s2=s3}
For any $1 \leq k \leq n$, the number of smooth decorated permutations
in $\SSn$ with exactly $k$ components in its SIF decomposition is
$\frac{b_{n,k}}{(n-k+1)!}$. 
\end{Th}

\begin{Cor} \label{conj:terms.integers}
The numbers $\frac{b_{n,k}}{(n-k+1)!}$ are integers for any
$1 \leq k \leq n$. 
\end{Cor}

\begin{figure}[h]
\begin{center}
\begin{tabular}{c|l|l|l|l|l|l|l|l|l|l|l}
\diagbox{$n$}{$k$} & 0 & 1 & 2 & 3 & 4 & 5 & 6 & 7 & 8 & 9 & 10\\
\hline
1 & 1 & 1 & & & & & & & & &\\
\hline
2 & 1 & 3 & 1 & & & & & & & &\\
\hline
3 & 1 & 7 & 7 & 1 & & & & & & & \\
\hline
4 & 1 & 15 & 29 & 15 & 1 & & & & & & \\
\hline
5 & 1 & 31 & 96 & 96 & 31 & 1 & & & & & \\
\hline
6 & 1 & 63 & 282 & 440 & 282 & 63 & 1 & & & & \\
\hline
7 & 1 & 127 & 771 & 1688 & 1688 & 771 & 127 & 1 & & & \\
\hline
8 & 1 & 255 & 2011 & 5803 & 8089 & 5803 & 2011 & 255 & 1 & & \\
\hline
9 & 1 & 511 & 5074 & 18520 & 33721 & 33721 & 18520 & 5074 & 511 & 1 & \\
\hline
10 & 1 & 1023 & 12488 & 55998 & 127698 & 166325 & 127698 & 55998 & 12488 & 1023 & 1
\end{tabular} 
\end{center}
\caption{Table of $s_1(n,k)$.}  \label{Fig:s1}
\end{figure}


\begin{figure}[h]
\begin{center}
\begin{tabular}{c|l|l|l|l|l|l|l|l|l|l}
\diagbox{$n$}{$k$} & 1 & 2 & 3 & 4 & 5 & 6 & 7 & 8 & 9 & 10\\
\hline
1 & 2 & & & & & & & & &\\
\hline
2 & 1 & 4 & & & & & & & &\\
\hline
3 & 2 & 6 & 8 & & & & & & & \\
\hline
4 & 3 & 18 & 24 & 16 & & & & & & \\
\hline
5 & 4 & 40 & 100 & 80 & 32 & & & & & \\
\hline
6 & 5 & 78 & 305 & 440 & 240 & 64 & & & & \\
\hline
7 & 6 & 140 & 798 & 1750 & 1680 & 672 & 128 & & & \\
\hline
8 & 7 & 236 & 1876 & 5838 & 8400 & 5824 & 1792 & 256 & & \\
\hline
9 & 8 & 378 & 4056 & 17136 & 34524 & 35616 & 18816 & 4608 & 512 & \\
\hline
10 & 9 & 580 & 8190 & 45480 & 122682 & 175896 & 137760 & 57600 & 11520 & 1024
\end{tabular}
\end{center}
\caption{Table of $s_2(n,k) = s_3(n,k)$.}  \label{Fig:s2}
\end{figure}




It is interesting to consider the asymptotic growth function for the
number $s(n)$ of smooth positroids on $[n]$.  We use the second
formula in \eqref{eq:smooth.en} to obtain the following data:


\begin{align*}
s(51)/s(50) &\approx  5.4489775,\\
s(101)/s(100)  &\approx  5.528236,\\
s(151)/s(150)  &\approx  5.555362, \\
s(201)/s(200)  &\approx 5.569062,\\
s(251)/s(250)  &\approx 5.5773263.  
\end{align*}
Based on this data, we conjecture the growth function is of the order $O(c^{n})$ for
some constant $c<6$.  

\section{Future Directions}\label{sec:future}

There are many parallels between Schubert varieties and positroid
varieties that have yet to be fully examined in the literature.  In
particular, the Johnson graph of a positroid $J(\M)$ is shown to play
a key role in the geometry of positroid varieties by \Cref{The: Main
theorem}, much like the Bruhat graph controls the geometry of Schubert
varieties \cite{carrell94}.

In 1998, Goresky, Kottwitz, and MacPherson \cite{GKM.1998} developed a
very general theory of equivariant cohomology for spaces with a torus
action and a moment map corresponding with a simple graph embedded in
the space with vertices given by fixed points of the torus action and
edges given by $1$-dimensional orbits connecting fixed points.  The
connection between these GKM graphs and equivariant cohomology has
been extensively developed over the past 24 years by Goldin
\cite{Goldin.99}, Guillemin-Holm-Zara \cite{Guillemin-Zara.2001,GHZ},
Harada-Henriques-Holm \cite{Harada-Henriques-Holm.2005}, and
Fukukawa-Ishida-Masuda \cite{Fukukawa-Ishida-Masuda.2014}.

For positroid varieties, we have a natural way to obtain an analog of
the Bruhat order and the Bruhat graph.  Recall the Johnson graph
$J(\M)$ has vertices indexed by bases in $I \in \M$, each of which
correspond bijectively with the $T$-fixed points in $\Pi_\M$.
If $I,J \in \M$ are two bases that are connected via an edge in
$J(\M)$, then by definition we can assume there is some exchange pair
$i<j$ such that $I\setminus \{i \} \cup \{j \} = J$.  Orient the edge
between them $I \longrightarrow J$. We call this the \textit{directed
Johnson graph} of $\M$.  This graph is closely related to the Gale
partial order on $k$-subsets in \Cref{sec:background}.

\begin{?}
Study the geometry and equivariant cohomology of a positroid variety
via the directed Johnson graph $J(\M)$.
\end{?}

\begin{?}
Identify the singular locus of a positroid variety.
\end{?}

\begin{?}
Study the enumeration and coset structure of the group operations of
flip, inverse, rotation for perms, or for derangements, or SIF perms.
\end{?}

\begin{?}
Study the combinatorics of the directed Johnson graphs of positroids.
\end{?}

\section{Acknowledgments}\label{Sec:Future}

We thank Herman Chau, Lauren Williams, Brendan Pawlowski, Stark Ryan, Joshua Swanson, and Allen Knutson for many insightful comments on this project.  Also, thanks to Yasutomo Kawashima for creating the Spirograph Maker app for the iphone (2017).


 \bibliographystyle{alpha} 
 \bibliography{positroid}

\end{document}